\documentclass[11pt,reqno]{amsart}

\usepackage[final]{epsfig}
\usepackage{amssymb,amsmath,amstext,amsthm,amsfonts}
\usepackage{euscript,a4wide}
\usepackage[ansinew]{inputenc}

\newtheorem{theorem}{Theorem}[section]
\newtheorem{proposition}[theorem]{Proposition}
\newtheorem{corollary}[theorem]{Corollary}
\newtheorem{lemma}[theorem]{Lemma}
\newtheorem{remark}[theorem]{Remark}
\newtheorem{maintheorem}{Theorem}
\newtheorem{maincorollary}[maintheorem]{Corollary}

\numberwithin{equation}{section}

\newcommand{\bR}{{\bf R}}
\newcommand{\bS}{{\bf S}}

\newcommand{\bN}{{\bf N}}

\newcommand{\cI}{{\mathcal I}}

\newcommand{\cP}{{\mathcal P}}

\newcommand{\cC}{{\mathcal C}}
\newcommand{\cD}{{\mathcal D}}

\newcommand{\cV}{{\mathcal V}}

\newcommand{\cT}{{\mathcal T}}

\newcommand{\de}{{\delta}}
\newcommand{\De}{{\Delta}}
\newcommand{\vfi}{{\varphi}}
\renewcommand{\epsilon}{\varepsilon}
\newcommand{\ep}{\varepsilon}

\DeclareMathOperator{\dist}{dist}

\theoremstyle{remark}

\title[Physical measures for infinite-modal maps] {
Physical measures for infinite-modal maps}

\begin{thanks} {V.A. was partially supported by CMUP-FCT
    (Portugal), CNPq (Brazil) and grants BPD/16082/2004 and
    POCI/MAT/61237/2004 (FCT-Portugal).  Part of this work
    was done while enjoying a post-doctorate leave from CMUP
    at PUC-Rio and IMPA. M.J.P. was partially supported by
    CNPq-Brazil/Faperj-Brazil/Pronex Dyn.  Systems.  }
\end{thanks}

\date{\today}

\author{V\'\i tor Ara\'ujo}
\address{Vítor Araújo, Instituto de Matemática,
  Universidade Federal do Rio de Janeiro, C. P. 68.530,
  21.945-970 Rio de Janeiro, RJ-Brazil \emph{and} Centro de
  Matemática da Universidade do Porto, Rua do Campo Alegre
  687, 4169-007 Porto, Portugal.}
\email{vitor.araujo@im.ufrj.br \text{and} vdaraujo@fc.up.pt}

\author{Maria Jos\'e Pacifico}

\address{Maria Jos\'e Pacifico, Instituto de Matematica, 
Universidade Federal do Rio de Janeiro, 
C. P. 68.530, CEP 21.945-970, 
Rio de Janeiro, R. J. , Brazil} 
\email{pacifico@impa.br and pacifico@im.ufrj.br}
\urladdr{http://www.dmm.im.ufrj.br}



\begin{document}

\begin{abstract}
We analyze certain parametrized families of one-dimensional
maps with infinitely many critical points from the
measure-theoretical point of view. We prove that such
families have absolutely continuous invariant probability
measures for a positive Lebesgue measure subset of
parameters. Moreover we show that both the densities of these
measures and their entropy vary continuously with the
parameter. In addition we obtain exponential rate of
mixing for these measures and also that they satisfy the
Central Limit Theorem.
\end{abstract}

\maketitle


\section{Introduction}
\label{sec:introduction}

One of the main goals of Dynamical Systems is to describe
the global asymptotic behavior of the iterates of most
points under a transformation of a compact manifold, either
from a topological or from a probabilistic (or ergodic)
point of view.  The notion of \emph{uniform hyperbolicity},
introduced by Smale in \cite{Sm}, and of \emph{non-uniform
  hyperbolicity}, introduced by Pesin \cite{Pe}, have been
the main tools to rigorously establish general results in
the field.

While uniform hyperbolicity is defined using only a finite
number of iterates of a given transformation, non-uniform
hyperbolicity is an asymptotic notion to begin with,
demanding the existence of non-zero Lyapunov exponents
almost everywhere with respect to some invariant probability
measure. 

On the one hand, the study of consequences of both notions in a
general setting has a long history, see
\cite{Mn,S,KH,B,BP,Y,BDV} for details and thorough
references.

On the other hand, it is rather hard in general to verify
non-uniform hyperbolicity, since we must take into account
the behavior of the iterates of the given map when time goes
to infinity. This was first achieved in the groundbreaking
work of Jakobson~\cite{Ja} on the quadratic family, which
was extended for more general one-dimensional families with
a unique critical point by many other mathematicians, see
e.g. \cite{BC1,R,MS,T,TTY}.  One-dimensional families with
two critical points were first considered in \cite{Ro} and
multimodal maps and maps with critical points and
singularities with unbounded derivative were treated in
\cite{LT,LV,BLS}. To the best of our knowledge, maps with
\emph{infinitely many critical points} were first dealt with in
\cite{PRV98}.

\begin{figure}[htb]
\centerline{\psfig{figure=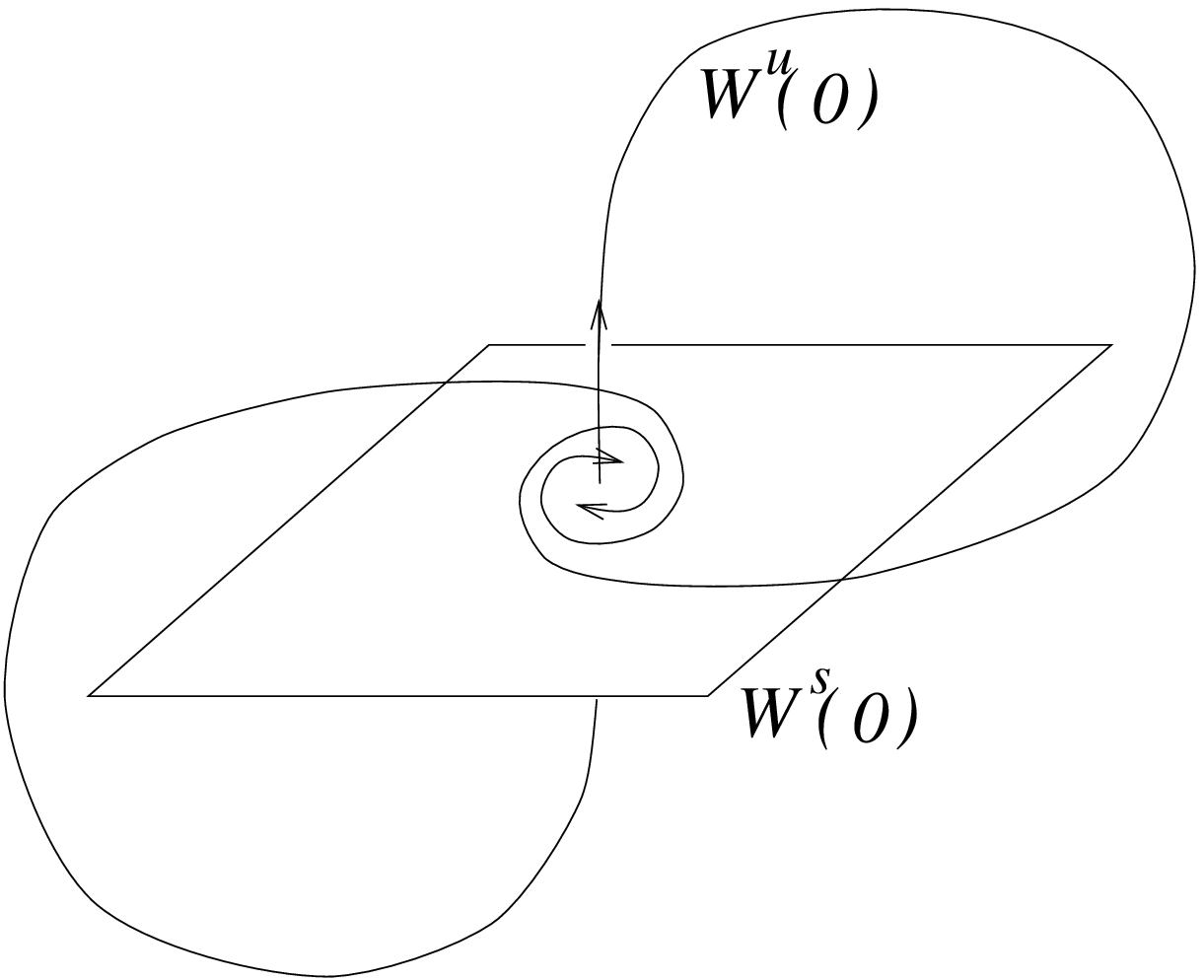,height=3cm}}
\caption{\label{Fig1} Double saddle-focus homoclinic connections}
\end{figure}

The aim of this paper is prove that the dynamics of the
family considered in \cite{PRV98}, for a positive Lebesgue
measure subset of parameters, is non-uniformly hyperbolic
and to deduce some consequences from the ergodic point of
view.  These families naturally appear as one-dimensional
models for the dynamical behavior near the unfolding of a
double saddle-focus homoclinic connection of a flow in a
three-dimensional manifold, see Figure~\ref{Fig1}
and~\cite{Sh}.  The main novelty is that we prove
\emph{global stochastic behavior} for a family of maps with
\emph{infinitely many regions of contraction}.

Roughly speaking, the family $f_\mu$ of one-dimensional
circle maps which we consider here is obtained from
first-return maps of the three-dimensional flow in
Figure~\ref{Fig1} to appropriate cross-sections and
disregarding one of the variables.  This reduction to a
one-dimensional model greatly simplifies the study of this
kind of unfolding and provides important insight to its
behavior. However as we shall see the dynamics of the
reduced model is still highly complex.

This family of maps is obtained translating the left-hand
side and right-hand side, vertically in opposite directions,
of the graph of the map $f=f_0$ described in
Figure~\ref{Fig4}.  This family approximates the behavior of
any generic unfolding of $f_0$.  Such unfolding was first
studied in~\cite{PRV98}, where it was shown that for a
positive Lebesgue measure subset $S$ of parameters the map
$f_\mu$, for $\mu\in S$, exhibits a chaotic attractor.  This
was achieved by proving that the orbits of the critical
values of $f_\mu$ have positive Lyapunov exponent and that
$f_\mu$ has a dense orbit.

Here we complement the topological description of the
dynamics of $f_\mu$ provided by~\cite{PRV98} for $\mu\in S$
with a probabilistic description constructing for the same
parameters a \emph{physical} probability measure $\nu_\mu$.
We say that an invariant probability measure $\nu$ is
\emph{physical} or \emph{Sinai-Ruelle-Bowen} (SRB) if there
is a positive Lebesgue measure set of points $x\in{\bf S}^1$
such that
\[
\lim_{n\to\infty}
\frac{1}{n}\sum_{k=0}^{n-1}\varphi\left(f_\mu^k(x)\right)
=\int\varphi\, d\nu,
\]
for any observable (continuous function) $\varphi:{\bf
  S}^1\to{\bf R}$. The set of points $x\in{\bf S}^1$ with
this property is called the \emph{basin} of $\nu$.  SRB
measures provide a statistical description of the asymptotic
behavior of a large subset of orbits. Combining this with
the results from~\cite{PRV98} we have that $f_\mu$ has
non-zero Lyapunov exponent almost everywhere with respect to
$\nu_\mu$, i.e. $f_\mu$ is non-uniformly hyperbolic for
$\mu\in S$.

The main feature needed for the construction of such
measures is to obtain positive Lyapunov exponent for
Lebesgue almost every point under the action of $f_\mu, \,
\mu\in S$. The presence of critical points is a serious
obstruction to achieve an asymptotic expansion rate on the
derivative of most points.  Therefore the control of
derivatives along orbits of the critical values is a central
subject in the ergodic theory of one-dimensional maps.

The crucial role of the orbits of the critical values on the
statistical description of the global dynamics of
one-dimensional maps was already present in the pioneer work
of Jakobson~\cite{Ja}, who considered quadratic maps and
obtained SRB measures for a positive Lebesgue measure subset
of parameters.

This was later followed by the celebrated papers of
Benedicks and Carleson~\cite{BC1,BC2}, where the parameter
exclusion technique was used to show that, for a positive
Lebesgue measure subset of parameters, the derivative along
the orbit of the unique critical value has exponential
growth and satisfies what is nowadays called a \emph{slow
  recurrence} condition to the critical point. This is
enough to construct SRB measures for those parameters.

Recently, in the unimodal setting it was established that
indeed the existence of SRB measures, and the exponential
growth of the derivative along the orbit of the critical
value, are equivalent conditions for Lebesgue almost every
parameter for which there are no sinks,
see~\cite{ALM,AM1,AM2}. See also \cite{BLS} for multimodal
maps.

In~\cite{PRV98} the technique of exclusion of parameters was
extended to deal with infinitely many critical orbits. Here
we refine this technique to obtain exponential growth of the
derivatives and slow recurrence to \emph{the whole critical
  set for Lebesgue almost every orbit}. By \cite{ABV00} this
ensures the existence of SRB measures for every parameter
$\mu\in S$, see Subsection~\ref{sec:exist-absol-cont} and
Theorem~\ref{thm:existacim}.

Moreover we are able to control the measure of the set of
points whose orbits are too close to the critical set during
the first $n$ iterates, showing that its Lebesgue measure is
exponential in $n$, see
Theorem~\ref{thm:approxestimate}. In addition, the Lebesgue
measure of the set of points whose derivative does not grow
exponentially fast in the first $n$ iterates decreases
exponentially fast with $n$, see
Theorem~\ref{thm:expansionestimate}.  By recent general
results on the ergodic theory of non-uniformly hyperbolic
systems \cite{ALP02,G}, both estimates above taken together
imply exponential decay of correlations for H\"older
continuous observables for $\nu_\mu$ and also that $\nu_\mu$
satisfies the Central Limit Theorem, for all $\mu\in S$, see
Subsection~\ref{sec:expon-decay} and
Corollary~\ref{cor:decayCLT}.
We remark that these properties are likewise satisfied by
uniformly expanding maps of $\textbf{S}^1$, which are the
touchstone of chaotic dynamics, see e.g. \cite{B,V}, in
spite of the presence of infinitely many points with
unbounded contraction (critical points).


Furthermore analyzing our arguments we observe that all the
estimates obtained do not depend on the choice of the
parameter $\mu\in S$. This shows after \cite{Al03,AKT} that
the density $d\nu_\mu/d\lambda$ of the SRB measure $\nu_\mu$
with respect to Lebesgue measure and its entropy
$h_{\nu_\mu}(f_\mu)$ vary continuously with $\mu\in S$, see
Subsection~\ref{sec:cont-vari-dens} and
Corollary~\ref{cor:contSRBEntropy}.  This type of result was
recently obtained in~\cite{F} for quadratic maps on the set
of parameters constructed in~\cite{BC1,BC2} using a similar
strategy.  Hence \emph{statistical properties of the maps
  $f_\mu$ for $\mu\in S$ are stable under small variations
  of the parameter}, i.e. this family is \emph{statistically
  stable} over $S$.

We emphasize that although the general strategy for proving
our results follows~\cite{BC1,BC2,PRV98,F} several new
difficulties had to be overcome. Indeed
unlike~\cite{BC1,BC2,PRV98} where the main purpose was to
obtain positive Lyapunov exponent along the orbits of
critical values, here we need to obtain \emph{positive
  Lyapunov exponents and slow recurrence to the critical set
  along almost every orbit}, which forces us to control the
distance to the critical set for far more iterates than
in~\cite{PRV98}. This demands at several places a bound on
the ratio between the second derivative at points nearby
the critical set. However there are inflection points which
impose extra restrictions on the arguments used in~\cite{PRV98}

Moreover with infinitely many critical points the derivative
of the smooth maps we consider here is not globally bounded,
(unlike any smooth unimodal family, see~\cite{BC1,BC2,F})
which demanded a proof of an exponential bound for the
derivative along the orbits of critical values.  In order to
obtain such a bound for a positive Lebesgue measure set of
parameters we changed the construction presented in
\cite{PRV98} adding a new constraint in the exclusion of
parameters algorithm.

The paper is organized as follows.  We first state precisely
our results in Subsections~\ref{sec:exist-absol-cont}
to~\ref{sec:cont-vari-dens}. We sketch the proof in
Section~\ref{sec:idea-construction}.  In
Section~\ref{sec:concstr-part-bound-distortion} we explain
how a sequence $(\cP_n)_{n\ge0}$ of partitions of ${\bf
  S}^1$ whose atoms have bounded distortion under action of
$f_\mu^n$ is constructed. Basic lemmas are stated and proved
in Section~\ref{sec:auxiliary-lemmas}. These are used to
obtain the main estimates in
Section~\ref{sec:fundamental-lemma}.  In
Sections~\ref{sec:slow-recurrence}
and~\ref{sec:fast-expansion-most} we use the main estimates
to deduce slow recurrence to the critical set and fast
expansion for most points.  In
Section~\ref{sec:extra-excl-param} we explain how an
exponential upper bound on the growth of the derivatives
along critical orbits can be obtained through an extra
condition imposed on the construction performed
in~\cite{PRV98} without loss.  Finally in
Section~\ref{sec:const-depend-unif} we keep track of the
estimates obtained during our constructs and show that they
do not depend on the parameter $\mu\in S$.

\subsection{Statement of the results}
\label{sec:statement-results}

\begin{figure}[htb]
\centerline{\psfig{figure=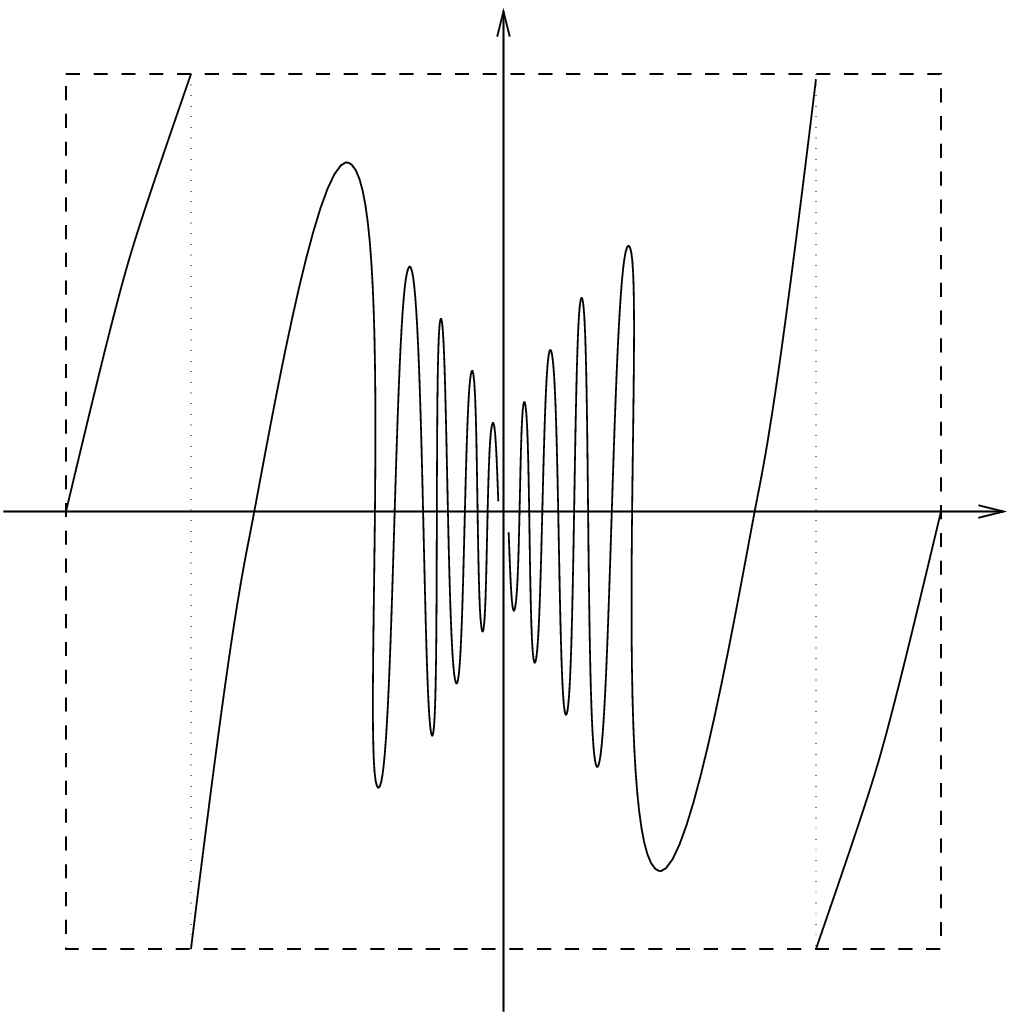,height=3cm}}
\caption{\label{Fig4} Graph of the circle map $f$.}
\end{figure}

Let $\hat f$ be the interval map
$\hat{f}:[-\epsilon_1,\epsilon_1]\to[-1,1]$ given by
\begin{equation}
\label{e2.3}
\hat{f}(z)=\left\{
\begin{array}{ll}
az^{\alpha}\sin(\beta\log(1/z)) & \mbox{   if   }z>0\\
-a|z|^{\alpha}\sin(\beta\log(1/|z|)) & \mbox{   if   }z<0,
\end{array}
\right.
\end{equation}
where $a>0$, $0<\alpha<1, \, \beta>0$ and $\ep_1>0$, see
Figure~\ref{Fig4}.

Maps $\hat{f}$ as above  have infinitely many critical points,
of the form
\begin{equation}
\label{e2.4}
x_{k}=\hat{x}\exp(-k\pi/\beta)
\mbox{    and    }x_{-k}=-x_{k}
\mbox{    for each large   }k>0
\end{equation}
where
$\hat{x}=\exp\big(-\frac1{\beta}\tan^{-1}\frac{\beta}{\alpha}\big)>0$
is independent of $k$.  Let $k_0\ge 1$ be the smallest
integer such that $x_k$ is defined for all $|k|\ge k_0$, and
$x_{k_0}$ is a local minimum.

We extend this expression to the whole circle $S^1=I/\{-1\sim
1\}$, where $I=[-1,1]$, in the following way.  Let $\tilde{f}$ be an
orientation-preserving expanding map of $S^1$ such that
$\tilde{f}(0)=0$ and $\tilde{f}'>\tilde\sigma$ for some constant
$\tilde\sigma>>1$. We define
$\epsilon=2\cdot x_{k_0}/(1+e^{-\pi/\beta})$,
so that $x_{k_0}$ is the middle point of the interval
$(e^{-\pi/\beta}\epsilon,\epsilon)$ and
fix two points $x_{k_0}<\hat{y}<\tilde{y}<\epsilon$, with
\begin{equation}\label{eq:condhaty}
|\hat{f}'(\hat{y})|>>1\quad\mbox{and also}\quad
2 \frac{1-\epsilon^\tau}{1+e^{-\pi/\beta}} x_{k_0}
> \hat y >  x_{k_0},
\end{equation}
where $\tau$ is a small positive constant to be defined in
what follows and we take $k_0=k_0(\tau)$ sufficiently big
(and $\epsilon$ small enough) in order that
\eqref{eq:condhaty} holds.  Then we take $f$ to be any
smooth map on $S^1$ coinciding with $\hat{f}$ on
$[-\hat{y},\hat{y}]$, with $\tilde{f}$ on
$S^1\setminus[-\tilde{y},\tilde{y}]$, and monotone on each
interval $\pm[\hat{y},\tilde{y}]$.

Finally let $f_\mu$ be the following one-parameter family of circle
maps unfolding the dynamics of $f=f_0$

\begin{equation}
\label{e2.4,8}
f_{\mu}(z)=\left\{
\begin{array}{ll}
f(z)+\mu & \mbox{   for   } z\in (0,\epsilon]\\
f(z)-\mu & \mbox{   for   } z\in [-\epsilon,0)
\end{array}
\right.
\end{equation}
for $\mu\in(-\epsilon,\epsilon)$. For
$z\in\mathbf{S}^1\setminus[-\epsilon,\epsilon]$ we assume
only that $\big|\frac{\partial}{\partial z}
f_\mu(z)\big|\ge2$.  In what follows we write
$z_k^\pm(\mu)=f_\mu(x_k)$ for $|k|\ge k_0$.

\begin{theorem}{\cite[Theorem A]{PRV98}}
\label{th:ref-annals}
For a given $\sigma\in(1,\sqrt{\tilde\sigma})$ there exists
an integer $N$ such that taking $k_0>N$ in the construction
of $(f_\mu)_{\mu}$, we can find a small positive constant $\tilde\rho$
such that for $0<\rho<\tilde\rho$ there exists a positive
Lebesgue measure subset
$S\subset[-\ep,-\ep^2]\cup[\ep^2,\ep]$ satisfying for every
$\mu\in S$
\begin{enumerate}
\item for all $n\ge 1$ and all $k_0\le |k|\le\infty$
  \begin{enumerate}
  \item $\left|\left(f_\mu^n\right)'(z_k^\pm(\mu))\right|
  \ge \sigma^n$;
  \item either $|f_{\mu}^{n}(f_{\mu}(x_l))|>\epsilon$
  or $|f_{\mu}^{n}(f_{\mu}(x_l))-x_{m(n)}|\geq
  e^{-\rho n}$;
  \end{enumerate}
 where $x_{m(n)}$ is the critical point
  nearest $f_{\mu}^{n}(f_{\mu}(x_l))$.
\item
  $\liminf_{n\to+\infty}n^{-1}\log|(f_\mu^n)'(z)|\ge\log\sigma/3$
  for Lebesgue almost every point $z\in S^1$;
\item there exists $z\in S^1$ whose orbit $\{f_\mu^n(z):n\ge 0\}$ is
  dense in $S^1$.
\end{enumerate}
\end{theorem}

The statement of Theorem~\ref{th:ref-annals} is slightly
different from the main statement of \cite{PRV98} but the
proof is contained therein.

\subsection{Existence of absolutely continuous invariant
  probability measures}
\label{sec:exist-absol-cont}

The purpose of this work is to prove that for parameters
$\mu\in S$ the map $f_\mu$ admits a unique absolutely continuous
invariant probability measure $\nu_\mu$, whose basin covers
Lebesgue almost every point of $\mathbf{S}^1$, and to study some
of the main statistical and ergodic properties of these measures.

In what follows we write $\lambda$ for the normalized
Lebesgue measure on ${\bf S}^1$.  Our first result shows the
existence of the $SRB$ measure.

\begin{maintheorem}
\label{thm:existacim}
Let $\mu\in S$ be given. Then there exists a
$f_\mu$-invariant probability measure $\nu_\mu$ which is
absolutely continuous with respect to $\lambda$ and such that for
$\lambda$-almost every $x\in {\bf S}^1$ and every continuous
$\vfi:{\bf S}^1\to\bR$
\begin{equation}
  \label{eq:1}
  \lim_{n\to+\infty}\frac1n\sum_{j=0}^{n-1}\vfi(f_\mu^j(x))
  = \int \vfi\, d\nu_\mu.
\end{equation}
\end{maintheorem}

The proof is based on the technique of parameter exclusion
developed in~\cite{PRV98} to prove
Theorem~\ref{th:ref-annals} and on recent results on
hyperbolic times for non-uniformly expanding maps with
singularities and criticalities, from~\cite{ABV00}.

In our setting \emph{non-uniform expansion} means the same as item
(2) of Theorem~\ref{th:ref-annals}. However due to the
presence of (infinitely many) criticalities and the
singularity at $0$, an extra condition is needed to
construct the $SRB$ measure: we need to control the average
distance to the critical set
along most orbits.

We say that $f_\mu$ has \emph{slow recurrence to the
  critical set} $\cC=\{x_k: |k| \ge k_0 \}\cup\{0\}$ if, for
every $\de>0$, there exists $\flat>0$ such that
\begin{equation}
  \label{eq:slow_recurrence_def}
  \limsup_{n\rightarrow\infty} \, 
\frac1n\sum_{k=0}^{n-1}-\log\dist_\flat\left(f_\mu^k(x),\cC\right)
< \de
\quad\mbox{for Lebesgue almost every}\quad x\in {\bf S}^1,
\end{equation}
where $\flat$ is a small positive value, and
$\dist_{\flat}(x,y)=  |x-y| $ if $ |x-y|\leq \flat$ and
$1$ otherwise.

Let $f:I\setminus\cC\to I$ be a $C^2$ map.  We say that
$\cC$ is a \emph{non-flat critical set}
if
there exist constants $B>1$ and $\beta>0$ such that
\begin{description}
\item[S1]\quad
$\displaystyle{
\frac{1}{B} \dist(x,\cC)^{\beta}
\le |f'(x)|
\le B \dist(x,\cC)^{-\beta}}$ \,;

\item[S2] \quad
$\displaystyle{
|\log|f'(x)| - \log |f'(y)|\,|
\le B \frac{ |x-y|}{\dist(x,\cC)^{\beta}}}$ \,;
\end{description}
for every  $x, y \in I\setminus\cC$
with $|x-y| < \dist(x,\cC)/2$.

The following result ensuring the existence of finitely many
physical probability measures is proved in~\cite{ABV00}.

\begin{theorem}
\label{t.singular}
If $f$ satisfies (S1), (S2), is non-uniformly expanding and
has slow recurrence to the critical set $\cC$, then there
are finitely many $\mu_1,\dots, \mu_l$ ergodic absolutely
continuous $f$-invariant probability measures such that
Lebesgue almost every point in $I$ belongs to the basin of
$\mu_i$ for some $i\in\{1,\dots,l\}$.
\end{theorem}

The maps $f_\mu$ satisfy conditions (S1)-(S2) above. 
Indeed we define $y_{k}=2\cdot x_k/(1+e^{-\pi/\beta})$, for
each $k\geq k_{0}$, so that $x_{k}$ is the middle point of
the interval $(y_{k+1},y_{k})$. 
We also use a similar notation for $k\le -k_0$.  We will
argue using the following lemmas, which correspond to Lemmas
3.2 and 3.3 proved in~\cite{PRV98}.

\begin{lemma}
\label{l3.1}
There exists $C>0$ depending on $\hat f$ only (not depending on
$\ep$ nor $\mu$) such that,
for every $x\in(y_{l+1},y_{l})$ and $l\geq k_0$, respectively,
$x\in(y_{l},y_{l-1})$ and $l\leq - k_0$, we have
\begin{enumerate}
\item $C^{-1}|x_{l}|^{\alpha-2} \cdot |x-x_l|^2
 \le |f(x)-f(x_{l})|
 \le C |x_{l}|^{\alpha-2} \cdot |x-x_l|^2$;
\item $
C^{-1}|x_{l}|^{\alpha-2}\cdot |x-x_{l}|
 \le |f'(x)|
 \le C |x_{l}|^{\alpha-2}\cdot |x-x_{l}|.
$
\end{enumerate}
\end{lemma}

\begin{lemma}
\label{l3.1,5}
Let $s,t\in[y_{l+1},y_{l}]$ with $l\geq k_{0}$,
respectively, $s,t\in[y_{l},y_{l-1}]$ with $l\leq -k_{0}$.
Then
\[
\left|\frac{f_\mu'(s)-f_\mu'(t)}{f_\mu'(t)}\right|
\leq 
K_1\frac{|s-t|}{|t-x_{l}|}
\]
where $K_1>0$ is independent of $l,s,t,\epsilon$ and $\mu$.
\end{lemma}

On the one hand since $0<\alpha<1$, 
$x\in(y_{l+1},y_{l})$ and $|x_l|<1$, then from item 2 of
Lemma~\ref{l3.1}
\[
C |x_{l}|^{\alpha-2}|x-x_{l}| = \big(C
|x_{l}|^{\alpha-2}|x-x_{l}|^2\big) |x-x_{l}|^{-1}
\le
\big(C |x_{l}|^{\alpha-2} |x_l|^2 \big) |x-x_{l}|^{-1}
\le C |x-x_{l}|^{-1}.
\]
On the other hand since $\alpha-2<0$ and $|x_l|<1$ we get $
C^{-1}|x_{l}|^{\alpha-2}|x-x_{l}|\ge C^{-1}|x-x_{l}|, $
showing that (S1) holds for $f_\mu$ with $B=C$ and
$\beta=1$, whenever $x\in(y_{k+1},y_k)$ and $x_k$ is the
closest critical point to $x$. Otherwise, if
$x\in(y_{k+1},y_k)$ and $x_{k+1}$ is the closest critical
point to $x$, then we have $|x-x_k|>|x-x_{k+1}|$ and so
by the above calculations we get
\begin{align*}
  |f_\mu'(x)|&\le
  C |x-x_k|^{-1} = C |x-x_{k+1}|^{-1}\cdot
  \left|\frac{x-x_k}{x-x_{k+1}} \right|^{-1}
  \le C |x-x_{k+1}|^{-1}, \quad\mbox{and}
  \\
  |f_\mu'(x)|&\ge
  \frac1C |x-x_k| = \frac1C |x-x_{k+1}| \cdot
  \left|\frac{x-x_k}{x-x_{k+1}} \right|
  \ge \frac1C |x-x_{k+1}|.
\end{align*}
This shows that (S1) is true for $f$ in all cases.

To check that (S2) also holds we write
\[
\frac{|f_\mu'(x)|}{|f_\mu'(y)|} = 
\frac{|f_\mu'(x)-f_\mu'(y)+f_\mu'(y)|}{|f_\mu'(y)|} \le
1+ \frac{|f_\mu'(x)-f_\mu'(y)|}{|f_\mu'(y)|} 
\]
and then because $\log(1+z)\le z$ for $z>-1$ we get
\[
|\log|f_\mu'(x)| - \log
 |f_\mu'(y)|\,|\le\frac{|f_\mu'(x)-f_\mu'(y)|}{|f_\mu'(y)|}
\le  K_1\frac{|x-y|}{|x-x_{l}|},
\]
which, by the same observation during the proof of (S1), is
enough to prove (S2) in all cases.

Thus according to Theorem~\ref{t.singular} and after
Theorem~\ref{th:ref-annals}, we only need to show that
$f_\mu$ has slow recurrence to the critical set for $\mu\in
S$ to achieve the result stated in
Theorem~\ref{thm:existacim}. This is done in
Sections~\ref{sec:auxiliary-lemmas}
to~\ref{sec:slow-recurrence}, where a stronger result is
obtained, as explained in what follows.


\subsection{Exponential decay of correlations and Central
  Limit Theorem}
\label{sec:expon-decay}

Using some recent developments on the statistical behavior
of non-uniformly expanding maps \cite{ALP02,G} we are able
to obtain exponential bounds on the decay of correlations
between H\"older continuous observables for $\nu_\mu$ with
$\mu\in S$. In addition it follows from standard techniques
that $\nu_\mu$ also satisfies the Central Limit Theorem. In
order to achieve this we refined the arguments
in~\cite{PRV98} using strong conditions on the exclusion of
parameters extending the estimates obtained therein for
critical orbits to get a exponential upper bound on the
growth of the derivative along orbits of critical values, as
explained in Section~\ref{sec:extra-excl-param}. Moreover we
where able to extend most of the estimates from~\cite{PRV98}
for Lebesgue almost every orbit, yielding an exponential
bound on the Lebesgue measure of the set of points whose
average distance to the critical set during the first $n$
iterates is small, as follows.

We first define the average distance to the critical set
\begin{equation}\label{eq:averagedist}
\cC_n^\flat(x)=
\frac1n\sum_{j=0}^{n-1}-\log\dist_\flat
\left(f_\mu^j(x),\cC\right).
\end{equation}
for a given $\flat>0$.
Then we are able to prove the following.

\begin{maintheorem}
  \label{thm:approxestimate}
  Let $\mu\in S$ and $\de>0$ be given. Then there are
  constants $C_1,\xi_1,\flat>0$ dependent on $\hat f$,
  $\sigma$, $k_0$ and $\delta$ only such that 
$\mathcal R(x)=\min\{ N\ge1 : \cC_n^\flat(x)<\de, \forall n\ge N\}$
satisfies
\[
\lambda\Big(
\{x\in {\bf S}^1: \mathcal R(x)>n \} \Big)
\le C_1\cdot e^{-\xi_1\cdot n}.
\]
\end{maintheorem}

We note that in particular this shows that $f_\mu$ has slow
recurrence to the critical set
and ensures the existence of the $SRB$ measure $\nu_\mu$ for
$\mu\in S$ by Theorem~\ref{t.singular}. 

We are also able to obtain, using the same techniques, an
exponential bound on the set of points whose expansion rate
up to time $n$ is less than the one prescribed by item (2)
of Theorem~\ref{th:ref-annals}. This is detailed in
Section~\ref{sec:fast-expansion-most}.

\begin{maintheorem}
  \label{thm:expansionestimate}
  Let $\mu\in S$ be given. Then there exist constants
  $C_2,\xi_2>0$ dependent on $\hat f$, $\rho$ and
  $k_0$ only such that $\mathcal E(x)=\min\{ N\ge1:
  \big|(f_\mu^n)'(x)\big|>\sigma^{n/3}, \forall n\ge N\}$
  satisfies
\[
\lambda\Big(\{x\in \mathbf{S}^1: 
\mathcal E(x)> n \}
\Big)\le C_2 \cdot e^{-\xi_2\cdot n}.
\]
\end{maintheorem}

In particular we obtain a new proof of item (2) of
Theorem~\ref{th:ref-annals}, which \emph{does not follow
directly from Theorem~\ref{thm:existacim} plus the Ergodic
Theorem} since it is not obvious whether $\log|f'|$ is
$\nu_\mu$ integrable.

Theorems~\ref{thm:approxestimate}
and~\ref{thm:expansionestimate} together ensure that for
$\mu\in S$ there are constants $C_3>0$ and $\xi_3\in(0,1)$
such that $\Gamma_n= \{x\in \mathbf{S}^1: \mathcal E(x)> n
\quad\mbox{or}\quad \mathcal R(x)>n \}$ satisfies
\begin{equation}
  \label{eq:strechedexp}
\lambda(\Gamma_n)\le C_3\cdot e^{-\xi_3\cdot n}
\end{equation}
for all $n\ge1$. This fits nicely into the following statements.

\begin{theorem}
  \label{thm:alp-g}
  Let $g:{\bf S}^1\to {\bf S}^1$ be a transitive $C^2$ local
  diffeomorphism outside a non-flat critical set $\cC$ such
  that~\eqref{eq:strechedexp} holds. Then
  \begin{enumerate}
  \item \cite[Theorem 1]{ALP02} there exists an absolutely
    continuous invariant probability measure $\nu$ and some
    finite power of $g$ is mixing with respect to $\nu$;
  \item \cite[Theorem 1.1]{G} there exist constants $C,c>0$
    such that the correlation function $
    \mathrm{Corr}_{n}(\varphi, \psi) = \left|\int (\varphi
      \circ g^{n})\cdot \psi \, d\nu - \int \varphi \, d\nu
      \int \psi \, d\nu\right|, $ for H\"older continuous
    observables $\varphi,\psi:{\bf S}^1\to\mathbf{R}$,
    satisfies for all $n\ge1$
\[
\mathrm{Corr}_{n}(\varphi, \psi)\le C \cdot e^{-c\cdot n}.
\]
\item \cite[Theorem 4]{ALP02} $\nu$ satisfies the Central
  Limit Theorem: given a H\"older continuous function \(
  \phi:{\bf S}^1\to\mathbf{R} \) which is not a coboundary
  (\( \phi\neq \psi\circ g - \psi \) for any \( \psi:{\bf
    S}^1\to\mathbf{R} \)) there exists \( \theta>0 \) such
  that for every interval \( J\subset \mathbf{R} \)
\[
\lim_{n\to\infty}
\nu \Big(\Big\{x\in \mathbf{S}^1: 
\frac{1}{\sqrt n}\sum_{j=0}^{n-1}\Big(\phi(g^{j}(x))-\int\phi d\nu
\Big)\in J \Big\} \Big)=
\frac{1}{\theta \sqrt{2\pi} }\int_{J} e^{-t^{2}/ 2\theta^{2}}dt.
\]
  \end{enumerate}
\end{theorem}

It is then straightforward to deduce the following conclusion.

\begin{maincorollary}
  \label{cor:decayCLT}
For every $\mu\in S$ the map $f_\mu$ has exponential
decay of correlations for H\"older continuous observables
and satisfies the Central Limit Theorem with respect to the
$SRB$ measure $\nu_\mu$.
\end{maincorollary}


\subsection{Continuous variation of densities and of
  entropy}
\label{sec:cont-vari-dens}

We note that during the arguments in Sections
\ref{sec:idea-construction} to~\ref{sec:fast-expansion-most}
the constants used in every estimation depend uniformly on
the values of $\rho,\sigma$ and $\ep$ which can be set right
from the start of the construction that proves
Theorems~\ref{thm:approxestimate}
and~\ref{thm:expansionestimate}.  This enables us to use
recent results of \emph{statistical stability} and
\emph{continuity of the $SRB$ entropy} from \cite{Al03,AKT},
showing that \emph{both the densities of the $SRB$ measures
  $\nu_\mu$ and the entropy vary continuously with $\mu\in
  S$}.

Let $\mathcal F$ be a family of $C^2$ maps of $\mathbf{S}^1$
outside a fixed non-flat critical set $\cC$ such that for
any given $f\in\mathcal F$ and $\delta_1>0$ there exists $\de_2>0$
satisfying for every measurable subset $E\subset
\mathbf{S}^1$
\[
\lambda(E)<\de_2\implies \lambda(f^{-1}(E))<\delta_1,
\]
that is $f_*(\lambda)\ll\lambda$. We say that a family
$\mathcal F$ as above is a \emph{non-degenerate family of maps}.

\begin{theorem}
  \label{thm:contvar}
  Let a non-degenerate family $\mathcal F$ of $C^2$ maps of
  $\mathbf{S}^1$ outside a fixed non-flat critical set $\cC$
  be given such that for every $f\in\mathcal F$ the
  corresponding functions $\mathcal E,\mathcal R:
  \mathbf{S}^1\to\mathbf{N}$ define a family
  $(\Gamma_n)_{n\ge1}$ satisfying~\eqref{eq:strechedexp}
  with constants $C_3,\xi_3$ not depending on $f\in\mathcal
  F$. Then
  \begin{enumerate}
  \item \cite[Theorem A]{Al03} the map $(\mathcal
    F,d_{C^2})\to (L^1(\lambda),\|\cdot\|_1),
    f\mapsto \frac{d\nu_f}{d\lambda}\in L^1(\lambda)$ is
    continuous, where $d_{C^2}$ is the $C^2$ distance and
    $\|\cdot\|_1$ the $L^1$-norm;
  \item \cite[Corollary C]{AKT} the map $(\mathcal
    F,d_{C_2})\to\mathbf{R}, f\mapsto h_{\nu_f}(f)$ is
    continuous.
  \end{enumerate}
\end{theorem}

We observe that $\mathcal F=\{f_\mu: \mu\in S\}$ satisfies
all the above conditions since 
\begin{itemize}
\item $\hat f$ is a $C^\infty$ map
whose non-zero singularities, albeit infinitely many, are of
quadratic type, and near zero $\hat f$ is bounded by
$|z|^\alpha$;
\item $f_\mu$ is obtained from $\hat f$ through a local
  diffeomorphism extension plus two translations (or rigid
  rotations when viewed on $\mathbf{S}^1$);
\item the values of $\beta,\ep,\sigma,\rho$ can be
  chosen so that
  \begin{itemize}
  \item $S$ is given by Theorem~\ref{th:ref-annals} with
  positive Lebesgue measure;
\item $f_\mu$ for $\mu\in S$ satisfies
  \eqref{eq:strechedexp} with $C_3,\xi_3>0$ depending only
  on $\ep,\sigma,\rho$ --- this is detailed in
  Section~\ref{sec:const-depend-unif}.
  \end{itemize}
\end{itemize}
Thus we deduce the following corollary which shows that
statistical properties of $f_\mu$ are stable
under small variations of the parameter $\mu$ within the set $S$.

\begin{maincorollary}
\label{cor:contSRBEntropy}
The following maps are both continuous:
\[
\begin{array}[l]{lll}
S &\to &(L^1(\lambda),\|\cdot\|_1)
\\
\mu & \mapsto & \frac{d\nu_\mu}{d\lambda}
\end{array}
\quad\mbox{and}\quad
\begin{array}[l]{lll}
S &\to & \mathbf{R}
\\
\mu & \mapsto & h_{\nu_\mu}(f)
\end{array}.
\]
\end{maincorollary}

\noindent
\textbf{Acknowledgments:} We are thankful to M. J. Costa and
M. Viana for helpful conversations about this work. We also
thank the referee for very useful comments and the
encouragement to correct a previous version of the proof.


\section{Idea of the proof}
\label{sec:idea-construction}

From now on we fix a parameter $\mu\in S$ and write
$\cC_\infty=\cup_{n=0}^{\infty}(f^n)^{-1}(\cC)$ for the set
of pre-orbits of the critical set $\cC$.  We also write
$f=f_\mu$ in what follows.  

Following \cite{PRV98} we consider a convenient partition
$\{I(l,s,j)\}$ of the phase space into subintervals, with a
bounded distortion property: trajectories with the same
itinerary with respect to this partition have derivatives
which are comparable, up to a multiplicative constant.  This
is done as follows.  Let $l\geq k_0$ and $y_{l}\in
(x_{l},x_{l-1})$ be as defined in
Subsections~\ref{sec:statement-results}
and~\ref{sec:exist-absol-cont}: $x_l$ is the middle point of
$(y_{l+1},y_{l})$.  We partition $(x_l,y_l)$ into
subintervals $ I(l,s)=(x_l+e^{-(\pi/\beta) s}\cdot(y_l-x_l),
x_l+e^{-(\pi/\beta)(s-1)}\cdot(y_l-x_l))$, $s\ge 1$.
\begin{figure}[htbp]
  \centerline{\psfig{figure=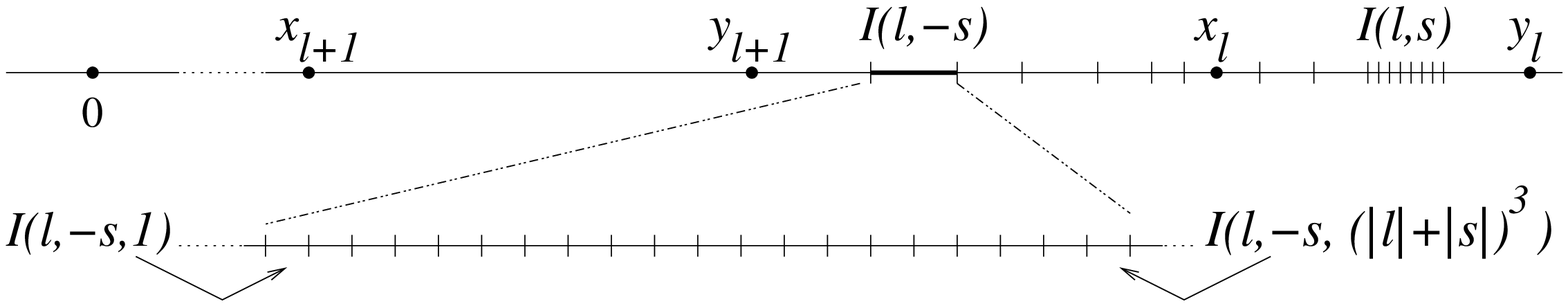,width=14cm}}
  \caption{The initial partition $\cP_0$.}
  \label{fig:partition}
\end{figure}
We denote by $I(l,-s)$ the subinterval of $(y_{l+1},x_l)$
symmetrical to $I(l,s)$ with respect to $x_l$.  We subdivide
$I(l,\pm s)$ into $(l+|s|)^3$ intervals $I(l,\pm s,j)$,
$1\le j\le (l+|s|)^3$ with equal length and $j$ increasing
as $I(l,\pm s, j)$ is closer to $x_l$, see
Figure~\ref{fig:partition}.  We also perform entirely
symmetric constructions for $l\leq -k_0$.  Let $I(\pm
k_0,1,1)$ be the intervals having $\pm\epsilon$ in their
boundaries.  Clearly we may suppose that $I(\pm k_0,1,1)$
are contained in the region
$\bS^1\setminus[-\tilde{y},\tilde{y}]$ where $f$ coincides
with $\tilde{f}$, and so $|f'|>\sigma_0>1$.  Finally, for
completeness, we set
$I(0,0,0)=I(0,0)=\bS^1\setminus[-\epsilon,\epsilon]$.

\begin{remark}
  \label{rmk:partitiondimension}
By the definition of $I(l,s,j)$ 
\[
|I(l,s,j)|=a_{1}\frac{e^{-(\pi/\beta)(|l|+|s|)}}{(|l|+|s|)^3}
\quad\text{and}\quad
a_{2}e^{-(\pi/\beta)(|l|+|s|)} \le
\dist(I(l,s,j),x_{l})\le a_{2}e^{-(\pi/\beta)(|l|+|s|-1)}
\]
where $|I|$ denotes the length of the interval $I$,
$a_{1}=\hat{x}\frac{(e^{(\pi/\beta)}-1)^2}{e^{(\pi/\beta)}+1}$
and
$a_{2}=\hat{x}\frac{e^{(\pi/\beta)}-1}{e^{(\pi/\beta)}+1}<1.$
Moreover for any $m\ge1$ we have $|x_m-x_{m+1}|=\hat x \cdot
(1-e^{-\pi/\beta})\cdot e^{-\frac{\pi}{\beta} m }$.  In
addition we have $\dist\big(I(l,s,j),0\big)=|x_l|\pm
\dist\big(I(l,s,j),x_{l}\big)$ according to the sign of $s$
and consequently $|\hat x-a_2|\cdot e^{-(\pi/\beta)|l|} \le
\dist(I(l,s,j),0)\le (a_2+\hat x)\cdot e^{-(\pi/\beta)|l|}$.
\end{remark}

We will separate the orbit of a point $x_0\in
I\setminus\cC_\infty$ into sequences of consecutive iterates
according to whether the point is near $\cC$ or is in the
expanding region $I(0,0,0)$. When $x_n=f^n(x_0)$ is near
$\cC$, we say that $n$ is a \emph{return time} and the
expansion may be lost.  But since we know that for $\mu\in
S$ the derivatives along the critical orbits grow
exponentially fast, we shadow the orbit of $x_n$ during a
\emph{binding period} by the orbit of the nearest critical
point and borrow its expansion. At the end of this binding
period, the expansion is completely recovered, which will be
explained precisely in Section~\ref{sec:auxiliary-lemmas}.

This picture is complicated by the infinite number of
critical points and by the possible returns near another
critical point during a binding period.
Iterates outside binding periods and return times
are \emph{free iterates}, where the derivative is uniformly
expanded.

Our main objective is to obtain slow recurrence to $\cC$,
which means that the returns of generic orbits are not too
close to $\cC$ on the average.  However \emph{even at a free
  iterate the orbits may be very close to the critical set},
by the geometry of the graph of $f_0$, which demands a
deeper analysis to achieve slow recurrence to the critical
set.
Moreover since $|f'|$ is not
bounded from above in our setting, we do not automatically
have an exponential bound on the derivative along orbits of
critical values, which is needed to better control the
recurrence to $\cC$ and must be proved by a separate argument
involving a stronger exclusion of parameters than in the
algorithm presented in~\cite{PRV98}.

Using the slow recurrence we show that the derivative along
the orbit of Lebesgue almost every point grows exponentially
fast.  Using the estimates from
Sections~\ref{sec:concstr-part-bound-distortion}
to~\ref{sec:fundamental-lemma} we are able to obtain more:
we deduce the exponential estimates on
Theorems~\ref{thm:approxestimate}
and~\ref{thm:expansionestimate} in
Sections~\ref{sec:slow-recurrence}
and~\ref{sec:fast-expansion-most}.

Finally the refinement on the parameter exclusion in
\cite{PRV98} and the dependence of the constants on the
choices made during the entire construction are detailed in
Sections~\ref{sec:extra-excl-param}
and~\ref{sec:const-depend-unif} respectively, where we
show that the estimates are uniform on $\mu\in S$.

\section{Refining the partition}
\label{sec:concstr-part-bound-distortion}

We are going to build inductively a sequence of partitions
$\mathcal{P}_0,\mathcal{P}_1,\ldots$ of $I$ (modulus a zero
Lebesgue measure set) into intervals.  We will define
inductively the sets
$R_n(\omega)=\left\{r_1,\ldots,r_{\gamma(n)}\right\}$ which
is the set of the return times of $\omega\in \cP_n$ up to
$n$ and a set
$Q_n(\omega)=
\left\{(l_1,s_1,j_1),\ldots,
(l_{\gamma(n)},s_{\gamma(n)},j_{\gamma(n)})\right\}$,
which records the indexes of the intervals such that
$f^{r_i}(\omega)\subset I(l_i,s_i,j_i)$,
$i=1,\ldots,r_{\gamma(n)}$.

In the process we will show inductively that for all $n\in \bN_0$
\begin{equation}
  \label{eq:fn-diffeo-omega}
  \forall \omega \in \mathcal{P}_n\quad 
f^{n+1}|_{\omega} \mbox{ is a diffeomorphism,}
\end{equation}
which is essential for the construction itself.
For $n=0$ we define
\[
\cP_0= \left\{I(0,0,0)\right\}\cup \left\{I(l,s,j):|l|\geq
  k_0, \, |s|\ge1 , \, 1\leq j \leq (|l|+|s|)^3\right\}.
\]
It is obvious that $\cP_0$ satisfies
\eqref{eq:fn-diffeo-omega} for $n=0$.  We set
$R_0\left(I(0,0,0)\right)=\emptyset$ and
$R_0(I(l,s,j))=\{0\}, Q_0(I(l,s,j))=\{(l,s,j)\}$ for all
possible $(l,s,j)\neq(0,0,0)$.

\begin{remark}
  \label{rmk:everybodyreturns}
This means that \emph{every} $I(l,s,j)$ with $|l|\ge k_0$,
$|s|\ge1$ and $j=1,\dots,(|l|+|s|)^3$ has a return at time
$0$, by definition. This will be important in
Section~\ref{sec:slow-recurrence}.
\end{remark}

For each $(l,s)$ with $|l|\ge k_0$ and $|s|\ge1$  such that
\begin{equation}
\label{e.afastado}
e^{-(\pi/\beta)|s|} \cdot
\frac{1-e^{-(\pi/\beta)}}{1+e^{-(\pi/\beta)}}
<  {\tau}, \quad\mbox{i.e.}\quad
|s|> s(\tau)=-\frac{\beta}{\pi}\log\left(
\tau\cdot \frac{1+e^{-(\pi/\beta)}}{1-e^{-(\pi/\beta)}}
\right),
\end{equation}
we define the \emph{binding period} $p(x)$ of $x\in I(l,s)$
to be the largest integer $p> 0$ such that
\begin{eqnarray}
|f^{h}(x_l)|\le\epsilon
&\text{and}&
|f^{h}(x)-f^{h}(x_{l})|
\leq |f^{h}(x_{l})-x_{m(h-1)}| e^{-\tau h}\nonumber\\
\label{BCeq}
& \text{or}\\
|f^{h}(x_{l})|>\epsilon
&\text{and}&
|f^{h}(x)-f^{h}(x_{l})|
\leq \epsilon^{1+\tau} e^{-\tau h}\nonumber
\end{eqnarray}
for all $1\leq h\leq p$, where $x_{m(h)}$ is the critical
point nearest $f^{h}(f(x_l))$ and $\tau>0$ is
a small constant to be specified during the construction.


Failing condition~\eqref{e.afastado} means that $I(l,s)$ is
not close enough to $\cC$ since
$$
|x-x_l|
 \ge e^{-(\pi/\beta)|s|} \cdot (y_l-x_l)
 \ge e^{-(\pi/\beta)|s|}\cdot 
\frac{1-e^{-(\pi/\beta)}}{1+e^{-(\pi/\beta)}} \cdot |x_l|
 \ge \tau|x_l|
 $$
 for all $x\in I(l,s)$, and in this case there is no
 expansion loss at the point $x$.  Indeed by Lemma~\ref{l3.1} and
 using the definition of $x_l$ from \eqref{e2.4} we get
 \begin{equation}
   \label{eq:expansaolivre}
   |f'(x)|
   \ge
   C^{-1}\cdot|x_{l}|^{\alpha-2}\cdot|x - x_{l}|
   \ge C^{-1}\cdot|x_{l}|^{\alpha-2}\cdot \tau |x_l|
   = \frac{\tau \hat x^{\alpha-1}}{C} 
   \cdot e^{(1-\alpha)\frac{\pi}{\beta}|l|}
 \end{equation}
 Since $1-\alpha>0$ and $|l|\ge k_0$, this ensures that
 $|f'(x)|>1$ if we take $k_0=k_0(\tau)$ big enough.

 \begin{remark}
   \label{rmk:k_0BIG}
   As we will explain along the proof, the values of $k_0$
   and $\tau^{-1}$ will both need to be taken sufficiently
   big.  We note that 
   $k_0\to\infty$ when $\tau\to0^+$. For more on these
   dependencies see Section~\ref{sec:const-depend-unif}.
 \end{remark}

We define \emph{the binding period $p(l,s)$ of the interval
$I(l,s)$ to be the smallest binding period of all points of
this interval}, that is $p(l,s)=\inf \{p(x) : x\in I(l,s) \}.$

\medskip

For $(l,s,j)$ with $|l|\ge k_0$, $|s|\ge 1$ and $1\leq j
\leq (|l|+|s|)^3$, write $I(l,s,j)^+$
for \emph{the union of $I(l,s,j)$ with its two adjacent intervals
in $\cP_0$.}

\medskip

Now we assume that $\cP_{n-1}$ is defined, satisfies
\eqref{eq:fn-diffeo-omega} and $R_{n-1},\,Q_{n-1}$ are also
defined on each element of $\mathcal{P}_{n-1}$. Fixing an
interval $\omega\in\mathcal{P}_{n-1}$ there are three
possible situations.

\begin{enumerate}
  
\item If $R_{n-1}(\omega)\neq\emptyset$ and
  $n<r_{\gamma(n-1)}+p(l_{\gamma(n-1)},s_{\gamma(n-1)})$
  then we say that $n$ is a \emph{bound time} for $\omega$,
  put $\omega\in\mathcal{P}_n$ and set
  $R_n(\omega)=R_{n-1}(\omega)$,
  $Q_{n}(\omega)=Q_{n-1}(\omega)$.
  
\item If either $R_{n-1}(\omega)=\emptyset$, or $n\geq
  r_{\gamma(n-1)}+p(l_{\gamma(n-1)},s_{\gamma(n-1)})$ and
  $f^n(\omega)\subset I(0,0,0)\cup I(\pm k_0,1,1)$ and
  $f^n(\omega)$ does not contain any $I(\pm k_0,1,1)$,
  then we say that $n$ is a \emph{free time} for $\omega$,
  put $\omega\in\mathcal{P}_n$ and set
  $R_n(\omega)=R_{n-1}(\omega)$,
  $Q_{n}(\omega)=Q_{n-1}(\omega)$.

\item If the two conditions above fail then
we consider two cases:

\begin{enumerate}
  
\item $f^n(\omega)$ does not cover completely any
  $I(l,s,j)$, with $|l|\geq k_0, |s|\ge 1$ and
  $j=1,\ldots,(|l|+|s|)^3$.  Because $f^n$ is continuous
  and $\omega$ is an interval, $f^n(\omega)$ is also an
  interval and thus it is contained in some $I(l,s,j)^+$, for a
  certain $|l|\geq k_0, |s|\ge1$ and
  $l=1,\ldots,(|l|+|s|)^3$, which is called the \emph{host
    interval}. 

\smallskip  

  If $|s|>s(\tau)$, then this $n$ is an \emph{inessential
    return time} for $\omega$ and we set
  $R_n(\omega)=R_{n-1}(\omega)\cup\{n\}$,
  $Q_{n}(\omega)=Q_{n-1}(\omega)\cup\{(l,s,j)\}$ and put
  $\omega\in\mathcal{P}_n$.

\smallskip  

  Otherwise, $|s|\le s(\tau)$ and there is no expansion
  loss, thus $n$ is again  a \emph{free time}, 
put $\omega\in\mathcal{P}_n$ and set
  $R_n(\omega)=R_{n-1}(\omega)$,
  $Q_{n}(\omega)=Q_{n-1}(\omega)$.

\item $f^n(\omega)$ contains at least an interval
  $I(l,s,j)$, with $|l|\geq k_0, |s|\ge1$ and
  $j=1,\ldots,(|l|+|s|)^3$, in which case we partition
  $\omega$ as follows.  We consider the sets
  \begin{align*}
\quad\qquad \omega_{\,l,s,j}'
  &= 
  f^{-n}\big(I(l,s,j)\big)\cap\omega 
\mbox{  for  } \left\{
  \begin{array}{l}
  |l|\ge k_0,\, |s|\ge1, 1\le j\le (|l|+|s|)^3 
\\
\mbox{and  } (|l|,s,j)\neq(k_0,1,1) 
\end{array}
\right.;
\\
\quad\qquad\omega_{\,0,0,0}'
  &=
  f^{-n}\big(I(0,0,0)\cup I(\pm k_0,1,1)\big)\cap\omega.
  \end{align*}

\medskip
Denoting by $\cI$ the set of indexes $(l,s,j)$ such that
$\omega_{\,l,s,j}'\neq\emptyset$ we have
\begin{equation}
  \label{eq:omega-decomp}
  \omega \setminus f^{-n}(\cC)=
  \bigcup_{(l,s,j)\in\cI} \omega_{\,l,s,j}'.
\end{equation}
By the induction hypothesis $f^n|_\omega$ is a
diffeomorphism and then each $\omega_{\,l,s,j}'$ is an
interval.  Moreover $f^n(\omega_{\,l,s,j}')$ covers the
whole $I(l,s,j)$ for $|l|\ge k_0,\, |s|\ge1,
j=1,\dots,(|l|+|s|)^3$, except eventually for one or two end
intervals. When $f^n(\omega_{\,l,s,j}')$ does not cover
entirely $I(l,s,j)$ we enlarge $\omega_{\,l,s,j}'$ gluing it
with its adjacent intervals in \eqref{eq:omega-decomp},
getting a new decomposition of $\omega\setminus
f^{-n}(\cC)$ into intervals
$\omega_{\,l,s,j}$ such that
\[
 \qquad\qquad
I(l,s,j)\subset f^n(\omega_{\,l,s,j})\subset I(l,s,j)^+
\mbox{  for  } |l|\geq k_0,\, |s|\ge1,\, 
1\le j\le (|l|+|s|)^3.
\]
We put $\omega_{\,l,s,j}\in
\cP_n$ for all $(l,s,j)$ such that
$\omega_{\,l,s,j}\neq\emptyset$, with $|l|\geq k_0$. This
results in a refinement of $\cP_{n-1}$ at $\omega$.

We set $Q_n(\omega_{\,l,s,j})=
Q_{n-1}(\omega)\cup\{(l,s,j)\}$ for every non-empty interval
$\omega_{\,l,s,j}$.  The interval $I(l,s,j)^+$ is the
\emph{host interval} of $\omega_{\,l,s,j}$.  We define
$R_n(\omega_{\,l,s,j})\cup\{n\}$ and we say that $n$ is an
\begin{enumerate}
\item \emph{escape time} for $\omega_{\,l,s,j}$ if $|l|\ge
  k_0$ and $s\le s(\tau)$.
\item \emph{essential return time} for $\omega_{\,l,s,j}$ if
  $|l|\ge k_0$ and $s> s(\tau)$.
\end{enumerate}

\end{enumerate}

\end{enumerate}

\begin{remark}
  \label{rmk:expansaolivre}
  We note that if $n$ is a free time or an escape time for
  $z$, then $x=f^n(z)$ either is in the region
  $\bS^1\setminus[-\ep,\ep]$ and thus $|f'(x)|\gg 1$, or
  satisfies the inequality~\eqref{eq:expansaolivre}. Hence
  on free times and escape times we always have expansion of
  derivatives bounded from below by some uniform constant
  $\sigma_0>1$.  We stress that we may and will assume that
  $\sigma_0>\max\{e,\sqrt{\tilde\sigma}\}$ in what follows.
\end{remark}

To complete the induction step all we need is to check that
\eqref{eq:fn-diffeo-omega} holds for $\cP_n$. Since
for any interval $J\subset \bS^1$
\[
\left.
\begin{array}{l}
  f^n|_J \mbox{ is a diffeomorphism }\\
  \cC\cap f^n(J)=\emptyset
\end{array}
\right\}\Rightarrow  f^{n+1}|_J \mbox{ is a diffeomorphism},
\]
all we are left to prove is that $\cC\cap
f^n(\omega)=\emptyset$ for all $\omega\in\cP_n$.

Let $\omega\in\mathcal{P}_n$. If $n$ is a free time for
$\omega$ then we are done. If $n$ is either a return time
for $\omega$, essential or inessential, or an escape time,
then by construction we have that $f^n(\omega)\subset
I(l,s,j)^+$ for some $|l|\geq k_0$, $|s|\ge 1$,
$j=1,\ldots,(|l|+|s|)^3$ (or for $(l,s,j)=(0,0,0)$) and thus
$\cC\cap f^n(\omega)=\emptyset$.  For the binding case
we use the following estimate.

\begin{proposition}
  \label{pr:bindingaway}
  Let $n\ge1$ and $\omega\in\cP_n$ be such that $n$ is a
  binding time for $\omega$. Then either $|f^n(x)|>
  x_{k_0}$ or $\dist(f^n(x), \cC) \ge \rho_0 \cdot
  e^{-\rho (n-r)}$ for all $x\in\omega$, where
  $r=r_{\gamma(n-1)}$ is the last return time for $\omega$
  with $n < r + p( f^r(\omega) )$ and
  $\rho_0=1-e^{-\rho}$. Moreover there is no element of
  $\cC$ between $f^n(x)$ and $f^{n-r}(x_{l_r})$, where
  $x_{l_r}$ is the critical point associated to the return
  at time $r$.
\end{proposition}

This result  is enough to conclude that $\cC\cap
f^n(\omega)=\emptyset$, completing the induction step.

\begin{proof}
  We know, from item 1b of Theorem~\ref{th:ref-annals}, that
  for $\mu\in S$, every $h\ge1$ and for all $|l|\ge k_0$
  either
\begin{equation}
  \label{eq:ou}
      |f^{h}(f(x_l))|>\epsilon
      \quad\text{or}\quad
      |f^{h}(f(x_l))-x_{m(h)}|\geq e^{-\rho h},
 \end{equation}
 where $x_{m(h)}$ is the critical point closest to
 $f^{h}(f(x_l))$ as before.  In the former case, if $n$ is a
 binding time for $\omega$, by the the definition of binding
 period, we get for all $x\in\omega$ that
\begin{eqnarray*}
  |f^n(x)|
&\geq& 
|f^{n-r_{\gamma(n-1)}}(x_{l_{\gamma(n-1)}})|
-
|f^n(x)-f^{n-r_{\gamma(n-1)}}(x_{l_{\gamma(n-1)}})|
\\
&\geq& 
\epsilon-\epsilon^{1+\tau} e^{-\tau (n-r_{\gamma(n-1)})}
\ge
\frac{\epsilon}{x_{k_0}} (1-\epsilon^\tau) x_{k_0}
=
2 \frac{1-\epsilon^\tau}{1+e^{-\pi/\beta}} x_{k_0}
> x_{k_0},
\end{eqnarray*}
according to condition~\eqref{eq:condhaty} on the choice of
$k_0$ as a function of $\tau$.


In the latter case in~\eqref{eq:ou}, by definition of
binding~\eqref{BCeq} and because we assume that $\rho<\tau$,
setting $m=m(n-1-r_{\gamma(n-1)})$ for simplicity, we get
that $| f^n(x) - x_{m} |$ is bounded by
\begin{eqnarray}
|f^{n-r_{\gamma(n-1)}}(x_{l_{\gamma(n-1)}})
- x_{m } |
-
| f^n(x) -
  f^{n-r_{\gamma(n-1)}}(x_{l_{\gamma(n-1)}})|
\nonumber
\\
\ge
e^{-\rho(n-1-r_{\gamma(n-1)}) }
- |f^{(n-r_{\gamma(n-1)})}(x_{l_{\gamma(n-1)}})-x_{m}|
\cdot e^{-\tau (n-r_{\gamma(n-1)})} \nonumber
\\
\ge
e^{-\rho(n-r_{\gamma(n-1)}) }\cdot
( 1- \ep )>0.
\qquad\qquad\qquad\qquad  \label{eq:1menosepsilon}
\end{eqnarray}

To complete the proof we consider the case when
$x_{l_{\gamma(n-1)}}$ is not  the closest critical point
to $f^{r_{\gamma(n-1)}}(x)$.  We first argue that no $x'\in\cC$ is
between $f^n(x)$ and
$f^{n-r_{\gamma(n-1)}}(x_{l_{\gamma(n-1)}})$. For
otherwise using \eqref{BCeq} and the definition of
$x_{l_{\gamma(n-1)}}$ we would have
\begin{eqnarray*}
\frac12 \cdot |x' - x_{m}|
&<&
|f^n_\mu(x) - f^{n-r_{\gamma(n-1)}}(x_{l_{\gamma(n-1)}})|
\\
&\le&
|f^{n-r_{\gamma(n-1)}}(x_{l_{\gamma(n-1)}}) -
x_{m}|\cdot e^{-\tau (n-r_{\gamma(n-1)})} 
\\
&\le&
\frac{e^{-\tau (n-r_{\gamma(n-1)})}}2\cdot
|x' - x_{m}|,
\end{eqnarray*}
a contradiction because $e^{-\tau (n-r_{\gamma(n-1)})} <1$.
Hence \emph{there exists $x'\in\cC$ such that $x'$ and
$x_{m}$ are consecutive critical points
in $\cC$ and both $f^n_\mu(x)$ and
$f^{n-r_{\gamma(n-1)}}(x_{l_{\gamma(n-1)}})$ are between
$x'$ and $x_{m}$}. But then
\begin{eqnarray*}
|x' - f^n(x)|
&\ge&
|x' - f^{n-r_{\gamma(n-1)}}  (x_{l_{\gamma(n-1)}}) | -
|f^n(x) - f^{n-r_{\gamma(n-1)}}  (x_{l_{\gamma(n-1)}}) |
\\
&\ge&
\frac12 |x' - x_{m}  |
-
|f^{n-r_{\gamma(n-1)}}(x_{l_{\gamma(n-1)}})
-x_{m}|
 e^{-\tau (n-r_{\gamma(n-1)})}
\\
&\ge&
\frac12 |x' - x_{m} |
-
\frac12\cdot |x' - x_{m}|\cdot
e^{-\tau (n-r_{\gamma(n-1)})}
\\
&\ge&
\frac12 |x' - x_{m} |\cdot
\big(
1- e^{-\tau (n-r_{\gamma(n-1)})}
\big).
\end{eqnarray*}
We observe that since $x'$ and $ x_m $ are consecutive
critical points we have that $x'$ is either $x_{m+1}$ or
$x_{m-1}$, thus
\[
|x'- x_m|
\ge
2 |f^{n-r_{\gamma(n-1)}}(x_{l_{\gamma(n-1)}})  -x_m |
\ge 2 e^{-\rho (n-r_{\gamma(n-1)})}.
\]
Combining the two last inequalities and taking into account
that $\rho<\tau$ gives
\[
|x' - f^n(x)| \ge 
e^{-\rho (n-r_{\gamma(n-1)})}
\cdot
\big(
1- e^{-\tau (n-r_{\gamma(n-1)})}
\big)
\ge
e^{-\rho (n-r_{\gamma(n-1)})}\cdot
\big(
1- e^{-\rho (n-r_{\gamma(n-1)})}
\big).
\]
Choosing $\rho$ and $\epsilon$ close to $0$ such that
$e^{-\rho}>\epsilon$ ($\rho<\log2$ and $\epsilon<1/2$ is
enough) we get $1-\epsilon>1-e^{-\rho}$ and we may then
replace $1-\epsilon$ by $\rho_0=1-e^{-\rho}$
in~\eqref{eq:1menosepsilon}, finishing the proof since $1-
e^{-\rho(n-r_{\gamma(n-1)})} \ge 1-e^{-\rho}=\rho_0$.
\end{proof}



\section{Auxiliary lemmas}
\label{sec:auxiliary-lemmas}

Here we collect some intermediate results needed for the
proofs of the main estimates. In all that follows we write
$C$ for a constant depending only on the initial map $\hat
f$ or $f_0$.


\begin{lemma}{\cite[Lemma 3.1]{PRV98}}
\label{l3.0,5}
Given $\alpha_{1},\alpha_{2},\beta_{1},\beta_{2}$ with
$\frac{\alpha_{1}}{\alpha_{2}}\neq\frac{\beta_{1}}{\beta_{2}}$,
there exists $\delta>0$ such that, for every $x$, at least one of
the following assertions hold:
$|\alpha_{1}\sin x+\beta_{1}\cos x|\geq \delta$
or
$|\alpha_{2}\sin x+\beta_{2}\cos x|\geq \delta.$
\end{lemma}
Using this we obtain the following property of bounded
distortion for the second derivative near critical
points. Observe that since there are inflection points
between consecutive critical points, the lower bound cannot
hold in general, so we restrict this bound to a small scaled
neighborhood of the critical set, reducing the value of
$\tau>0$ if needed.

\begin{lemma}
  \label{le:bobo}
  There exists $\tau>0$ small enough and a constant $C>0$
  depending only on $\hat f$ such that for every $k\ge k_0$
  and $t\in[y_{k+1},y_k]$ we have $|f''(t)|/|f''(x_k)| \le
  C$.
  
  Moreover we can find a constant $C$, depending only on
  $\hat f$ and $\tau$, such that for $|t-x_k|<\tau|x_k|$ we
  have $|f''(t)|/|f''(x_k)|\ge C^{-1}$.
\end{lemma}

Note that \emph{the condition on the lower bound above is
  satisfied by all points in a return situation}, either
essential or inessential, during the construction of the
sequence of partitions $\cP_n$, as detailed in the previous
Section~\ref{sec:concstr-part-bound-distortion}.

\begin{proof} Since $f$ is symmetric we can assume without
  lost that $z>0$ in what follows. We compute
  \begin{align*}
    f'(z)&=-a z^{\alpha-1}[\alpha\sin(\beta\log z)+
    \beta\cos(\beta\log z)] \quad\mbox{and}
    \\
    f''(x)&=-a
    z^{\alpha-2}[A\sin(\beta\log z)+  B\cos(\beta\log z)],
  \end{align*}
  where $A=\alpha(\alpha-1)-\beta^2$ and
  $B=\beta(2\alpha-1)$. Note that for $z=x_l$ we
  have
  \begin{align*}
    0=f^\prime(x_l)&=\alpha\sin(\beta\log x_l)+
    \beta\cos(\beta\log x_l)
  \end{align*}
  and
  \begin{align*}
    \frac{\alpha}{A} \neq \frac{\beta}{B} \iff
    \alpha^2+\beta^2\neq0 \quad\mbox{which is true,
      since}\quad \alpha,\beta>0.
  \end{align*}
  Applying the Lemma~\ref{l3.0,5} we get, because
  $f'(x_k)=0$, that
\[
\frac{|f''(t)|}{|f''(x_k)|}
\le
\Big|\frac{x_{k+1}}{x_k}\Big|^{\alpha-2}
\cdot\frac{|A|+|B|}{\de}
\le
e^{-\frac{\pi}{\beta}(\alpha-2)}
\cdot\frac{|A|+|B|}{\de},
\]
where the last inequality is a direct consequence
of~\eqref{e2.4}.

Now for the lower bound we compute the zeroes of $f^{''}$
and obtain
\begin{align*}
  z_k=\hat z \exp( -k\frac{\pi}{\beta}) \quad\mbox{for}\quad
  k\ge k_0, \quad \mbox{where}\quad \hat z=\exp\Big(
  \frac1{\beta}\tan^{-1}\frac{\beta(2\alpha-1)}{\beta^2-\alpha^2+\alpha}
  \Big).
\end{align*}
By the expressions for the zeroes of $f^{'}$ and the zeroes
of $f^{''}$ we see that their distance is scaled by a common
factor. Thus if we consider a scaled neighborhood of each
$x_l$ specified by a small enough $\tau>0$, that is,
considering only $t$ with $|t-x_l|<\tau|x_l|$, then we
ensure that $t$ is far away from the inflection point.

Moreover the value of the quotient in the statement of the
lemma is invariant under the scaling: let $l>k\ge k_0$ and
$|t-x_k|<\tau|x_k|$. Then on the one hand
\begin{align*}
  e^{-(l-k)\frac{\pi}{\beta}}|t-x_k|=|s-x_l|<\tau|x_l|
  \quad\mbox{where}\quad
  s=t\cdot e^{-(l-k)\frac{\pi}{\beta}}.
\end{align*}
On the other hand
\begin{align*}
  \frac{|f''(s)|}{|f''(x_l)|}
  &=
  \frac{\big|a s^{\alpha-2}[A\sin(\beta\log s)+
    B\cos(\beta\log s)]\big|}
  {\big|a x_l^{\alpha-2}[A\sin(\beta\log x_l)
    +B\cos(\beta\log x_l)]\big|}
  \\
  &=
  \frac{\big| 
    a\big(t\cdot e^{-(l-k)\frac{\pi}{\beta}}\big)^{\alpha-2}
    [A\sin(\beta\log t + (k-l)\pi)+
    B\cos(\beta\log t + (k-l)\pi)]\big|}
  {\big|a \big(x_k \cdot
    e^{-(l-k)\frac{\pi}{\beta}}\big)^{\alpha-2}
    [A\sin(\beta\log x_k + (k-l)\pi)
    +B\cos(\beta\log x_k+ (k-l)\pi)]\big|}
\end{align*}
is the same value of $|f''(t)|/|f''(x_k)|$. Therefore the
lower bound is given by
\begin{align*}
  \inf\Big\{  \frac{|f''(t)|}{|f''(x_k)|} : |t-x_k|<\tau|x_k|\Big\}
\end{align*}
for any $k\ge k_0$. Taking $C>0$ big enough we conclude the proof.
\end{proof}

The next result guarantees that orbits of points in
$[-C\epsilon^\alpha,-\tilde y]\cup[\tilde
y,C\epsilon^\alpha]$ remain expanding during a number $m_0$
of iterates that can be fixed arbitrarily large by reducing
$\epsilon$ and $\mu$.  Recall that $|\tilde{f}'| > \sigma_0
\gg 1$ and that $f$ is $C^1$-close to $\tilde{f}$
outside $[-\tilde{y},\tilde{y}]$ if $\mu$ is small.

\begin{lemma}{\cite[Lemma 6.1]{PRV98}}
\label{le:emezero}
There exist $c,C>0$ and $m_0 \ge c\,\log(1/\epsilon)$ such
that if $\tilde y<|x|\le C\epsilon^\alpha$, then
$f^{i}(x) \notin [-\tilde{y},\tilde{y}]$ and
$|f'(f^{i}(x))|\geq\sigma_0 ,$ for all $1\le
i\le m_0$ and all $\mu\in[-\epsilon,\epsilon]$.
\end{lemma}

Next we establish some results of bounded distortion and
uniformly bounded expansion during binding periods.

\begin{lemma}[Bounded distortion on binding periods]
  \label{le:distorcao}
  There exists $A=A(C,\tau)>1$ such that for all $x\in
  I(l,s)$ we have
\[
\frac{1}{A}\leq\left|\frac{(f^{j})'(\xi)}
{(f^{j})'(f(x_{l}))}\right|\leq A
\]
for every $1\leq j\leq p(l,s)$ and every
$\xi\in[f(x_l),f(x)]$.
\end{lemma}

\begin{proof}
We let $\eta=f(x_{l})$ and consider $0\le i <j$.
There are two cases to treat, corresponding to the two possibilities
in (\ref{BCeq}).
If $|f^i(\eta)|\le\epsilon$ then, by Lemma \ref{l3.1,5},
\[
\left|\frac{f'(f^{i}(\xi))-f'(f^{i}(\eta))}
           {f'(f^{i}(\eta))}\right|
\leq  C\frac{|f^{i}(\xi)-f^{i}(\eta)|}
           {|f^{i}(\eta)-x_{m(i)}|}
\leq C e^{-\tau i}.
\]
If $|f^i(\eta)|>\epsilon$, then
$|f^i(\xi)-f^i(\eta)|\le\epsilon^{1+\tau} e^{-\tau i}\ll\epsilon$
and so the interval bounded by $f^i(\xi)$ and $f^i(\eta)$ is
contained in the region $S^1\setminus[-\tilde{y},\tilde{y}]$, where
$f=\tilde{f}$.
Thus,
\[
\left|\frac{f'(f^{i}(\xi))-f'(f^{i}(\eta))}
           {f'(f^{i}(\eta))}\right|
\leq C |f^{i}(\xi)-f^{i}(\eta)|
\leq C \epsilon^{1+\tau} e^{-\tau i}
\leq C e^{-\tau i}.
\]
Putting together all the above we get
$
\sum_{i=0}^{j-1}\left|\frac{f'(f^{i}(\xi))-f'(f^{i}(\eta))}
                           {f'(f^{i}(\eta))}\right|
\leq C \sum_{i=0}^{j}e^{-\tau i}
\leq C.
$
Thus
\[
  \log\left|\frac{(f^{j})'(\xi)}
{(f^{j})'(\eta)}\right| 
\le
\sum_{i=0}^{j-1}\log\left(1+
\left|
\frac{f'(f^i(\xi))}{f'(f^i(\eta))} -1
\right|
\right)
 \le
\sum_{i=0}^{j-1} \left|
\frac{f'(f^i(\xi))}{f'(f^i(\eta))} -1
\right|
\le C,
\]
and the statement of the lemma follows.
\end{proof}

Now we prove an exponential bound on the derivative along
the orbit of each critical value $z_k, |k|\ge k_0$.  This is
needed to get a lower bound for the binding time $p$ in
terms of the position of the return interval given by
$(l,s)$ in the following Lemma~\ref{le:maisdistorcao}.  This
demands a proof since the derivative of $f$ is unbounded
due to the presence of infinitely many critical points,
unlike the quadratic family where we have this property for
free.

In our setting we will obtain this by imposing an extra
condition, besides item (1b) from Theorem~\ref{th:ref-annals},
in the construction of the set $S$ of parameters $\mu$ which
we shall consider in the proof of
Theorems~\ref{thm:approxestimate} and~\ref{thm:expansionestimate}.
This condition is expressed by the inequality
\begin{equation}
  \label{eq:CPR0}
  \sum_{j=0}^{n-1} -\log \dist\big( f^j_\mu(z_k)
  ,\cC\big)\le 
  \hat M\cdot n,
\end{equation}
for all $n\ge1$ which are \emph{not bound times} for the
orbit of the critical value $z_k$, for every $|k|\ge k_0$
and for some big fixed constant $\hat M>0$.  Since for all
$x\in I$ we have $| x | \ge \dist\big( x ,\cC\big)$ this
implies
\begin{equation}
  \label{eq:CPR1}
  \sum_{j=0}^{n-1} -\log \big| f^j_\mu(z_k) \big|\le \hat M\cdot n
\end{equation}
and we are able to deduce the following.

\begin{lemma}
  \label{le:CPR1}
Assume that \eqref{eq:CPR1} holds for some $\mu\in S$. Then
there exists a constant $M>0$ such that 
$\big| (f^n)'(z_k) \big| \le M^n$
for all $n\ge1$ and
$|k|\ge k_0$.
\end{lemma}

\begin{proof}
  Note that just by taking the derivative of $f$ we get
  that there exists a constant $C>0$ such that $\big|
  f'(x) \big| \le C \cdot |x|^{\alpha-1}$.  On the one
  hand, for $n\ge1$ such that $n$ is \emph{not a bound time}
  for $z_k$ and $|k|\ge k_0$ we use \eqref{eq:CPR1} to get
  \begin{align*}
    \big| (f^n)'(z_k) \big| & 
\le \prod_{j=0}^{n-1} C\cdot\big|
f^j_\mu(z_k)\big|^{\alpha-1}
 \le
\exp\Big( \sum_{j=0}^{n-1} \big( \log C +
(\alpha-1)\cdot\log \big|
f^j_\mu(z_k)\big| \big) \Big)
\\
& \le
\exp\Big( n\cdot \log C + (1-\alpha)\cdot \hat M \cdot n
\Big)
= \tilde M^n,
  \end{align*}
  where $\tilde M=\exp\big( \log C + (1-\alpha)\cdot \hat M
  \big)$.  We can assume without loss that $\tilde M>A$
  where $A>1$ is given by Lemma~\ref{le:distorcao}.

  On the other hand, if $n$ is a \emph{bound time} for
  $z_k$, let $t_1<n$ be the return time for $z_k$ which
  originated the binding and $z_{l_1}$ be the corresponding
  bound critical value. Then by Lemma~\ref{le:distorcao}
\[
\big|
(f^n)'(z_k)
\big|
\le 
A\cdot  \big| (f^{n-t_1})'(z_{l_1})\big|
\cdot\big| (f^{t_1})'(z_k)\big|.
\]
\emph{If $t_1$ is not a bound time for $z_{l_1}$}, then by
the bound just proved for all critical values on free times
and return situations we bound the last expression by $A
\tilde M^{n-t_1} \tilde M^{t_1}=A\cdot \tilde M^n$.
Otherwise there are $t_2<n-t_1$ the return time for
$z_{l_1}$ which originated the binding and $z_{l_2}$ the
corresponding bound critical value and, as above, we get
\[
\big|
(f^n)'(z_k)
\big|
\le
A ^2 \cdot \big| (f^{n-t_1-t_2})'(z_{l_2})\big|
\cdot \big| (f^{t_2})'(z_{l_1})\big|
\cdot\big| (f^{t_1})'(z_k)\big|.
\]
\emph{If $n-t_1-t_2$ is not a bound time for $z_{l_2}$} we
bound this expression by $A^2 \tilde M^{n-t_1-t_2} \tilde
M^{t_2} \tilde M^{t_1}\le A^2 \tilde M^n$. Otherwise we
repeat the argument.  Knowing that \cite[Beginning of Section
4, page 450]{PRV98} the orbit of every critical value has an
initial number $j_0\gg 1$ of free iterates, we see that this
argument must end in a free time and we arrive at $\big|
(f^n)'(z_k) \big|\le A^l \tilde M^n$ where $l<n$ is the
number of nested binding periods obtained. Since $\tilde
M>A$ we conclude the proof of the lemma by setting $M=\tilde
M^2$.
\end{proof}

Now we obtain the estimates
for the binding time assuming that \eqref{eq:CPR1} holds for
$\mu\in S$. In Section~\ref{sec:extra-excl-param} we explain
how to obtain \eqref{eq:CPR1} for a positive measure subset
of parameters $S$.

\begin{lemma}[Expansion during binding periods]
  \label{le:maisdistorcao}
  There are constants $A_0=A_0(\ep,\rho,\tau)>1$, $\iota=\iota(M)>0$ and
  $\theta=\theta(M,\epsilon,\rho,\tau)\in\bN$ such that for
  $n\ge1$ and $\omega\in\cP_n$ with
  $R_n(\omega)\neq\emptyset$, if $r$ is the last return time
  for $\omega$ and $f^{r}(\omega)\subset I(l,s,j)$, then
  setting $p=p(l,s)>0$ and $\zeta=2(\rho+\tau)/\log\sigma$ we
  have for $\tau$ small enough
\begin{enumerate}
\item[(a)] $\iota(M)\cdot (|l|+|s|) \le
  p\leq \frac{2\pi}{\beta\log\sigma}(|l|+|s|)$;
\item[(b)] $|(f^{p+1})'(f^{r}(x))| \ge
  \frac1C\cdot \ep^{1+\tau}\cdot
  \exp\Big((1-\zeta)\frac{\pi}{\beta}( |l|+|s|)\Big)$ and if
  $|l|+|s|\ge\theta$, then
  $|(f^{p+1})'(f^{r}(x))|\ge
  \frac1C\cdot\exp\Big(\frac{\pi}{\beta}(
  |l|+|s|)\Big)$, for every $x\in\omega$;
\item[(c)]
$|(f^{p+1})'(f^{r}(x))|
\geq A_0\cdot \sigma^{(p+1)/3}>2
$
for every $x\in\omega$.
\end{enumerate}
\end{lemma}

\begin{proof}
To prove item (a), we use the definition of the partition
and the construction of the refinement.
As $p>0$, 
$(l,s,j)\neq(\pm k_0,1,1)$ and so 
$|f^{r}(x)-x_{l}|
\geq a_{2}\cdot e^{-(\pi/\beta)(|l|+|s|)},
$
by Remark~\ref{rmk:partitiondimension},
where $x$ is any given point in $\omega$.
Using the second-order Taylor approximation and
Lemma~\ref{le:bobo} we get
\[
|f^{r+1}(x)-f(x_{l})| 
\geq
\frac{1}{C}|f''(x_{l})| (a_{2}e^{-(\pi/\beta)(|l|+|s|)})^2
\geq 
\frac{1}{C} \cdot 
e^{-\frac{\pi}{\beta}|l|(\alpha-2)} \cdot
e^{-2(\pi/\beta)(|l|+|s|)},
\]
where $|f''(x_l)|\ge C^{-1}|x_l|^{\alpha-2}= C^{-1}\cdot
{\hat x}^{\alpha-2}\cdot
e^{-\frac{\pi}{\beta}|l|(\alpha-2)}$ by Lemma~\ref{l3.1}(2).
Then for each $0\le j\leq p$, there is some $\xi$ between
$f(x_l)$ and $f^{r+1}(x)$ such that
\begin{eqnarray}
|f^{j+r+1}(x)-f^{j+1}(x_{l})|
&=&
|(f^{j})'(\xi)|\cdot|f^{r+1}(x)-f(x_{l})|
\nonumber
\\
&\geq&
C^{-1}  \cdot 
e^{-\frac{\pi}{\beta}|l|(\alpha-2)} \cdot
e^{-2(\pi/\beta)(|l|+|s|)}\cdot|(f^{j})'(\xi)|.\label{eq:mil}
\end{eqnarray}
Now since $\alpha-2<0$ and we can take $|l|\ge k_0$ very
big, as a consequence of Lemma~\ref{le:distorcao} and of the
exponential growth of the derivative at the critical orbits,
we get the following bound
\[
2 \cdot e^{-2(\pi/\beta)(|l|+|s|)}\sigma^{j}
\le
|f^{j+r+1}(x)-f^{j+1}(x_{l})|
\leq 2
\]
Hence
$e^{-2(\pi/\beta)(|l|+|s|)}\sigma^{j}\leq 1$
for all $1\leq j\leq p$.
In particular, 
$$
-2(\pi/\beta)(|l|+|s|)+ p \log(\sigma)\leq 0,
\quad\text{implying}\quad
p\leq\frac{2(\pi/\beta)(|l|+|s|)}{\log\sigma},
$$
thus proving the upper bound in (a).

For the lower bound in (a), we note that by the definition
of binding period, we have that
\begin{eqnarray}
|f^{r+p+1}(x)-f^{p+1}(x_{l})|
\ge |f^{p+1}(x_{l})-x_{m(p)}| \cdot e^{-\tau (p+1)}
&\text{if}&
|f^{p+1}(x_l)|\le\epsilon \nonumber
\\
& \text{and} \label{eq:nomeio}
\\
|f^{p+r+1}(x)-f^{p+1}(x_{l})|
\geq \epsilon^{1+\tau} \cdot e^{-\tau (p+1)}
&\text{if}&
|f^{p+1}(x_{l})|>\epsilon.   \nonumber
\end{eqnarray}
So in either case using Lemmas~\ref{le:distorcao}
and~\ref{le:CPR1} we get
\begin{align*}
|f^{p+r+1}(x)-f^{p+1}(x_{l})|
& =
\big| (f^p)'(\xi) \big| \cdot | f^{r+1}(x)-f(x_l) | 
 \le
A\cdot \big| (f^p)'(z_l) \big| \cdot |f^{r+1}(x)-f(x_l)| 
\\
&\le
A\cdot M^p \cdot
C |x_l|^{\alpha-2}|f^r(x)-x_l|^2
\le
AM^p C e^{-\frac{\pi}{\beta}(\alpha|l|+2|s|)}
\\
&\le
AM^p e^{-\alpha\frac{\pi}{\beta}(|l|+|s|)},
\end{align*}
where $\xi$ is some point between $f^{r+1}(x)$ and $f(x_l)$
and in the second line above we used Lemma~\ref{l3.1}.

On the one hand, if $|f^{p+1}(x_{l})|>\epsilon$ then
we arrive at $\epsilon^{1+\tau} e^{-\tau (p+1)} \le A
M^p  e^{-\alpha\frac{\pi}{\beta}(|l|+|s|)}$.  On the other hand, if
$|f^{p+1}(x_l)|\le\epsilon$ then by condition (1b)
from Theorem~\ref{th:ref-annals} we arrive at $ e^{-\rho
  \cdot p} e^{-\tau (p+1)} \le A M^p
e^{-\alpha\frac{\pi}{\beta}(|l|+|s|)}. $ In both cases if we take $M$ big
enough, then we get a bound of the form $p\ge \iota(M)\cdot(
|l|+|s|)$ concluding the proof of (a).

Now we prove (b). Since $p+1$ is not a binding time we
must have (\ref{eq:nomeio}) as before.
Using Theorem~\ref{th:ref-annals}(1b), setting
$\Delta_{p+1}=\min\{\ep^{1+\tau}\cdot e^{-\tau(p+1)},
e^{-(\rho+\tau)(p+1)} \}$ we get for some
$\xi\in[f(x_{l}),f^{r+1}(x)]$, by Lemma~\ref{le:distorcao},
and using second-order Taylor expansion of $f$ near $x_l$
together with the upper bound from Lemma~\ref{le:bobo}
\begin{eqnarray}
\Delta_{p+1}
&\le& |f^{p+1}(x_{l})-f^{p+r+1}(x)|
=
|(f^{p})'(\xi)|\cdot|f(x_{l})-f^{r+1}(x)|\nonumber
\\
& \leq&
 C \cdot |(f^{p})'(f^{r+1}(x))|\cdot
|f''(x_{l})|
\cdot |f^{r}(x)-x_l|^{2}\label{eq:mileum}.
\end{eqnarray}
Note that the second order Taylor expansion near $x_l$ gives
\begin{align*}
  \big|f(x_l)-f^{r+1}(x)\big|=
  \big|f^{'}(x_l)\cdot(x_l-f^r(x)) +
  \frac{f^{''}(\xi)}2\cdot (x_l-f^r(x))^2\big| =
  \big|\frac{f^{''}(\xi)}2\cdot (x_l-f^r(x))^2\big|
\end{align*}
for some $\xi$ between $x_l$ and $f^r(x)$. Together
with Lemma~\ref{le:bobo} we obtain the bound
in~\eqref{eq:mileum}.

On the other hand, using Lemma~\ref{l3.1} and again
Lemma~\ref{le:distorcao} we get
\[
|(f^{p+1})'(f^{r}(x))|
=
| f'(f^{r}(x)) | \cdot
|(f^{p})'(f^{r+1}(x))|
\ge
C^{-1} |x_{l}|^{\alpha-2}\cdot  |f^{r}(x)
-x_{l}| \cdot |(f^{p})'(f^{r+1}(x))|.
\]
Hence by the previous expression together
with~(\ref{eq:mileum}) we deduce
\[
|(f^{p+1})'(f^{r}(x))|
\ge
\frac{\ep^{1+\tau}\cdot e^{-(\rho+\tau)(p+1)}}{C\cdot |f^{r}(x)
-x_{l}|}
\ge
\frac1C\cdot \ep^{1+\tau}\cdot e^{-(\rho+\tau)(p+1)}
\cdot e^{\frac{\pi}{\beta}( |l|+|s|)}
\]
Now we have two possibilities, either
$\epsilon>e^{-(\rho+\tau)(p+1)}$ or not. In the former case
we obtain by the upper bound in item (a) and because we can
assume that $2+\tau\le3$ and $\rho+\tau\le1$
\[
|(f^{p+1})'(f^{r}(x))| \ge
e^{-(2+\tau)(\rho+\tau)(p+1)} \cdot e^{\frac{\pi}{\beta}(
  |l|+|s|)} \ge
\frac1C\exp\Big(
\frac{\pi}{\beta}(
|l|+|s|)\Big).
\]
In the latter case we get $|(f^{p+1})'(f^{r}(x))| \ge C^{-1}
\ep^{1+\tau}\cdot
\exp\Big(\big(1-\frac{\rho+\tau}{\log\sigma}\big)\frac{\pi}{\beta}(
|l|+|s|)\Big).$ Observe that in this case
$\log\epsilon\le-(\rho+\tau)(p+1)$, thus by the lower bound
in item (a) and since $k_0=C\log(1/\epsilon)$ we have
\[
(\rho+\tau)\cdot\iota(M)\cdot(|l|+|s|+1)\le
(\rho+\tau)\cdot(p+1) \le C\cdot k_0.
\]
So if $|l|+|s|\ge\theta$ with
$(\rho+\tau)\cdot\iota(M)\cdot(\theta+1)>C\cdot k_0$ then
only the first alternative can happen.
This concludes the proof of item (b).


In order to prove (c) we use Lemma~\ref{le:distorcao} once
again and the Mean Value Theorem applied to $f^{'}$ near
$x_l$ together with the lower bound from Lemma~\ref{le:bobo}
to get
\begin{align}
  |(f^{p+1})'(f^{r}(x))|^{2}
  & = 
  |(f^{p})'(f^{r+1}(x))|^{2}\cdot
  |f'(f^{r}(x))|^{2} \nonumber
  \\ 
  &\geq \big(C^{-1} |(f^{p})'(f(x_{l}))|\cdot
  |(f^{p})'(f^{r+1}(x))|\big)\cdot\big(
  C^{-1} |f''(x_{l})|^{2} \cdot
  |f^{r}(x)-x_{l}|^{2}\big)\nonumber
  \\
  &= C^{-1} |(f^{p})'(f(x_{l}))|
\cdot|f''(x_{l})|
\cdot \big(|(f^{p})'(f^{r+1}(x))|\cdot
|f''(x_{l})| \cdot |f^{r}(x)-x_{l}|^{2} \big).\label{eq:nova}
\end{align}
Note that we can use the lower bound from
Lemma~\ref{le:bobo} because $r$ is a return time, not an
escape nor a free time.  Comparing the last
expression~\eqref{eq:nova} with~\eqref{eq:mileum} we see
that~\eqref{eq:nova} is bounded below by
\begin{align*}
&\ge C^{-1} |(f^{p})'(f(x_{l}))|
   \cdot|f''(x_{l})|
\cdot \big(
   |(f^{p})'(\xi)|\cdot|f(x_{l})-f^{r+1}(x)|
 \big)
\\
&= C^{-1} |(f^{p})'(f(x_{l}))| \cdot |f''(x_{l})|
        \cdot |f^{p+1}(x_{l})-f^{r+p+1}(x)|
      \end{align*}
Finally from  Lemma~\ref{l3.1}(2) we deduce
\begin{align}      
\geq C^{-1} \sigma^{p} \cdot C^{-1} |x_l|^{\alpha-2}\cdot 
       |f^{p+1}(x_{l})-f^{r+p+1}(x)|\label{eq:miledois}.
  \end{align}

  Now we consider two cases. 

  On the one hand, if $|f^{p+1}(x_l)|>\epsilon$ then, by the
  definition of $p$ in (\ref{BCeq}), we must have $
  |f^{p+1}(x_{l})-f^{r+p+1}(x)|
  >\epsilon^{1+\tau}e^{-\tau(p+1)}$.
  Equation~\eqref{eq:miledois} together with $\alpha-2<0$
  and $|x_l|\le\ep$ imply that $|(f^{p+1})'(f^{r}(x))|^{2}\ge
  C^{-2}\sigma^p \cdot \epsilon^{\alpha-2}
  \cdot|f^{p+1}(x_{l})-f^{r+p+1}(x)|$. Now we can write $
  |(f^{p+1})'(f^{r}(x)|^{2} \ge
  C^{-1}\sigma^{p}\epsilon^{\alpha-1+\tau}e^{-\tau(p+1)} \ge
  A_0^2\cdot \sigma^{2(p+1)/3}, $ if we fix
  $\tau<\min\{1-\alpha,\log\sigma/3\}$ and take
  $\epsilon$ small enough.

  On the other hand, if $|f^{p+1}(x_l)|\le\epsilon$ then $
  |f^{p+1}(x_{l})-f^{r+p+1}(x)| >|f^{p+1}(x_l)-x_{m(p)}|
  e^{-\tau(p+1)} $ by (\ref{BCeq}).

  Finally we note that there is only one possibility
  according to item (1b) of Theorem~\ref{th:ref-annals},
  that is $|f^{p+1}(x_{l})-x_{m(p)}|\geq e^{-\rho p }$, and
  thus $ |f^{p+1}(x_{l})-f^{r+p+1}(x)| > C^{-1}\epsilon\,
  e^{-(\rho+\tau)(p+1)}.  $ Hence
\[
|(f^{p+1})'(f^{r}(x)|^{2}
\ge  C^{-1} \cdot \sigma^{p}\cdot \epsilon^{\alpha-2}
\cdot e^{-(\rho+\tau)(p+1)}
\ge A_0^2\cdot \sigma^{2(p+1)/3}
\]
as long as we take $\epsilon$, $\rho$ and $\tau$ small enough.
This concludes the proof of the lemma.
\end{proof}

\begin{remark}
  \label{rmk:zetaM}
Note that $\iota(M)\to0$ when $M\to\infty$.
\end{remark}

Now we will obtain estimates of the length of $|f^n(\omega)|$.

\begin{lemma}[Lower bounds on the length at return times]
\label{lem:fn-omega-estimate}
Let $r$ be an essential or an inessential return time for
$\omega\in\mathcal{P}_{n-1}$, with host interval
$I(l,s,j)^+$. Let $p=p(l,s)$ denote the length of its
binding period and set $Q=Q(l,s,\tau,\rho)=
\frac{\epsilon^{1+\tau}}{(|l|+|s|)^3}\cdot
e^{\zeta\frac{\pi}{\beta}(|l|+|s|)}$.
Then the following holds:
\begin{enumerate}
\item Assuming that $r^*\leq n-1$ is the next return
  situation for $\omega$ (either essential, inessential or
  an escape) we have $\left|f^{k}(\omega)\right|\geq Q\cdot
  \sigma_0^q \cdot e^{(1-2\zeta)
    \cdot\frac{\pi}{\beta}(|l|+|s|)} \cdot |f^{r}(\omega)|$,
  for $q=k-(r+p+1)$ and every $k$ such that $r+p+1\le k \le
  r^*$.

  Moreover we also have $\left|f^{k}(\omega)\right|\geq
  \sigma_0^q \cdot A_0\cdot\sigma^{(p+1)/3} \cdot
  |f^{r}(\omega)| \ge A_0 \cdot |f^{r}(\omega)| >
  |f^{r}(\omega)|$.
        
\item If $r$ is the last return time for $\omega$ up to
  iterate $n-1$ and also an \emph{essential} return, and
  $r^*$ is a return time for $\omega$, then setting
  $q=k-(r+p+1)$ we have
  \[
  \left|f^{k}(\omega)\right|\geq a_1\cdot
  Q\cdot\sigma_0^q \cdot
  e^{-2\zeta\frac{\pi}{\beta}(|l|+|s|)}\quad
  \mbox{for all}\quad r+p+1\le k \le r^*.
  \]
\end{enumerate}
Suppose that $r$ is an escape time for
$\omega\in\mathcal{P}_{n-1}$.
\begin{enumerate}
\item[(3)] If $r^*\le n-1$ is the next return situation for
  $\omega$, then $\left|f^{k}(\omega)\right|\geq
  \sigma_0^{k-r}\cdot \left|f^{r}(\omega)\right|$, for every
  $k$ such that $r< k \le r^*$.
\end{enumerate}
\end{lemma}

\begin{remark}
  \label{rmk:Qmaiorque1}
Note that $Q=Q(l,s,\tau,\rho)\to\infty$ when $|l|+|s|\to\infty$.
\end{remark}

\begin{proof}
  We start by assuming $r^*\le n-1$. 
  By the mean value theorem we have $|f^{r^*}(\omega)|\geq
  |(f^{r^*-r})'(f^r(\xi))|\cdot |f^{r}(\omega)|$ for some
  $\xi\in\omega$.  Using Remark~\ref{rmk:expansaolivre} and
  Lemma~\ref{le:maisdistorcao} we get by setting
  $q=r^*-(r+p+1)$
\begin{eqnarray}
  \left|f^{r^*}(\omega)\right|
  &\geq&
  \left|\left(f^{q}\right)'\left(f^{r+p+1}(\xi)\right)\right|
\cdot
  \left|\left(f^{p+1}\right)'\left(f^r(\xi)\right)\right|
\cdot
  \left|f^{r}(\omega)\right| \label{eq:outra}
  \\
  &\geq&  \sigma_0^q \cdot 
  \frac1C\cdot \ep^{1+\tau}\cdot
  e^{(1-\zeta)\frac{\pi}{\beta}(|l|+|s|)}
  \cdot |f^{r}(\omega)| \nonumber
  \\
  &\ge&  \sigma_0^q \cdot
  \frac{\ep^{1+\tau}}{C}\cdot
  e^{\zeta\frac{\pi}{\beta}(|l|+|s|)}
  \cdot
  e^{(1-2\zeta)\frac{\pi}{\beta}( |l|+|s|)}
  \cdot |f^{r}(\omega)|. \label{eq:valor2}
\end{eqnarray}
Taking into account the definition of $Q$ in
\eqref{eq:valor2} above we obtain the first part of item
(1).

If $r$ is an essential return time for $\omega$, then
$I(l,s,j)\subset f^r(\omega)$ and $|f^r(\omega)| \geq a_{1}
\frac{e^{-(\pi/\beta)(|l|+|s|)}} {(|l|+|s|)^3}$, hence
\[
\left|f^{r^*}(\omega)\right| \geq \sigma_0^q \cdot
\frac{\ep^{1+\tau}\cdot
  e^{\zeta\cdot\frac{\pi}{\beta}(|l|+|s|)}} {C\cdot
  (|l|+|s|)^3} \cdot a_{1} e^{-2\zeta\frac{\pi}{\beta}(|l|+|s|)},
\]
and by the definition of $Q$ this proves item (2).

For the second part of item (1) just use the inequality from
Lemma~\ref{le:maisdistorcao}(c) in~(\ref{eq:outra}).

To get item (3) observe that between the iterate $r$ and
$r^*-1$ there are only free iterates for $\omega$, thus the
estimate follows by the uniform expanding rate on free
times.  

Altogether this concludes the proof of the lemma for the
case $k=r^*$. For the other cases $r+p<k\le r^*$ observe
that only the number of free iterates after the last bound
iterate until $k$ is affected, and this number equals
$q=k-(r+p+1)$.
\end{proof}


\begin{lemma}[Bounded Distortion]
\label{lem:bounded-distortion}
There is a constant $D_0=D_0(\rho,\tau,\sigma)>0$ such that for
$\omega\in\mathcal{P}_{n-1}$, $n\in\bN$, and for
every $x,y\in \omega$ we have
$\big|(f^n)'(x)\big|/\big|(f^n)'(y)\big|
\leq D_0.$
\end{lemma}

\begin{proof}
  Let $R_{n-1}(\omega)=\left\{r_1,\ldots,r_{\gamma}\right\}$
  and $Q_{n-1}(\omega)= \left\{(l_1,s_1,j_1),\ldots,
    (l_{\gamma},s_{\gamma},j_{\gamma})\right\}$, be the sets
  of return situations (essential returns, inessential
  returns and escapes) and indexes of host intervals of $\omega$,
  respectively, as defined during the construction of the
  partition.  Let $\omega_i=f^{r_i}(\omega)$,
  $p_i=p(l_i,s_i)$ for $i=1,\ldots,\gamma$ and, for
  $y,z\in\omega$, let $y_k=f^k(y)$ and $z_k=f^k(z)$
  for $k=0,\dots,n-1$.  Observe that $\omega_i\subset
  I(l_i,s_i,j_i)^+$ for all $i$ and
\begin{equation}
  \label{eq:Distortion}
\left| \frac{(f^n)^\prime(z)}{(f^n)^\prime (y)}
\right| 
= \prod^{n-1}_{k=0} \left|
\frac{f^\prime (z_k)}{f^\prime (y_k)}
 \right|
\le 
\prod^{n-1}_{k=0} \left( 1+
\left| \frac{f^\prime(z_k)-f^\prime(y_k)}{f^\prime(y_k)}
\right| \right).
\end{equation}
On free iterates, if $y_k\in[-\ep,\ep]$, then by
Lemma~\ref{l3.1,5}
\begin{equation}
  \label{eq:quotient1}
\left|
\frac{f^\prime(z_k)-f^\prime(y_k)}{f^\prime(y_k)}
\right| 
\le K_1\cdot\left| \frac{z_k-y_k}{y_k-\tilde x_{k}} \right|
\le K_1\cdot\frac{|f^k(\omega)|}{\De_k(\omega)},
\end{equation}
where we define $\De_k(\omega)=\dist(f^k(\omega),\cC)=
\inf_{x\in\omega}\dist(f^k(x),\cC)$ and $\tilde x_{k}$
is the critical point closest to $y_k$.  We observe that in
this case the interval $f^k(\omega)$ is between two
consecutive critical points, $x_{l_k}$ and $x_{l_k+1}$, and
is contained in some $I(l_k,s_k,j_k)^+$ with $s_k\le s(\tau)$.
Note that by the exponential character of the
initial partition, we have
\begin{equation}
  \label{eq:tosse}
|f^k(\omega)|\le C\cdot |I(l_k, s_k,1)^+|  \quad
\mbox{and}\quad \De_k(\omega)\ge C^{-1}\cdot |I( l_k,  s_k)|
\end{equation}
for some constant $C>0$ depending only on $\hat f$ (see
Remark~\ref{rmk:partitiondimension}). 

Otherwise for free iterates $y_k\in S^1\setminus[-\ep,\ep]$
we get
\begin{eqnarray}
  \label{eq:cha}
\sum_{\genfrac{}{}{0pt}{}{r_i+p_i<k<r_{i+1}}{|y_k|>\ep}}
\left|
\frac{f^\prime(z_k)-f^\prime(y_k)}{f^\prime(y_k)}
\right| 
\le
\frac{L}{\tilde\sigma}
\sum_{\genfrac{}{}{0pt}{}{r_i+p_i<k<r_{i+1}}{|y_k|>\ep}}|z_k - y_k|
\le
\frac{L}{\tilde\sigma}
 \sum_{\genfrac{}{}{0pt}{}{r_i+p_i<k<r_{i+1}}{|y_k|>\ep}}
|f^k(\omega)| \nonumber
\\
\le
\frac{L}{\tilde\sigma}\sum_{r_i+p_i<k<r_{i+1}}
\sigma_0^{k-r_{i+1}}\cdot |f^{r_{i+1}}(\omega)|
\le
K_2 \cdot \frac{|\omega_{i+1}|}{\De_{r_{i+1}}(\omega)}
\end{eqnarray}
by definition of $f$ on $S^1\setminus[-\ep,\ep]$, since
$|f' \mid S^1\setminus[-\ep,\ep]|>\tilde\sigma$ and
$|f''\mid S^1\setminus[-\ep,\ep]|\le L$ for some
constant $L$. We recall also that $\De_{r_{i+1}}(\omega)<1$ by
definition.

\medskip

For an escape time $k=r_i$ with $i\in\{1,\dots,\gamma\}$ we
have either $|y_k|\le\epsilon$ and then we have the
inequality \eqref{eq:quotient1}, or $|y_k|>\epsilon$ and we
get as in \eqref{eq:cha}
\begin{equation}
  \label{eq:naoescapas!}
\left|
\frac{f^\prime(z_k)-f^\prime(y_k)}{f^\prime(y_k)}
\right| 
\le
\frac{L}{\tilde\sigma} \cdot |z_k - y_k|
\le
\frac{L}{\tilde\sigma} \cdot |f^k(\omega)|
\le 
K_3 \cdot \frac{|\omega_{i}|}{\De_{r_{i}}(\omega)}.
\end{equation}
Hence up to now all cases are bounded by the same type of
expression.

\medskip

Next we find a bound for iterates during binding times.  Let
us fix $i=1,\dots,\gamma$ such that $r_i$ is \emph{not an
  escape time for $\omega$}, i.e. it is either an essential
or inessential return time. Then for $k=r_i$ we have the
same bound~\eqref{eq:quotient1}. For $r_i<k\le r_i+p_i$ we
get for some $\xi\in\omega$
\begin{align*}
  |z_k-y_k|
  &=
  |(f^{k-r_i})'(f^{r_i}(\xi))|
  \cdot|z_{r_i}-y_{r_i}|
  \le
  |(f^{k-r_i})'(f^{r_i}(\xi))|
  \cdot |f^{r_i}(\omega)|
  \\
  &=
  |(f^{k-r_i-1})'(f^{r_i+1}(\xi))|\cdot
  |f^\prime(f^{r_i}(\xi)) - f^\prime(x_{l_i})|\cdot
  |\omega_i|
  \\
  & \leq  C\cdot |(f^{k-r_i-1})'(f^{r_i+1}(\xi))|\cdot
  |f''(x_{l_i})| \cdot |f^{r_i}(\xi)-x_{l_i}|
  \cdot |\omega_i|,
\end{align*}
where we have used the Mean Value Theorem applied to $f'$
near the critical point $x_{l_i}$ together with the upper
bound from Lemma~\ref{le:bobo}.

We now have two possibilities by definition of $p_i$. On the
one hand, for the first case in~\eqref{BCeq} there exists
$w\in [f(x_{l_i}),f^{r_i+1}(\xi)]$ such that, using second
order Taylor expansion and the lower bound from
Lemma~\ref{le:bobo}
\begin{eqnarray}
\label{eq:umlado}
|f^{k-r_i}(x_{l_i})-x_{m(k-r_i-1)}|e^{-\tau (k-r_i)}
& \geq &
|f^{k}(\xi)-f^{k-r_i}(x_{l_i})|
\\
&=& 
|(f^{k-r_i-1})'(w)|
\cdot
|f^{r_i+1}(\xi)-f(x_{l_i})| \nonumber
\\
& \geq &  
C^{-1} |(f^{k-r_i-1})'(w)|\cdot
|f''(x_{l_i})| \cdot
|f^{r_i}(\xi)-x_{l_i}|^{2} \nonumber
\\
& \geq &
(AC)^{-1} 
|(f^{k-r_i-1})'(f^{r_i+1}(\xi))|\cdot
|f''(x_{l_i})|
|f^{r_i}(\xi)-x_{l_i}|^{2},\nonumber
\end{eqnarray}
where we have used Lemma~\ref{le:distorcao} in the last
inequality. Note that we can use the lower bound from
Lemma~\ref{le:bobo} since $r_i$ is an essential or inessential
return time for $\omega$.

The last two expression together show that
\begin{align*}
  |z_k-y_k|\cdot|f^{r_i}(\xi)-x_{l_i}| \le
  (AC^2)\cdot|f^{k-r_i}(x_{l_i})-x_{m(k-r_i-1)}| e^{-\tau
    (k-r_i)} \cdot |\omega_i|.
\end{align*}
This and Lemma~\ref{l3.1,5} provide
\begin{align*}
\left|
\frac{f^\prime(z_k)-f^\prime(y_k)}{f^\prime(y_k)}
\right| 
&\le
K_1\left| \frac{z_k-y_k}{y_k-\tilde x_{k}} \right|
\le
AC^2 K_1 e^{-\tau(k-r_i)}\cdot
\frac{|\omega_i|\cdot
|f^{k-r_i}(x_{l_i})-x_{m(k-r_i-1)}|}
{|f^{r_i}(\xi)-x_{l_i}|\cdot|y_k-\tilde x_{k}|}.
\end{align*}
In order to bound the denominator from below, we note that
clearly $|f^{r_i}(\xi)-x_{l_i}|\ge\Delta_{r_i}(\omega)$
since $\xi\in\omega$ and $x_{l_i}\in\cC$. Moreover, from
Proposition~\ref{pr:bindingaway}, the closest critical point
$\tilde x_k$ to $y_k$ and the closest critical point
$x_{m(k-r_i-1)}$ to $f^{k-r_i}(x_{l_i})$ are either equal or
else consecutive critical points of $f$ and, in the latter
case, both $\omega_k$ and $f^{k-r_i}(x_{l_i})$ lie between
these consecutive critical points. In the case $\tilde
x_k=x_{m(k-r_i-1)}$ we can bound the previous expression by
\begin{align*}
&D\cdot e^{-\tau(k-r_i)}\cdot
\frac{|\omega_i|}{\De_{r_i}(\omega)}
\cdot \frac{|f^{k-r_i}(x_{l_i})-x_{m(k-r_i-1)}|}
{|f^{k-r_i}(x_{l_i})-x_{m(k-r_i-1)}|-
|y_k- f^{k-r_i}(x_{l_i})|}
\\
&\le
D\cdot \frac{e^{-\tau(k-r_i)}}{1-e^{-\tau(k-r_i)}}
\cdot
\frac{|\omega_i|}{\De_{r_i}(\omega)}
\le
D_1\cdot e^{-\tau(k-r_i)}\cdot
\frac{|\omega_i|}{\De_{r_i}(\omega)}.
\end{align*}
But when $\tilde x_k\neq x_{m(k-r_i-1)}$ we have, since
$y_k$ and $f^{k-r_i}(x_{l_i})$ are between these critical points
\begin{align*}
  |\tilde x_k-y_k|
  &\ge
  |\tilde x_k-f^{k-r_i}(x_{l_i})| - |y_k-f^{k-r_i}(x_{l_i})|
  \\
  &\ge
  |f^{k-r_i}(x_{l_i}) - x_{m(k-r_i-1)}|
  - e^{-\tau(k-r_i)}|f^{k-r_i}(x_{l_i}) - x_{m(k-r_i-1)}|
  \\
  &=
  (1-e^{-\tau(k-r_i)})|f^{k-r_i}(x_{l_i}) - x_{m(k-r_i-1)}|
\end{align*}
and we arrive at the same bound as before.

On the other hand, for the second case in~\eqref{BCeq} we
get a similar inequality in~\eqref{eq:umlado} providing
\[
|z_k-y_k|\cdot|f^{r_i}(\xi)-x_{l_i}|
\le
(AC^2)\cdot \ep^{1+\tau}
e^{-\tau (k-r_i)}
\cdot |\omega_i|
\]
and thus by definition of $\tilde f$ we get
\begin{eqnarray*}
\left|
\frac{f^\prime(z_k)-f^\prime(y_k)}{f^\prime(y_k)}
\right| 
&\le&
 \frac{L\cdot |z_k-y_k|}{\tilde \sigma}
\le
\frac{ACL}{\tilde\sigma} e^{-\tau(k-r_i)}
\frac{|\omega_i|\cdot
  \ep^{1+\tau}}{|f^{r_i}(\xi)-x_{l_i}|}
\\
&\le&
D_2 \cdot e^{-\tau(k-r_i)} \cdot
\frac{|\omega_i|}{\De_{r_i}(\omega)}.
\end{eqnarray*}
This shows that for every $i=1,\dots,\gamma$ which is
\emph{not} an escape time, we have
\begin{equation}
  \label{eq:espirro}
\sum_{k=r_i}^{r_i+\ell} 
\left|
\frac{f^\prime(z_k)-f^\prime(y_k)}{f^\prime(y_k)}
\right| 
\le D_3 \cdot \frac{|\omega_i|}{\De_{r_i}(\omega)}
\le\frac1{C}\cdot\frac{|I(l_i,s_i,j_i)^+|}{|I(l_i,s_i)|},
\end{equation}
for all $\ell=1,\dots,p_i$, where we have used the
definition of $\omega_i$ and of host interval, together with
the same estimate as in~\eqref{eq:tosse}.

Considering
\eqref{eq:quotient1},~\eqref{eq:cha},~\eqref{eq:naoescapas!}
and~\eqref{eq:espirro} and summing over all iterates we
obtain
\begin{equation}
  \label{eq:somaF}
\sum_{k=0}^{n-1} 
\left|
\frac{f^\prime(z_k)-f^\prime(y_k)}{f^\prime(y_k)}
\right| 
\le D_4
\sum_{k\in F_1}\frac{|f^k(\omega)|}{\Delta_{k}(\omega)}
+
\Big(\frac{L}{\tilde\sigma}+K_2\Big)
\sum_{k\in F_2} |f^k(\omega)|.
\end{equation}
Here the left hand side sum is over the set $F_1$ of free
iterates together with return situations (essential and
inessential returns and escapes) from $k=0$ to $k=n-1$.  The
right hand side sum is over the set $F_2$ of free iterates
which are not followed by any return, from
$r_\gamma+p_\gamma$ to $n$.  Moreover $D_4$ is a constant
depending only on $\ep,\tau$ and $\tilde\sigma$.  So if we
can bound~\eqref{eq:somaF} uniformly we then find a uniform
bound to~\eqref{eq:Distortion} also and complete the proof
of the lemma.

The right hand side sum in~\eqref{eq:somaF} can is bounded $
\sum_{k\in F_2} |f^k(\omega)|\le \sum_{k=0}^{n-1}
\tilde\sigma^{k-n}|f^n(\omega)|\le C, $ since
$|f^n(\omega)|$ is always less than $1$.
Now we bound the left hand side sum
\begin{align*}
\sum_{k\in F_1}\frac{|f^k(\omega)|}{\Delta_{k}(\omega)}
& \le 
\sum_{\genfrac{}{}{0pt}{}{k\in F_1}{
  ( l_k, s_k)=(0,0)}}
\frac{|f^k(\omega)|}{\Delta_{k}(\omega)} 
+
\sum_{|l|\ge k_0} \sum_{|s|\ge 1} \sum_{
\genfrac{}{}{0pt}{}{k\in F_1}{
  ( l_k, s_k)=(l,s)}
}
\frac{|f^k(\omega)|}{\Delta_{k}(\omega)} 
\\
&\le
C\cdot \frac{\sigma_1}{\sigma_1-1}
+
\sum_{|l|\ge k_0} \sum_{|s|\ge 1}
\frac{\sigma_1}{\sigma_1-1}\cdot \frac1{C}\cdot
\frac{|f^{q(l,s)}(\omega)|}{|I(l,s)|},
\end{align*}
   by \eqref{eq:tosse}, where $q(l,s)=\max\{ 0\le k \le n-1
   : (\hat l_k, \hat s_k) =(l, s) \}$ and we convention that
   whenever $\{ 0\le k \le n-1 : (\hat l_q, \hat s_q) =(l,
   s) \}=\emptyset$ we have
   $\frac{|f^{q(l,s)}(\omega)|}{|I(l,s)|}=0$.  We have
   used the following estimate for any given fixed value of
   $(l,s)$
   \[
   \sum_{\{k:\hat s_k=s\}} |f^k(\omega)|
   \le
   |f^{q(l,s)}(\omega)| \sum_{\{k:(\hat l_k, \hat
     s_k)=(l,s)\}} \sigma_1^{k-q(l,s)}
   \le
   \frac{\sigma_1}{\sigma_1-1}\cdot |f^{q(l,s)}(\omega)|
   \le C \cdot |I(l,s,j)^+|,
   \]
   because writing $\{k:\hat s_k=s\}=\{k_1<k_2<\dots<k_h\}$
   we have $|f^{k_i}(\omega)|\le
   \sigma_1^{-1}\cdot|f^{k_{i+1}}(\omega)|$ for
   $i=1,\dots,h$, where
   $1<\sigma_1=\min\{\sigma_0,e^{(1-2\zeta)\cdot\frac{\pi}{\beta}(k_0+1)}\}
   \le\min\{\sigma_0,e^{(1-2\zeta)\cdot\frac{\pi}{\beta}(|l|+|s|)}\}$,
   after Lemma~\ref{le:maisdistorcao}(b) together with
   Remark~\ref{rmk:expansaolivre}.


   We observe that by construction we must have
   $|I(l,s,j)^+|/|I(l,s)|\le 9(|l|+|s|)^{-3}$
   and so we arrive at
   $
    \sum_{k\in
    F_1}\frac{|f^k(\omega)|}{\Delta_{k}(\omega)}
  \le C 
   \sum_{|l|\ge k_0} \sum_{|s|\ge 1 }
   \frac{9} {(|l|+|s|)^3}
   < \infty,
   $
   finishing the proof of the lemma.

\end{proof}



\section{Probability of deep essential returns}
 \label{sec:fundamental-lemma}


Here we use the results from
Section~\ref{sec:auxiliary-lemmas} to 
%
estimate the probability of having an orbit with a given
sequence of host intervals at essential return situations.

Recall that $\cC_\infty=\cup_{n=0}^{\infty}(f^n)^{-1}(\cC)$
is the set of pre-orbits of the critical set $\cC$.  For
each $x\in I\setminus\cC_\infty$ let $\omega$ be the element
of $\cP_n$ which contains $x$. 

Consider the sets $R_n(\omega)$ and $Q_n(\omega)$ of return
situations (essential and inessential returns and escapes)
and indexes of host intervals for $\omega$, during the
iterates $0$ to $n$.  Let $u_n(\omega)$ denote, for
$\omega\in\cP_n$, the number of essential returns or
escapes associated to $\omega$ between $0$ and $n$, let
$0\leq t_1(\omega)\leq \ldots \leq t_{u_n}(\omega)\leq n$ be the
instants of occurrence of the essential returns or escapes
and let $(l_1,s_1,j_1),\ldots, (l_{u_n},s_{u_n},j_{u_n})$ be
the corresponding critical points and indexes of the
respective host intervals.  We say that the sum
$|s_i|+|j_i|$ is the \emph{depth} of the corresponding
return of $\omega$.

Note that by construction $t_1(\omega)=0$ for all $\omega\in
\cP_0\setminus I(0,0,0)$ and $t_1(\omega)=1$ for $\omega=
I(0,0,0)$, see Remark~\ref{rmk:everybodyreturns}.

\begin{lemma}[No return probability]
\label{le:escape-ae}
For every $n\ge0$ there exists \emph{no} non-degenerate
interval $\omega\in\cP_n$ such that $\omega\in\cP_{n+k}$ for
all $k\ge1$. Moreover there exist constants $0<\xi_0<1$ and
$K_0>0$ (depending only on $\sigma,\sigma_0$ and on $\zeta$
from Lemma~\ref{le:maisdistorcao}), and
$n_0\ge1$ such that for every $n>n_0$ we have
$
\lambda\big(
\bigcup \{\omega\in\cP_n: u_n\mid\omega=1 \}
\big) \le K_0\cdot e^{-\xi_0 n}.
$
\end{lemma}

\begin{proof}
  If $\omega\in\cP_{n+k}$ for all $k\ge0$, then $\omega$ is
  not refined in all future iterates. This means that
  $f^{n+k}(\omega)$ has no essential returns nor escapes for
  $k\ge1$. Hence every iterate is either free or a binding
  time associated to a inessential return.  Let
  $p_0,p_1,p_2,\dots$ be the length of the binding period
  associated to every (if any) inessential return time
  $t<r_0<r_1<r_2<\dots$ for $\omega$ after $t$, where $0\le
  t\le n$ is the last essential return time or escape before
  $n$.  Let $k\ge0$ and $r_i+p_i\le n+k<r_{i+1}$ for some
  $i\ge0$, where we set $r_{i+1}=+\infty$ if $r_i$ is the
  last inessential return time after $t$ (it may happen that
  $r_1=+\infty$ in which case $i=0$ and there are no
  inessential returns after
  $t$). Lemmas~\ref{le:maisdistorcao}
  and~\ref{lem:fn-omega-estimate} ensure that
  \begin{equation}
    \label{eq:pequenomega1}
    2\ge|f^{n+k}(\omega)|
    \ge 
    2^{i}\cdot \sigma_0^{n+k-t-\sum_{j=0}^{i}(p_j+1)} |f^t(\omega)|,
  \end{equation}
  for arbitrarily big values of $k\ge0$, where the exponent
  of $\sigma_0$ counts the number of free iterates between
  the times $t$ and $n+k$.  

  Note that if there are no inessential returns, then $i=0$
  and so
  \begin{equation}
    \label{eq:pequenomega2}
    |f^t(\omega)|\le 2\cdot\sigma_0^{t+p_0+1-n-k}.
  \end{equation}
  Otherwise 
$|f^t(\omega)|\le
2^{1-i}\cdot\sigma_0^{t+\sum_{j=0}^i(p_j+1) -n-k}$ from
\eqref{eq:pequenomega1}.

We conclude that $|f^t(\omega)|=0$ which is not possible
for a non-degenerate interval. This proves the first part of
of the Lemma.

Now let $\omega\in\cP_n$ be such that $u_n(\omega)=1$.  Then
by the construction of the refinement, we see that
$\omega\in\cP_i$ for all $0\le i \le n$.  Hence
\begin{itemize}
\item either $\omega=I(0,0,0)$ with $t_1=1$ the unique
  essential return up to iterate $n$ and
  $f(\omega)=I(l,s,j)$ with $(l,s,j)\neq(0,0,0)$;
\item or $\omega=I(l,s,j)$ with $|l|\ge k_0,|s|\ge1$ and
  $j=1,\dots,(|l|+|s|)^3$, having a single essential return
  $t_1=0$ up to iterate $n$.
\end{itemize}
We concentrate on the latter case and write
$p_0,p_1,p_2,\dots$ and $0=t_1=r_0<r_1<r_2<\dots$ the
binding periods associated to their respective inessential
return times of the orbit of $\omega$ as before. Then
by~\eqref{eq:pequenomega1} and~\eqref{eq:pequenomega2}, and
since during binding periods the intervals of the partition
are not subdivided (see
Section~\ref{sec:concstr-part-bound-distortion}), we have
either $n\le r_0+p_0=p_0$ or
\begin{align*}
|\omega|\le 
2\cdot \sigma_0^{p_0-n} 
&
\quad\mbox{if  } p_0(\omega)< n\le r_1(\omega)+p_1(\omega) \le
+\infty, \quad\mbox{or} 
\\
|\omega|\le 2^{1-i}\cdot\sigma_0^{-(n-\sum_{j=0}^{i} (p_j(\omega)+1))} 
&
\quad\mbox{for }i\ge1
  \mbox{  such that  } r_i(\omega)+p_i(\omega) < n 
  \le r_{i+1}(\omega)+p_{i+1}(\omega).
\end{align*}
We note that in the case of $f(\omega)=I(l,s,j)$, that is,
when $t_1(\omega)=1$, we can repeat the arguments for the
interval $f(\omega)$, arriving at the same bounds for
$|\omega|$ except for an extra factor of $\tilde\sigma$
since $|f(\omega)|\ge\tilde\sigma|\omega|$.  Hence we may
write according to three cases above
\begin{eqnarray*}
\lambda\big(
\cup\{ \omega\in\cP_n: u_n\mid\omega=1 \}\big)
&\le& 
\sum_{n\le p_0(\omega)}
|\omega|
\quad + \hspace{-0.5cm}
\sum_{
{p_0(\omega)< n\le r_1(\omega)+p_1(\omega)\le+\infty}}
\hspace{-0.8cm}|\omega| 
\\
& & +
\sum_{
  {r_i(\omega)+p_i(\omega)< n \le r_{i+1}(\omega)+p_{i+1}(\omega)}
} \hspace{-1.7cm}|\omega| 
\quad+\quad
S_3
\\
& = & S_0 + S_1 + S_2+ S_3,
\end{eqnarray*}
where, by the above comment, we may assume that every sum
ranges over $\omega\in\cP_0\setminus I(0,0,0)$ and $S_3$
corresponds to the sum over the partition elements in
$I(0,0,0)\cap\cP_1$, which is bounded by
$(S_0+S_1+S_2)/\tilde\sigma$. 

For $S_0$ we use Lemma~\ref{le:maisdistorcao}(a) to deduce
that the summands in $S_0$ are the elements of $\cP_0$ such
that
$n\le p_0(\omega)\le \frac{2\pi}{\beta\log\sigma}(|l|+|s|)$,
that is $|l|+|s|\ge\frac{\beta\log\sigma}{2\pi}\cdot n$,
thus by Remark~\ref{rmk:partitiondimension}, setting
$C_0=\frac{\beta\log\sigma}{2\pi}$
\begin{align}
S_0
&\le&\hspace{-0.5cm}
\sum_{|l|+|s|\ge C_0\cdot  n} \hspace{-0.5cm}
a_{1}\frac{e^{-(\pi/\beta)(|l|+|s|)}}{(|l|+|s|)^3}
\le
a_1 \sum_{k\ge C_0\cdot n} k\cdot
\frac{e^{-(\pi/\beta)k}}{k^3} 
\le
\frac{a_1}{(C_0\cdot n)^2}\cdot
\frac{e^{-(\pi/\beta) C_0\cdot n}}{( 1-e^{-C_0\cdot n})}
\le K_0' \cdot \sigma^{-n/2}.\label{eq:serie}
\end{align}
We write $S_1=S_{11}+S_{12}$ where
\[
S_{11}= \sum_{
  {p_0(\omega)\le
    n/2}
} |\omega|
\quad\mbox{and}\quad S_{12}=\sum_{n/2<p_0(\omega)}
|\omega|\le K_0'\cdot \sigma^{-n/4}
\]
and we have used the bound~\eqref{eq:serie}.  We also split
$S_2=S_{21}+S_{22}$ according to whether
$n-\sum_{k=0}^{i}(p_k(\omega)+1)\ge n/2$ or not, obtaining
\[
S_{21}= \sum_{|\omega|\le 2^{1-i}\cdot
\sigma_0^{-n/2}} |\omega|
\quad\mbox{and}\quad
S_{22}=\sum_{
{n-\sum_{j=0}^{i}(p_j(\omega)+1) < n/2}} 
|\omega|.
\]
Since $2^{1-i}\le2$ we get $S_{11}+S_{21}\le 2\cdot S_{11}$
and the summands $\omega\in\cP_0$ satisfy $|\omega| \le
2\cdot\sigma_0^{-n/2}$, thus by
Remark~\ref{rmk:partitiondimension} we get
\begin{equation*}
|l|+|s|\ge \frac{n}2\cdot\frac{\log\sigma_0}{3+\pi/\beta}
+ \frac{-\log(2/a_1)}{3+\pi/\beta}
\ge C_1 \cdot \frac{n}4,
\end{equation*}
where $C_1=\log\sigma_0/(3+\pi/\beta)$ for every big enough
$n$. Then $S_{11}+S_{21}\le 2 K_0'' \cdot \sigma_0^{-n/5}$ by
the same calculations as in~\eqref{eq:serie} with slightly
different constants.

For $S_{22}$ we note that $n-\sum_{k=0}^{i}(1+p_k(\omega)) <
n/2$ implies $\sum_{k=0}^{i}(1+p_k(\omega))>n/2$ and so, by
Lemma~\ref{le:maisdistorcao}(b) we get
\[2\ge|f^{r_i+p_i+1}(\omega)| >
e^{(1-2\zeta)\frac{\pi}{\beta}\sum_{j=0}^{i}(1+p_j(\omega))}
\cdot |\omega| > e^{n\pi(1-2\zeta)/(2\beta)}\cdot |\omega|
\]
and hence again by Remark~\ref{rmk:partitiondimension}, for
every big enough $n$, we must have $|l|+|s|>C_2\cdot n/4$
where $C_2=\pi(1-2\zeta)/(3\beta+\pi)$.  We deduce that
$S_{22}\le K_0''' \cdot e^{-(\pi/\beta) C_2\cdot n/4}$
following the same calculations as in~\eqref{eq:serie}.

Putting all together we see that there are constants
$0<\xi_0<1$, $K_0>0$ and $n_0\ge1$ such that
$S_0+S_1+S_2+S_3\le K_0\cdot e^{-\xi_0 n}$ for all $n\ge
n_0$, with $\xi_0$ and $K_0$ dependent on $\sigma,\sigma_0$
and $\zeta$, as stated.
\end{proof}

Let $v\le u\leq n$ and $v$ pairs of positive integers
$(\eta_1,\upsilon_1),\ldots,(\eta_v, \upsilon_v)$ be given,
$\eta_i\ge k_0$, $\upsilon_i\ge1$ and
$\eta_i+\upsilon_i\ge \Theta$, where
\begin{equation}
  \label{eq:Theta}
\Theta=\Theta(\flat)=\frac{\beta}{\pi}\cdot \log\frac{a_2}{\flat}.
\end{equation}
This value of $\Theta$ is chosen in such a way that
$\dist(I(\eta,\upsilon),\cC)<\flat$ if, and only if,
$\eta+\upsilon\ge\Theta$.  We assume that $\flat$ is so
small that $\Theta\ge\theta$ (recall that $\theta$ was
defined in Lemma~\ref{le:maisdistorcao}) and, following
Remark~\ref{rmk:Qmaiorque1}, that $\Theta$ is big enough in
order that $Q(l,s,\tau,\rho)>1$ for all $|l|+|s|\ge\Theta$.

For $\omega\in\cP_n$ with $u_n(\omega)=u$, let
$0=t_1<\dots<t_u\le n$ be the return situations (essential
returns or escapes) of $\omega$, $(l_i,s_i)$ the
corresponding indexes of the host intervals and
$d_n(\omega)=v$ the number of pairs $(l_i,s_i)$ such
that $|l_i|+|s_i|\ge\Theta$.  We say that a return situation
whose depth satisfies $|l_i|+|s_i|\ge\Theta$ is a \emph{deep
  return}.  Denote by $1\le r_1 < \dots < r_v \le u$ the
indexes of the return situations corresponding to deep
returns.

We now define the subset
$A_{(\eta_1,\upsilon_1),\ldots,(\eta_v,
  \upsilon_v)}^{u,v}(n)$ as
\[
\bigcup\big\{\omega\in\cP_n : \, u_n(\omega)=u, \, d_n(\omega)=v
\,\mbox{ and }\,
 \big| f^{t_{r_i}}(\omega)\big| \subset I(\eta_i,\upsilon_i),
\,
 i=1,\ldots,v\big\}.
\]

\begin{proposition}[Probability of essential returns with
  specified depths]
\label{prop:depth-probability-1} 
If $\Theta$ is large enough (depending on $D_0$ from
Lemma~\ref{lem:bounded-distortion}), then for every big
$n\ge u\ge v \ge 1$ and $\vartheta=\max\{\alpha,3\zeta\}$
\[
\lambda\left(
A_{(\eta_1,\upsilon_1),\ldots,(\eta_v, \upsilon_v)}^{u,v}(n)
\right)
\leq
\binom{u}{v}
\exp\Big[
( 2\vartheta -1) \frac{\pi}{\beta}
\sum_{i=1}^v (\eta_i+\upsilon_i)
\Big].
\]
\end{proposition}

\begin{proof}
  We start by fixing $n\in\bN$, $u\in\{1,\dots,n\}$,
  $v\in\{1,\dots,u\}$ and taking $\omega^0\in\cP_0$.  Let
  $\omega\in\cP_n$ be such that $\omega\subset\omega^0$ and
  $u_n(\omega) =u$ and $1=t_1<\dots<t_u\le n$ be the
  return situations (essential returns or escapes) of
  $\omega$.
  
  For $m=1,\dots,u$ we write
  $\omega^m=\omega((l_1,s_1,j_1),\dots,(l_m,s_m,j_m))\in\cP_{t_m}$
  the subset of $\omega^00$ satisfying
  \begin{align}\label{eq:nested}
    f^{t_i}(\omega^m)\subset I(l_i,s_i,j_i)^+\,,\,
    i\in\{1,\ldots,m-1\} \,\mbox{ and }\,
    I(l_m,s_m,j_m)\subset f^{t_m}(\omega^m)\subset
    I(l_m,s_m,j_m)^+,
  \end{align}
by the definition of the sequence of partitions $\cP_n$.
Note that we get a nested sequence of sets
\[
\omega^0 \supsetneq \omega^1 \supsetneq \dots \supsetneq \omega^u=\omega.
\]
We define $\cT=\{ \omega\in\cP_n : \omega\subset\omega^0, \,
u_n(\omega) = u \}$ and consider the sequence of deep
return situations (essential or escapes) of each
$\omega\in\cT$: $1\le r_1 < \dots < r_v \le u$, that is, the
indexes of the return situations such that for $i=1,\dots,v$
\[
 t_{r_i} \mbox{  satisfies  } 
\big| f^{t_{r_i}}(x) \big|
 \in I(\eta_i,\upsilon_i) 
\mbox{ for all } 
x\in\omega \mbox{ and  } \eta_i+\upsilon_i \ge \Theta.
\]
Now we define by induction a sequence of partitions of $\cT$
which will enable us to determine the estimates we need.

Start by putting $\cV_0=\cup \{ \omega\in\cT\}$. We define
for $1\le i\le v$ and $h\ge0$ such that $r_i\le
r_i+h<r_{i+1}$ the subset
\begin{align*}
  \cV^{j^1,\dots,j^i}_{r_i+h}
  =\bigcup
  \{ \omega^{r_i+h}: \omega\in\cT \,\&\,
  \eqref{eq:nested} \text{  holds with  }
  m=r_i+h \,\&\, j^k=j_{r_k} \text{  for  } k=1,\dots,i\},
\end{align*}
where we make the convention $r_{v+1}=u$;
and for $1\le h <r_1$ we set
$\cV_{h}
  =\bigcup
  \{
  \omega^{h}: \omega\in\cT
  \}.$



Now we compare
$\lambda\big(\cV_{r_{i}}^{j^1,\dots,j^{i}}\big)$ with
$\lambda\big( \cV_{r_{i-1}}^{j^1,\dots,j^{i-1}}\big)$. We claim
\begin{align}
  \label{eq:compareV}
  \lambda\big(\cV_{r_i}^{j^1,\dots,j^i}\big)
  &\le
\frac{C}{(\eta_i+\upsilon_i)^3}\cdot
\frac{e^{-\frac{\pi}{\beta}(\eta_i+\upsilon_i)}}
{e^{-\vartheta\frac{\pi}{\beta}(\eta_{i-1}+\upsilon_{i-1})}}
\cdot \lambda\big( \cV_{r_{i-1}}^{j^1,\dots,j^{i-1}} \big),
\end{align}
where $\vartheta=\max\{2\zeta,\alpha\}\in (0,1)$.
Assuming this claim we deduce
\begin{align*}
  \lambda\big( \cV_u^{j^1,\dots,j^v} \big)
&\le
\dots
\le
\lambda\big(\cV_v^{j^1,\dots,j^v} \big)
\le
\frac{C}{(\eta_v+\upsilon_v)^3}\cdot
\frac{e^{-\frac{\pi}{\beta}(\eta_v+\upsilon_v)}}
{e^{-\vartheta\frac{\pi}{\beta}(\eta_{v-1}+\upsilon_{v-1})}}
\cdot \lambda\big( \cV_{r_v-1}^{j^1,\dots,j^{v-1}} \big)
\\
&\le
\left(
\prod_{i=1}^s \frac{C}{(\eta_i+\upsilon_i)^3}
\right)
\exp\left(
-\frac{\pi}{\beta}\sum_{i=1}^v (\eta_i+\upsilon_i)
+\vartheta\frac{\pi}{\beta}\sum_{i=1}^v (\eta_{i-1}+\upsilon_{i-1})
\right)
\cdot
\lambda(\cV_0).
\end{align*}
Finally we need to consider all possible combinations of the
events $\cV_u^{j^1,\dots,j^v}$ which are included in
$A_{(\eta_1,\upsilon_1),\ldots,(\eta_v,
  \upsilon_v)}^{u,v}(n)$. Note that for any given $v\le u$
there are $\binom{u}{v}$ ways of having $v$ deep returns
among $u$ return situations and, by symmetry, for any
sequence of deep returns with given depth
$(\eta_i,\upsilon_i)$ there are
$4\cdot(\eta_i+\upsilon_i)^3$ different possibilities of
falling in an element of the partition $\cP_0$. Thus since
$\lambda(\cV_0)\le\lambda(\omega_0)$ we arrive at
\begin{align*}
  \lambda\Big( A_{(\eta_1,\upsilon_1),\ldots,(\eta_v,
  \upsilon_v)}^{u,v}(n)\Big)
&\le
\binom{u}{v}\cdot \prod_{i=1}^v 4(\eta_i+\upsilon_i)^3
\cdot
\left(
\prod_{i=1}^v \frac{C}{(\eta_i+\upsilon_i)^3}
\right)
\cdot
\\
&
\cdot
\exp\left(
(\vartheta-1)\frac{\pi}{\beta}\sum_{i=1}^v (\eta_i+\upsilon_i)
\right)
\cdot
\sum_{\omega_0\in\cP_0}
  e^{\vartheta\frac{\pi}{\beta}(\eta_0+\upsilon_0)}\lambda(\omega_0)
\\
&\le
\binom{u}{v}
\cdot\exp\left(
(2\vartheta-1)\frac{\pi}{\beta}\sum_{i=1}^v (\eta_i+\upsilon_i)
\right),
\end{align*}
where we have used that $\sum_{\omega_0\in\cP_0}
e^{\vartheta\frac{\pi}{\beta}(\eta_0+\upsilon_0)}\lambda(\omega_0)<\infty$
and also that $\sum_{i=1}^v (\eta_i+\upsilon_i)\ge v\Theta$
and that $\Theta$ can be taken as large as needed. 

Thus to complete proof it is enough to prove the
claim~\eqref{eq:compareV}.  For this we proceed as follows.

Given $\omega\in\cT$ and $\omega^{r_{i-1}}\in
\cV_{r_{i-1}}^{j^1,\dots,j^{i-1}}$ we have
$\omega^{r_i}=\omega^{r_{i-1}}\cap
\cV_{r_{i}}^{j^1,\dots,j^{i}}$ which is the set of points in
$\omega^{r_{i-1}}$ which remain in the next level of the
partition. We divide the argument into the following cases.

\begin{enumerate}
\item $t_{r_{i}-1}$ is an essential return with depth $(l,s)$.
\end{enumerate}

In this case, since $\omega^{r_i}\subsetneq\omega^{r_i-1}
\subsetneq\cdots\subsetneq \omega^{r_{i-1}+1}\subsetneq
\omega^{r_{i-1}}$ and
$\omega^{r_{i}}\in\cP_{t_{r_{i}}-1}$ by the refinement
algorithm, we can use the Bounded Distortion
Lemma~\ref{lem:bounded-distortion} to write
\begin{align*}
  \frac{\big| \omega^{r_i}\big|}{\big|\omega^{r_{i-1}}\big|}
  &\le
  \frac{|\omega^{r_{i-1}+1}|}{|\omega^{r_{i-1}}|}
  \cdots
  \frac{|\omega^{r_i}|}{|\omega^{r_i-1}|}
  \le
  1\cdot 1\cdots
  D_0\cdot \frac{\Big|
    f^{t_{r_{i}}}\big(\omega^{r_{i}}\big) \Big|}
  {\Big| f^{t_{r_{i}}}\big(\omega^{r_{i}-1}\big)  \Big|}
  \\
  &\le
  D_0\cdot \frac{\big| I(\eta_i,\upsilon_i,j_i)^+ \big|}
  {a_1 Q e^{-2\zeta\frac{\pi}{\beta}(|l|+|s|)}}
  \le
  D_0\cdot\frac{9 a_1 e^{-\frac{\pi}{\beta}(\eta_i+\upsilon_i)}}
  {a_1 Q \cdot (\eta_i+\upsilon_i)^3 e^{-2\zeta\frac{\pi}{\beta}(|l|+|s|)} }
  \\
  &=
  \frac{C}{(\eta_i+\upsilon_i)^3}\cdot
  \frac{e^{-\frac{\pi}{\beta}(\eta_i+\upsilon_i)}}
  {e^{-2\zeta\frac{\pi}{\beta}(|l|+|s|)}}
  \le
  \frac{C}{(\eta_i+\upsilon_i)^3}\cdot
  \frac{e^{-\frac{\pi}{\beta}(\eta_i+\upsilon_i)}}
  {e^{-2\zeta\frac{\pi}{\beta}(\eta_{i-1}+\upsilon_{i-1})}},
\end{align*}
where $C=C(\rho,\tau,\epsilon,\sigma)$ and in the second and
third inequalities we used
Remark~\ref{rmk:partitiondimension} and
Lemma~\ref{lem:fn-omega-estimate}(2). In the last inequality
we argue as follows.

If $r_i=r_{i-1}+1$, then $|l|+|s|=\eta_{i-1}+\upsilon_{i-1}$
by definition. If $r_i>r_{i-1}+1$, then by definition of
deep returns $|l|+|s|<\Theta\le \eta_{i-1}+\upsilon_{i-1}$.

\begin{enumerate}
\item[(2)] $t_{r_{1}-1}$ is an escape time having host interval
  $I(l,s,j)^+$ with $(l,s,j)\neq (\pm k_0,1,1)$.
\end{enumerate}

Arguing as in the previous case, we only need to get a lower
bound for $| f^{t_{r_i}}\big(\omega^{r_{i}-1}\big)|$. By the
Mean Value Theorem, for some $\xi\in
f^{t_{r_{i-1}}}(\omega^{r_{i}-1})$ we have
\begin{align*}
  \Big| f^{t_{r_i}}\big(\omega^{r_{i}-1}\big) \Big| & =
  \Big| \big(f^{t_{r_i}-t_{r_{i}-1}}\big)'(\xi)\Big| \cdot
  \big|f^{t_{r_{i}-1}}(\omega^{r_{i}-1})\big|
  \\
  &= \Big|
  \big(f^{t_{r_i}-t_{r_{i}-1}-1}\big)'\big(f(\xi)\big)\Big|
  \cdot \Big| f'(\xi) \Big| \cdot
  \big|f^{t_{r_{i}-1}}(\omega^{r_{i}-1})\big|
  \\
  &\ge \big( \sigma_0^{t_{r_i}-t_{r_{i}-1}-1}\big)\cdot
  \frac{\tau\hat{x}^{\alpha-1}}{C}
  e^{(1-\alpha)\frac{\pi}{\beta}|l|} \cdot
  \big|f^{t_{r_{i}-1}}(\omega^{r_{i}-1})\big|
  \\
  &\ge C(\tau) \cdot e^{(1-\alpha)\frac{\pi}{\beta}|l|}\cdot
  a_1\cdot
  \frac{e^{-\frac{\pi}{\beta}(|l|+|s|)}}{(|l|+|s|)^3}
  \\
  &= C(\tau) \cdot
  \frac{e^{-\frac{\pi}{\beta}(\alpha|l|+|s|)}}{(|l|+|s|)^3}
  =
  C(\tau) \cdot
  e^{(\alpha-1)\frac\pi{\beta}|s|}\cdot
  \frac{e^{-\alpha\frac{\pi}{\beta}(|l|+|s|)}}{(|l|+|s|)^3}
  \\
  &\ge C(\tau)\cdot
  e^{(\alpha-1)\frac\pi{\beta}s(\tau)}\cdot
  \frac{e^{-\alpha\frac{\pi}{\beta}(|l|+|s|)}}{(|l|+|s|)^3}
  \ge 
  C(\tau)\cdot e^{-\alpha\frac{\pi}{\beta}(|l|+|s|)},
\end{align*}
where we have used the bound~\eqref{eq:expansaolivre} in the
first inequality and that $|s|\le s(\tau)$ in the last
inequalities.
Thus we arrive at a similar bound
\begin{equation}
  \label{eq:quotient}
   \frac{\big| \omega^{r_i} \big|}{\big|\omega^{r_{i-1}}\big|}
   \le
   \frac{C}{(\eta_i+\upsilon_i)^3}\cdot
   \frac{e^{-\frac{\pi}{\beta}(\eta_i+\upsilon_i)}}
   {e^{-\alpha\frac{\pi}{\beta}(\eta_{i-1}+\upsilon_{i-1})}}.
\end{equation}

\begin{enumerate}
\item[(3)] $t_{r_{i}-1}$ is an escape time having as host
  interval $I(k_0,1,1)$ or $I(-k_0,1,1)$.
\end{enumerate}

From Lemma~\ref{le:emezero} we have that points in
$f^{t_{r_{i}-1}}(\omega^{r_{i}-1})$ remain outside
$[-\epsilon,\epsilon]$ during a minimum number
$m_0\ge\log(1/\epsilon)$ of iterates with
derivative bigger than $\sigma_0$, which can be taken
larger than $e$. Hence
\begin{align*}
  \Big| f^{t_{r_i}}\big(\omega^{r_{i}-1}\big) \Big| 
  & =
  \Big|\big(f^{t_{r_i}-t_{r_{i}-1}}\big)'(\xi)\Big|
  \cdot \big|f^{t_{r_{i}-1}}(\omega^{r_{i}-1})\big|
  \\
  &\ge
  \sigma_0^{m_0}\cdot\big|f^{t_{r_{i}-1}}(\omega^{r_{i}-1})\big|
  \ge
  e^{\log(1/\epsilon)}\cdot\frac{e^{-\frac{\pi}{\beta}(k_0+1)}}{(k_0+1)^3}
  \\
  &\ge
  C\frac{e^{\log(1/\epsilon)}}{\big(\log(1/\epsilon)\big)^3}
  \cdot e^{-\frac{\pi}{\beta}(k_0+1)}
  \ge
  C\cdot e^{-\alpha\frac{\pi}{\beta}(k_0+1)},
\end{align*}
where we have used that $k_0\le C\cdot\log(1/\epsilon)$ and
assumed that $\epsilon$ is small enough so that
$e^{\log(1/\epsilon)}\ge
\big(\log(1/\epsilon)\big)^3$. Thus we arrive again at a
bound of the form \eqref{eq:quotient}. We remark that since 
$| f^{t_{r_i}}\big(\omega^{r_{i}-1}\big)|\le2$ the number of
free iterates $t_{r_i}-t_{r_{i}-1}$ is bounded from above by
a constant depending on $k_0$, $\beta$ and $\sigma_0$ only,
so that the lower bound we get is uniformly away from zero
for all intervals satisfying this third case.

Now we are ready to obtain \eqref{eq:compareV} as follows
for $1\le i\le v$
\[
\lambda\big(\cV_{r_i}^{j^1,\dots,j^i}\big)
=
\sum_{\omega^{r_{i-1}}\in\cV_{r_{i-1}}^{j^1,\dots,j^{i-1}}}
\frac{| \omega^{r_i} |}{|\omega^{r_{i-1}}|}\cdot
|\omega^{r_{i-1}}|
\le
\frac{C}{(\eta_i+\upsilon_i)^3}\cdot
\frac{e^{-\frac{\pi}{\beta}(\eta_i+\upsilon_i)}}
{e^{-\vartheta\frac{\pi}{\beta}(\eta_{i-1}+\upsilon_{i-1})}}
\cdot \lambda\big( \cV_{r_{i-1}}^{j^1,\dots,j^{i-1}} \big),
\]
where $\vartheta=\max\{2\zeta,\alpha\}\in (0,1)$ is obtained
comparing the bounds for each case.  

The proof is complete.
\end{proof}

We define using the same notations as before
with $\eta\ge k_0$, $\varsigma\ge1$ and
$\eta+\varsigma\ge\Theta$ 
\begin{align*}
A_{(\eta,\varsigma),j}^{v,u} (n) 
&= \bigcup\{ \omega\in \cP_n : (u_n(\omega),d_n(\omega))=(u,v)
\quad\text{and}\quad
|f^{t_{r_j}}(x)| \in I(\eta,\varsigma),\forall x\in\omega \},
\\
A_{(\eta,\varsigma)}^{v,u} (n) 
&=
\bigcup\{ \omega\in\cP_n : (u_n(\omega),d_n(\omega))=(u,v)
\quad\text{and}\quad
\\
&\qquad\qquad
\text{there exists  } 1\le j\le v \quad\text{such that}\quad
|f^{t_{r_j}}(x)| \in I(\eta,\varsigma) ,\forall x\in\omega \},
\\
A_{(\eta,\varsigma)}(n)
&=
\bigcup\{ \omega\in\cP_n : \text{there exists  } 
t\le n \text{ such that  $t(\omega)$ is an
  essential return and }
\\
&\qquad\qquad
|f^{t}(x)| \in I(\eta,\varsigma) \text{ for each } x\in\omega \},
\end{align*}
and derive the following corollary which will be used
during the final arguments.

\begin{corollary}
\label{cor:c.1}
We have for $1\le j\le v \le u \le n$ that 
\begin{enumerate}
\item $\lambda \big( A_{(\eta,\varsigma),j}^{v,u}(n) \big) \le
\binom{u}{v}e^{(2\vartheta-1)\frac{\pi}{\beta}(\eta+\varsigma)}$;
\item $\lambda \big(A_{(\eta,\varsigma)}^{v,u}(n)\big) \le
v \binom{u}{v} e^{(2\vartheta-1)\frac{\pi}{\beta}(\eta+\varsigma)}$;
\item $\lambda \big(A_{(\eta,\varsigma)}(n)\big) \le
n^3 e^{o(\Theta)n} e^{(4\vartheta-1)\frac{\pi}{\beta}(\eta+\varsigma)}$;
\end{enumerate}
if $\Theta$ is sufficiently big and $\zeta>\vartheta/3$ is
small enough, where $o(\Theta)\to0$ when $\Theta\to\infty$.
\end{corollary}

\begin{proof}
  For item (1) we note that since
\[
A_{(\eta,\varsigma),j}^{v,u} (n)
\subseteq 
\bigcup_{
\eta_i+\varsigma_i\ge\Theta, \eta_i\ge k_0, \varsigma_i\ge1, i\neq j}
A_{(\eta_1,\varsigma_1),\ldots,(\eta_{i-1}, \varsigma_{i-1}),
  (\eta,\varsigma),(\eta_{i+1}, \varsigma_{i+1}),\dots,
  (\eta_u,\varsigma_u)}^{v,u}(n)
\]
then 
$\lambda\big( A_{(\eta,\varsigma),j}^{v,u} (n) \big)
\le
\binom{u}{v}
\Big( \sum_{ l+s\ge\Theta}
e^{(2\vartheta-1)\frac{\pi}{\beta} (l+s)} \Big)^{v-1}
\cdot
e^{(2\vartheta-1)\frac{\pi}{\beta} (\eta+\varsigma)}
\le 
\binom{u}{v}
e^{(2\vartheta-1)\frac{\pi}{\beta} (\eta+\varsigma)},$ 
as long as $\Theta$ is big
enough in order that $\sum\limits_{l+s\ge\Theta}
e^{(2\vartheta-1)\frac{\pi}{\beta} (l+s)}\le1$.

From this we get item (2) since
$
A_{(\eta,\varsigma)}^{v,u}(n)\subset
\bigcup_{j=1}^v A_{(\eta,\varsigma),j}^{v,u}(n).
$

For item (3) we note that $ A_{(\eta,\varsigma)}(n)\subset
\bigcup_{v=1}^u \bigcup_{u=v}^n
A_{(\eta,\varsigma)}^{v,u}(n), $ but since there is a deep
essential return at iterate $t$ (which is \emph{not} an
escape situation) before $n$, we know that the corresponding
binding period $p$ is larger than
$\iota(M)(\eta+\varsigma)>\iota(M)\cdot\Theta$.  

Let us assume first that $t+p\le n$. In this situation the
maximum number $u$ of essential return situations in the
first $n$ iterates of such $\omega\in\cP_n$ is bounded by
$n/(\iota(M)\cdot\Theta\cdot \tilde v)$, where $\tilde v$ is
the number of deep essential returns among the $u$ essential
return situations. Since we know that $\tilde v\ge1$ we get
$u\le (\iota(M)\Theta)^{-1}\cdot m$. This alone enables us
to bound the measure of the subset
$A_{(\eta,\varsigma)}^*(n)$ of $A_{(\eta,\varsigma)}(n)$ of
those intervals $\omega\in\cP_n$ such that we can find
$t$ with $t+p\le n$ as follows:
\begin{align}
\lambda\big( A_{(\eta,\varsigma)}^* (n) \big)
&\le
\sum_{v=1}^u \sum_{u=v}^n  \lambda
\big(A_{(\eta,\varsigma)}^{v,u}(n)\big)\nonumber
\\
&\le
\sum_{v=1}^{n/(\iota(M)\Theta)}
\sum_{u=v}^n
v \binom{u}{v} e^{(2\vartheta-1)\frac{\pi}{\beta}(\eta+\varsigma)}
=
e^{(2\vartheta-1)\frac{\pi}{\beta}(\eta+\varsigma)}
\hspace{-0.4cm}
\sum_{v=1}^{(\iota(M)\Theta)^{-1}\cdot n} \hspace{-0.4cm} v\cdot
\sum_{u=v}^n \binom{u}{v}\nonumber
\\
&<
e^{(2\vartheta-1)\frac{\pi}{\beta}(\eta+\varsigma)}
\cdot n \cdot \hspace{-0.4cm}
\sum_{v=1}^{(\iota(M)\Theta)^{-1}\cdot n} 
\hspace{-0.4cm} v\binom{n}{v} \nonumber
\\
&<
n\cdot e^{(2\vartheta-1)\frac{\pi}{\beta}(\eta+\varsigma)}
\cdot \binom{n}{(\iota(M)\Theta)^{-1}\cdot n} \cdot
\sum_{v=1}^{(\iota(M)\Theta)^{-1}\cdot n} \hspace{-0.4cm} v \nonumber
\\
&<
\frac{n^3}{(\iota(M)\Theta)^2} 
\cdot
e^{(2\vartheta-1)\frac{\pi}{\beta}(\eta+\varsigma)}
\cdot
e^{o(\Theta) n}
\le
n^3 e^{o(\Theta)n}
e^{(2\vartheta-1)\frac{\pi}{\beta}(\eta+\varsigma)} \label{eq:part1}
\end{align}
where we used the bound
$\binom{n}{(\iota(M)\Theta)^{-1}\cdot n}\le e^{o(\Theta) n}$ which
  can be obtained by a straightforward application of
  Stirling's Formula, as long as $\Theta$ is big
  enough.

Let us now assume that the only deep essencial return $t$ of
$\omega$ satisfies $m=t+p> n$. By definition
$A_{(\eta,\varsigma)}(n)\setminus A_{(\eta,\varsigma)}^*(n) 
\subseteq A_{(\eta,\varsigma)}^*(m)$ and so
\begin{align}\label{eq:part2}
  \lambda\big( A_{(\eta,\varsigma)} (n) \big)
&\le
\lambda\big( A_{(\eta,\varsigma)}^* (n) \big) + 
\lambda\big( A_{(\eta,\varsigma)}^* (m) \big).
\end{align}
Hence we can apply the previous reasoning with $m$ in the
place of $n$ to obtain
\begin{align*}
\lambda\big( A_{(\eta,\varsigma)}^* (m) \big)
&\le
m^3 e^{o(\Theta)m} e^{(2\vartheta-1)\frac{\pi}{\beta}(\eta+\varsigma)}.
\end{align*}
Finally since
$m-n<p\le\frac{2\pi}{\beta\log\sigma}(\eta+\varsigma)$ (by
Lemma~\ref{le:maisdistorcao}(a)) we have for $\vartheta$
small enough
  \begin{align}
    \lambda\big( A_{(\eta,\varsigma)}^*(m) \big) &\le n^3
    e^{o(\Theta) n}
    e^{(2\vartheta-1)\frac{\pi}{\beta}(\eta+\varsigma)} \cdot
    e^{o(\Theta)(m-n)}\cdot \left(\frac{m}{n}\right)^3 \nonumber
    \\
    &\le n^3 e^{o(\Theta) n}
    \cdot e^{(\frac{2o(\Theta)}{\log\sigma}+2\vartheta-1)
      \frac{\pi}{\beta}(\eta+\varsigma)}
    \cdot \left( 1+
      \frac{2\pi(\eta+\varsigma)}{n\beta\log\sigma} \right)^3
    \nonumber
    \\
    &\le \label{eq:part3}
    n^3 e^{o(\Theta) n}
    e^{(3\vartheta-1)\frac{\pi}{\beta}(\eta+\varsigma)}.
  \end{align}
  Putting \eqref{eq:part1}, \eqref{eq:part2} and
  \eqref{eq:part3} together, we complete the proof of the
  lemma.
\end{proof}

\section{Slow recurrence to the critical set}
\label{sec:slow-recurrence}

Now we make use of the lemmas from
Section~\ref{sec:fundamental-lemma} to prove
Theorem~\ref{thm:approxestimate} and consequently also
Theorem~\ref{thm:existacim}.  We start by recalling the
definition of $\cC_n^\flat(\omega)$
from~\eqref{eq:averagedist} and that $u_n(\omega)$ is the
number of essential return situations or escape times of the
$f$-orbit of $\omega\in\cP_n$ between $0$ and $n$.

We let $0\leq t_1 <\ldots < t_{u_n}\leq n$ be the essential
return times or escape times of the orbit of each point of
$\omega$ and write $(l_1,s_1,j_1),\ldots,
(l_{u_n},s_{u_n},j_{u_n})$ for the corresponding critical
points and depths at each essential return situation, as in
Section~\ref{sec:fundamental-lemma}.  We recall also that
$d=d_n(\omega)$ is the number of pairs $(l_i,s_i)$ such that
$l_i+s_i\ge \Theta$, where $\Theta=\Theta(\flat)$ is defined
in~\eqref{eq:Theta}.

We consider the sequence of deep return situations of $\omega$:
$1\le r_1 < \dots < r_{d} \le u_n$ among the sequence of
return situations and then define
\[
\cD^{\flat}_n(x)=\sum_{k=1}^{d}
\big(|l_{r_k}|+|s_{r_k}|\big), \quad x\in\omega,
\]
which is constant on the elements of $\cP_n$, and get the
following bound.

\begin{proposition}
  \label{pr:depthbound}
  There exists $B_0=B_0(\sigma,\rho,\tau)>0$ such that for
  every $\omega\in\cP_n$ such that $d_n\mid\omega\ge2$ we
  have $\cC^{\flat}_n(x) \le \frac{B_0}n \cdot
  \cD^{\flat}_n(x)$ for all $x\in\omega$.
\end{proposition}

We start by proving the following.

\begin{lemma}
  \label{le:sumbinding}
  Let $1\le t_i<n$ be a essential return, inessential return
  or escape situation for $\omega\in\cP_n$, with binding
  time $p_i=p(l_i,s_i)$ (which we set to zero in the case of
  an escape situation). Then there exists a constant $B_1>0$
  such that for all $x\in\omega$
\[
\cC(t_i,t_i+p_i)
= \sum_{k=t_i}^{t_i+p_i}-\log\dist_\flat\left(f^k(x),\cC\right)
\le B_1\cdot(|l_i|+|s_i|).
\]
Moreover if $p_i>0$, then
$\dist\big(f^{t_i+k}(x),\cC\big)>
\dist\big(f^{t_i}(x),\cC\big)$ for $k=1,\dots,p$ as long
as $\rho>0$ is small enough.  In particular
$\cC(t_i,t_i+p_i)=0$ for escape situations or return
situations which are not deep, i.e.  the host interval
$I(l,s,j)$ is such that $|l|+|s|\le\Theta$.
\end{lemma}

\begin{proof}
  Let us fix $x\in\omega$ in what follows.  We consider
  first the case of $t_i$ being an essential return or
  escape situation with $|l_i|+|s_i|\le\Theta$,
  i.e. $t_i\neq t_{r_k}$ for all $k=1,\dots,d$. Thus $
  \log\dist_\flat\left(f^{t_i}(x),\cC\right)=0$. If $t_i$
  is an escape situation there is nothing else to prove
  because $p_i=0$.  Otherwise $p_i>0$ and we have two
  possibilities during the binding times $t_i+k$ with
  $k=1,\dots,p_i$: either $\big|
  f^{t_i+k}(x)\big|>\epsilon$, in which case we have again
  $\log\dist_\flat\left(f^{t_i+k}(x),\cC\right)=0$; or else
  $\big| f^{t_i+k}(x)\big|\le\epsilon$. In this last case by
  Proposition~\ref{pr:bindingaway} and the definition of
  binding time, if $ \dist\big(f^{t_i+k}(x),\cC\big)\le
  \dist\big(f^{t_i}(x),\cC\big), $ then
\[
\rho_0\cdot e^{-\rho k} \le
\dist\big(f^{t_i+k}(x),\cC\big)\le
\dist\big(f^{t_i}(x),\cC\big)
\le
a_2\cdot e^{-\frac{\pi}{\beta}(|l_i|+|s_i|)}
\]
and thus
\[
\frac{2\pi}{\beta\log\sigma}(|l_i|+|s_i|)\ge
p_i\ge k \ge
\frac1\rho\log\frac{\rho_0}2 +\frac{\pi}{\beta\rho}(|l_i|+|s_i|)
\]
which is impossible as long as
$\pi/(\beta\rho)>2\pi/(\beta\log\sigma)$.
Therefore choosing $k_0$ big enough and $\rho$ sufficiently
small we get $\dist\big(f^{t_i+k}(x),\cC\big)>
\dist\big(f^{t_i}(x),\cC\big)$ for all $k=1,\dots,p_i$
and so $\cC(t_i,t_i+p_i)=0$ for escapes and returns which
are not deep. This proves the last part of the statement of
the lemma.

On the other hand, if $t_i=t_{r_k}$ for some $1\le k \le d$,
we get
  \begin{equation}
    \label{eq:depth1}
    \dist_{\flat}(f^{t_{r_k}}(x),\cC)\ge 
    a_2\cdot e^{-\frac{\pi}{\beta}(|l_{r_k}+|+|s_{r_k}|)}
  \end{equation}
and so if $t_{r_k}$ is an escape or a return we have the
contribution 
\[
-\log \dist_{\flat}(f^{t_{r_k}}(x),\cC) \le -\log a_2 +
\frac{\pi}{\beta}(|l_{r_k}|+|s_{r_k}|)
\]
for the sum $\cC_n^\flat(x)$. For escape situations the
proof ends here. 

However, for return times we must consider the next binding
period. We stress that we are assuming \eqref{eq:CPR0}
holds.

Now we have for $h=1,\dots,p_{r_k}$, $l=l_{r_k}$ and
$i=r_k$
\begin{equation}
  \label{eq:dist0}
 \dist(f^{t_i+h}(x),\cC) \ge
\dist\big(f^h(x_l),\cC\big)-\dist\big( f^h(x_l), f^{t_i+h}(x)\big).
\end{equation}
Now we consider the two cases in the relations
\eqref{BCeq}.
\begin{description}
\item[Case $|f^h(x_l)|\le\epsilon$] 
  we obtain that \eqref{eq:dist0} is bounded from below by
  \begin{align}
    \label{eq:distC}
    (1-e^{-\tau h})\dist\big(f^h(x_l),\cC\big);
  \end{align}
\item[Case $|f^h(x_l)|>\epsilon$] we note that
  from \eqref{eq:condhaty} we have
  \begin{align*}
    \epsilon^{\tau}<\frac{1-e^{-\pi/\beta}}2
    \quad\text{so}\quad
    \epsilon^{1+\tau}<\epsilon\frac{1+e^{-\pi/\beta}}2
    =\epsilon-x_{k_0}<\dist(f^h(x_l),\cC);
  \end{align*}
  this ensures that \eqref{eq:dist0} is bounded from below
  by the same expression~\eqref{eq:distC} above.
\end{description}
Hence we deduce
\begin{equation}
  \label{eq:dist1}
  \sum_{h=1}^{p_i} -\log \dist(f^{t_i+h}(x),\cC) 
\le
C(\tau) + \sum_{h=1}^{p_i} - \log \dist\big(f^h(x_l),\cC\big).
\end{equation}
To bound the last sum we use the assumption \eqref{eq:CPR0}
on free times of the orbit of $x_l$.  We need to sum up to
the first \emph{free time} of the orbit of $x_l$ after
$p_i$. But if $p_i$ is a bound time for the orbit of $x_l$,
then by \cite[Lemma 5.3(a)]{PRV98} its binding period must
be smaller than $\frac{2\rho}{\log\sigma} p_i$. Thus there
exists a free time $n$ for the orbit of $x_l$ with $p_i\le n
\le \big(1+\frac{2\rho}{\log\sigma}\big) p_i$. Hence
\[
\sum_{h=1}^{p_i} - \log \dist\big(f^h(x_l),\cC\big)
\le \hat M \cdot \big(1+\frac{2\rho}{\log\sigma}\big)
\cdot p_i
\]
and so \eqref{eq:dist1} is bounded by $C(\tau)+\hat M \cdot
\big(1+\frac{2\rho}{\log\sigma}\big) \cdot p_i$.

For $\epsilon$ small enough we have that $k_0$ and $p_i$ are
very big and in both cases \eqref{eq:dist0} we get
$\cC(t_{r_k},t_{r_k}+p_{r_k})\le\tilde M \cdot p_{r_k}$ for a
constant $\tilde M$. By Lemma~\ref{le:maisdistorcao}(a) we
obtain
\[
\cC(t_{r_k},t_{r_k}+p_{r_k})\le  B_1\cdot(|l_{r_k}|+|s_{r_k}|)
\]
for a constant $B_1>0$, concluding the proof of the lemma.
\end{proof}

Next we show that the depth of an inessential return or free
time is not greater than the depth of the essential return
situation that precedes it.

\begin{lemma}
  \label{lem:essen-geq-inessen} Let $t_i$ be an essential
  return or an escape for $\omega\in \mathcal{P}_{n}$ with
  $I(l_i,s_i,j_i) \subset f^{t_i}(\omega)\subset
  I(l_i,s_i,j_i)^+$.  Then for each consecutive inessential
  return $t_i<t_{i}(1)<\dots<t_{i}(v)<n$ before the next
  essential return or escape and for each free time, i.e.
  for all iterates $t_i(k)+p_i(k)<j\le t_i(k+1)$ for
  $k=0,\dots,v-1$ with $p_i(k)$ the binding period of the
  $k$th return, the host interval $I(l,s,j)\supset
  f^{j}(\omega)$ is such that
  $|l|+|s|< |l_i|+|s_i|$.
\end{lemma}

\begin{proof}
  By items (1) and (3) of Lemma~\ref{lem:fn-omega-estimate}
  we have $|f^{j}(\omega)| > |f^{t_i}(\omega)| >
  |I(l_i,s_i,j_i)|$.  Thus because each $j$ is an
  inessential return or a free time we get
\[
a_{1}\frac{e^{-(\pi/\beta)(|l|+|s|)}}
{(|l|+|s|)^3}
>
|f^{t_{i}(k)}(\omega)|
>
a_{1}\frac{e^{-(\pi/\beta)(|l_i|+|s_i|)}}{(|l_i|+|s_i|)^3}.
\]
As $z^{-3}\cdot e^{-(\pi/\beta)z}$ is decreasing for $z>0$,
we conclude that $|l|+|s|< |l_i|+|s_i|$.
\end{proof}

If we have deep returns we can use sharper bounds to obtain
a relation between the logarithmic distance to the critical
set along the orbit between consecutive essential returns
and the depth of the first return. Note that by
(\ref{eq:expansaolivre}) on ``deep'' free iterates between
critical points we may write
\begin{equation}
  \label{eq:expansaolivre2}
  |f'(x)|
  \ge
  \left(\frac{\tau \hat x^{\alpha-1}}{C}
    \cdot e^{-(1-\alpha)\frac{\pi}{\beta}|s(\tau)|}\right)
    \cdot e^{(1-\alpha)\frac{\pi}{\beta}(|l|+|s|)}
    >  e^{(1-\hat\alpha)\frac{\pi}{\beta}(|l|+|s|)}
\end{equation}
for some $\alpha<\hat\alpha<1$, since on these times we have
$|s|\le s(\tau)$ and we assume that $|l|+|s|\ge\Theta$ is so
big that the expression in parenthesis in
(\ref{eq:expansaolivre2}) is larger than
$e^{-(\hat\alpha-\alpha)\Theta}$.

\begin{lemma}
  \label{le:logsoma}
  Let $t_i$ be a deep essential return or an escape for
  $\omega\in \mathcal{P}_{n}$ with $I(l_i,s_i,j_i) \subset
  f^{t_i}(\omega)\subset I(l_i,s_i,j_i)^+$. Let
  $t_i=t_i(0)<t_{i}(1)<\dots<t_{i}(v)<t_i(v+1)=t_{i+1}$ be
  the consecutive inessential returns before the next
  essential return or escape and let $p_i(k)$ be the
  corresponding binding times for $k=0,\dots,v+1$. Define
  $J=\cup_{k=0}^{v}\{t_i(k)+p_i(k)+1, \dots,t_{i}(k+1)\}$
  the set of free iterates and of inessential return times.
  Assume that every such iterate $h\in J$ is contained in a
  ``deep'' interval $I(l_h,s_h,j_h)$, i.e.
  $l_h+s_h\ge\Theta$. Then there exists a constant $B_2>0$
  such that
$
\sum_{h\in J} -\log\dist_\flat\big( f^h(x) , \cC) \le
B_2(|l_i|+|s_i|).
$
\end{lemma}

\begin{proof}
  By the choice of $\Theta$ we have $Q>1$ in
  Lemma~\ref{lem:fn-omega-estimate} and
  (\ref{eq:expansaolivre2}) also holds. Then we deduce
\[
\big| f^{t_{i+1}}(\omega) \big| 
\ge 
\prod_{k=0}^{v} e^{(1-3\zeta)\frac{\pi}{\beta}
( |l_{t_i(k)}| + |s_{t_i(k)}|)} 
\prod_{h=t_i(k)+p_i(k)}^{t_i(k+1)-1}
e^{(1-\hat\alpha)\frac{\pi}{\beta}(|l_h|+|s_h|)}
\cdot |f^{t_i+p_i}(\omega)|
\]
and then since we may assume that $3\zeta<\hat\alpha$, by
the definition of essential returns and using
Lemma~\ref{le:maisdistorcao}(b) we get
\[
\frac{e^{-\frac{\pi}{\beta}(|l_{i+1}|+|s_{i+1}|)}}
{(|l_{i+1}|+|s_{i+1}|)^3}
\ge
\exp\left(
\frac{\pi}{\beta}\Big(
(1-\hat\alpha)\sum_{h\in J} (|l_h|+|s_h|)
\Big)
\right)
\cdot
\frac{e^{-3\zeta\frac{\pi}{\beta}(|l_{i}|+|s_{i}|)}}
{(|l_{i}|+|s_{i}|)^3}
\]
or equivalently
\begin{equation}
  \label{eq:eraumavez}
3\log\left(
\frac{|l_{i}|+|s_{i}|}{|l_{i+1}|+|s_{i+1}|}
\right)
+
3\zeta\frac{\pi}{\beta}(|l_{i}|+|s_{i}|)
\ge
\frac{\pi}{\beta}(|l_{i+1}|+|s_{i+1}|)
+
\frac{\pi}{\beta}\Big(
(1-\hat\alpha)\sum_{h\in J} (|l_h|+|s_h|)
\Big).
\end{equation}
Since $|l_{i+1}|+|s_{i+1}|\gg1$ we also have that the left
hand side of (\ref{eq:eraumavez}) is smaller than
\[
\left(
3\zeta + \frac{3\beta \log (|l_{i}|+|s_{i}|)}{\pi(|l_{i}|+|s_{i}|)}
\right)
\frac{\pi}{\beta}(|l_{i}|+|s_{i}|)
\le
\left(
3\zeta + \frac{3\beta \log \Theta}{\pi\cdot\Theta}
\right)
\frac{\pi}{\beta}(|l_{i}|+|s_{i}|)
\le
5\zeta\frac{\pi}{\beta}(|l_{i}|+|s_{i}|)
\]
for $\Theta$ sufficiently big, because
$|l_{i}|+|s_{i}|\ge\Theta$ and $\log(z)/z$ is decreasing for
$z>e$. Thus we get
\begin{equation}
  \label{eq:eila}
5\zeta(|l_{i}|+|s_{i}|)
\ge
|l_{i+1}|+|s_{i+1}|
+
(1-\hat\alpha)\sum_{h\in J} (|l_h|+|s_h|).
\end{equation}
Now since every iterate is ``deep'' for all $x\in\omega$ we
have the bound
\[
\sum_{h\in J} -\log\dist_\flat\big(f^h(x),\cC\big)
\le
-\#J\cdot\log a_2
+
\frac{\pi}{\beta}\sum_{h\in J} (|l_h|+|s_h|)
\le
B_2\cdot (|l_{i}|+|s_{i}|)
\]
for a constant $B_2>0$ depending on $\zeta$ and $\hat\alpha$
from (\ref{eq:eila}), as long as $|l_h|+|s_h|\ge\Theta$ is
sufficiently big.  This completes the proof of the lemma.
\end{proof}

\begin{proof}[Proof of Proposition~\ref{pr:depthbound}:\,]
  Let us fix $x\in \omega\in\cP_n$ with $d_n\mid\omega\ge2$ and
  $i\in\{1,\dots,u_n(x)-1\}$.  According to
  Remark~\ref{rmk:partitiondimension} and the definition of
  deep essential return situation we have
  \begin{equation}
    \label{eq:depth01}
    \dist_{\flat}(f^{r_i}(x),\cC)\ge 
    a_2\cdot e^{-\frac{\pi}{\beta}(|l_i|+|s_i|)}.
  \end{equation}
  Note that the above truncated distance is $1$ on the
  return situations $t_i$ which are not deep, by the choice
  of $\Theta$. Moreover this distance is also $1$ for all
  iterates between such $t_i$ (not deep) and the next return
  situation $t_{i+1}$ by Lemma~\ref{lem:essen-geq-inessen}
  for the inessential return and free iterates, and by
  Lemma~\ref{le:sumbinding} for the bound iterates.

  Hence \emph{we only have to take care of the deep
    essential return or escape times} plus the next iterates
  before the following essential return situation. The sum
  of the logarithms of the truncated distance on binding
  periods, given by Lemma~\ref{le:sumbinding}, is bounded by
  a constant times the depth of the return which originated
  the binding. In addition, the same sum over the free
  and the inessential return iterates is likewise bounded by
  the depth of the essential return or escape $t_i$, by
  Lemma~\ref{le:logsoma}. If we keep the notations
  introduced in the statements of
  Lemmas~\ref{lem:essen-geq-inessen} and~\ref{le:logsoma},
  then we may write
  \begin{eqnarray*}
\cC(t_i,t_{i+1})
&\le&    
\sum_{k=0}^v\Big[
\cC\big(t_i(k),t_i(k)+p_i(k)\big)
+
\cC\big(t_i(k)+p_i(k)+1,t_i(k+1)\big)
\Big]
\\
&\le&
B_1 \sum_{k=0}^v \Big(\big|l_{t_i(k)}\big| +
\big|s_{t_i(k)}\big|\Big)
+
B_2(|l_i|+|s_i|)
\\
&\le&
B_1 \sum_{k=0}^v \cC\big(t_i(k)+p_i(k)+1,t_i(k+1)\big)
+
B_2(|l_i|+|s_i|)
\\
&\le&
B_2\cdot(1+B_1)\cdot(|l_i|+|s_i|) 
  \end{eqnarray*}

  Setting $B_0=B_2(1+B_1)$ this finishes the proof of
  Proposition~\ref{pr:depthbound}.
\end{proof}

\subsection{The expected value of the distance at return
  times}
\label{sec:expect-value-dist}

The statement of Proposition~\ref{pr:depthbound} together
with Lemma~\ref{le:escape-ae} and
Proposition~\ref{prop:depth-probability-1} ensure that, to
obtain slow recurrence to the critical set, we need to bound
$\cD_n^\flat(x)/n$ for Lebesgue almost every $x\in I$.
Indeed we have for every big enough $n$
\[
\big\{x\in I:  \cC^\flat_n(x)>\de \big\}
\subseteq
\bigcup\big\{
\omega\in\cP_n : u_n\mid\omega\equiv1
\big\}
\cup
\Big\{x\in I: u_n(x)\ge2\mbox{ and }
\cD^\flat_n(x)>\frac{n}{B_0}\cdot \de \Big\}
\]
and Lemma~\ref{le:escape-ae} shows that the left hand side subset
of the above union has exponentially small measure. We now
show that Proposition~\ref{prop:depth-probability-1} implies
a similar bound for the right hand subset.

\begin{lemma}
  \label{le:esperanca}
  For all $z\in\big(0,\frac{(1-2\vartheta)\pi}{2\beta}\big)$
  there is $\Theta_1$ so that $ \int
  e^{z\cdot\cD_n^\flat(x)}\, d\lambda(x) \le
  e^{o(\Theta)\cdot n}$ for
  $\Theta>\Theta_1=\Theta_1(z,\tau,\rho,\sigma)$, where
  $o(\Theta) \to 0$ when $\Theta\to\infty$.
\end{lemma}

\begin{proof}
The integral in the statement equals the following series
\[
\sum_{
\genfrac{}{}{0pt}{}{1\le v \le u \le n}
{(\eta_1,\upsilon_1),\ldots,(\eta_v,\upsilon_v)}}
\exp\left(
z \sum_{k=1}^{v} (\eta_k+\upsilon_k)
\right)
\cdot
\lambda\left(
A_{(\eta_1,\upsilon_1),\ldots,(\eta_v, \upsilon_v)}^{u,v}(n)
\right),
\]
where $\eta_k+\upsilon_k\ge \Theta$, $\upsilon_k\ge1$ and
$\eta_k\ge k_0$ for $k=1,\dots,v$.
Proposition~\ref{prop:depth-probability-1} provides the
bound
\[
\sum_{
\genfrac{}{}{0pt}{}{1\le v \le u \le n}
{(\eta_1,\upsilon_1),\ldots,(\eta_v,\upsilon_v)}
}
\hspace{-0.6cm}
\binom{u}{v}
e^{z\sum_{k=1}^v (\eta_k+\upsilon_k)+
( 2\vartheta -1) \frac{\pi}{\beta}
\sum_{k=1}^v (\eta_k+\upsilon_k)}
=
\hspace{-0.8cm}
\sum_{
\genfrac{}{}{0pt}{}{1\le v \le u \le n}
{(\eta_1,\upsilon_1),\ldots,(\eta_v,\upsilon_v)}
}
\hspace{-0.6cm}
\binom{u}{v}
e^{
( z+ ( 2\vartheta -1) \frac{\pi}{\beta})\cdot
\sum_{k=1}^v (\eta_k+\upsilon_k)}.
\]
Now setting $\Delta=\sum_{k=1}^v (\eta_k+\upsilon_k)$ and
\[
K(v,\Delta)=\#\Big\{ ((l_1,s_1),\dots,(l_v,s_v)):
\sum_{k=1}^v(l_k+s_k)=\Delta,\,
l_k\ge k_0,\,
s_k\ge1, \, l_k+s_k\ge\Theta,\,
1\le k \le v
\Big\}
\]
we may rewrite the last series as
$
\sum_{1\le v\le u\le n}
\sum_{\Delta\ge v\Theta}
\binom{u}{v}
K(v,\Delta)\cdot e^{(z+ ( 2\vartheta -1) \frac{\pi}{\beta}) \Delta}.
$
To estimate $K(v,\Delta)$ we observe that 
\[
K(v,\Delta)\le
\#\Big\{(n_1,\dots,n_{2v}): \sum_{k=1}^{2v} n_k = \Delta
\mbox{  and  } n_k\ge0, k=1,\dots,2v\Big\}
=
\binom{\Delta+2v-1}{2v-1},
\]
where $\binom{n}{k}=\frac{n!}{k!\cdot (n-k)!}$ is a binomial
coefficient.  By a standard application of
Stirling's Formula we get
\[
K(v,\Delta)\le \Big(
C^{1/\Delta}\cdot\big(1+\frac{2v-1}{\Delta}\big)
\cdot\big(1+\frac{\Delta}{2v-1}\big)^{(2v-1)/\Delta}
\Big)^\Delta
\le e^{z\cdot \Delta},
\]
since $\Delta\ge v\Theta$ ensures that the expression
in parenthesis can be made arbitrarily close to $1$ if $\Theta$
is taken bigger than some constant $\Theta_0=\Theta_0(z)$, where
$0<C<1$ is a constant independent of $\Theta$ and we assume
that $z>0$ is small. Hence we arrive at
\begin{align*}
\int e^{z \cD_n^\flat(x)}\, d\lambda(x)
&\le
\sum_{1\le v\le u\le n}
\sum_{\Delta\ge v\Theta}
\binom{u}{v} e^{(2z+ ( 2\vartheta -1) \frac{\pi}{\beta}) \Delta}
\le
\sum_{v=0}^n \binom{n}{v}
\cdot C \cdot e^{(2z+ ( 2\vartheta -1) \frac{\pi}{\beta})
  \Theta v}
\\
&\le
\Big( 1 +
C \cdot e^{(2z+ ( 2\vartheta -1) \frac{\pi}{\beta})\Theta }
\Big)^n = e^{o(\Theta)\cdot n},
\end{align*}
as long as $0< z < (1-2\vartheta)\pi/(2\beta)$ and $\Theta>
\Theta_1>\max\{\theta,\Theta_0\}$ is big enough so that
$Q>1$ in Lemma~\ref{lem:fn-omega-estimate} and
(\ref{eq:expansaolivre2}) holds, as stated.
\end{proof}

As a consequence of this bound we can use Tchebychev's
inequality with $z$ and $\Theta$ as in the statement of
Lemma~\ref{le:esperanca} to obtain
\[
\lambda\left(\Big\{ \cD^\flat_n \ge \frac{n}{B_0}\cdot \de \Big\}\right)
=
\lambda\Big(\big\{
e^{z \cD^\flat_n} > e^{z\cdot\de\cdot n/B_0}
\big\}\Big)
\le
e^{-z\cdot\de\cdot n/B_0}
\int  e^{z\cdot\cD^\flat_n}\, d\lambda
\le e^{-z\cdot\de\cdot n/B_0}\cdot e^{o(\Theta)\cdot n}.
\]
Now observe that we may take $\flat=\flat(\delta)>0$ so small that
$\Theta(\flat)$ becomes big enough to satisfy all constraints
and moreover $o(\Theta)<\delta\cdot\frac{z}{n}$.

As already explained, this together with
Lemma~\ref{le:escape-ae} implies that there are $C,\xi>0$,
where $\xi=\xi(\delta)$, such that $\lambda\big(\{ x\in I:
\cC_n(x)>\de \}\big)\le C e^{-\xi n}$ for all $n\ge1$.  Then
since
\[
\{ x\in I : \mathcal R(x)>n \}
\subseteq
\bigcup_{k>n}
\{x\in I: \cC_k(x)>\de\}
\] 
we conclude that there are
constants $C_1,\xi_1>0$ such that 
Theorem~\ref{thm:approxestimate} holds.



\section{Fast expansion for most points}
\label{sec:fast-expansion-most}

Here we use the results from the previous sections to prove
Theorem~\ref{thm:expansionestimate} and as a consequence
obtain Corollary~\ref{cor:decayCLT}. We start
by setting
\[
E_n=\Big\{ \omega\in\cP_n : \exists\, 1\le k\le n
\mbox{ s.t. }
\dist(f^k(\omega),\cC)<e^{-\rho\cdot n}\, ,
f^k(\omega)\subset I(l,s)
\mbox{ and }
|l|\ge k_0, |s|>s(\tau)
\Big\}
\]
and proving the following bound.
\begin{lemma}
\label{le:there-are-constants}
There are constants $C,\xi>0$ dependent on $\hat f$, $k_0$,
$\zeta$, $\rho$, and $\tau$ only such that 
$\lambda\big( \cup E_n \big)\le C \cdot e^{-\xi\cdot n}$
for all $n\ge1$.
\end{lemma}

\begin{proof}
  Let us take $\omega\in E_n$ and let $k\in\{1,\dots, n\}$ be the
  iterate which is very close to the critical set.  Observe
  that by Remark~\ref{rmk:partitiondimension} the constraint
  on the distance implies 
\[
-\rho n > \log\dist(f^k(\omega),\cC) \ge \log a_2
 -\frac{\pi}{\beta}(|l|+|s|)
\quad\mbox{and so}\quad
|l|+|s|\ge \frac{\beta\rho}{20\pi}\cdot n.
\]
Since this iterate is in the binding region, there must be
an essential return $t<k$, $t\in R_n(\omega)$, whose depth
is at least as larger as $|l|+|s|$, by the results of
Lemmas~\ref{le:sumbinding} and~\ref{lem:essen-geq-inessen}.

Hence, according to the definition of
$A_{(\eta,\upsilon)}(n)$ from
Section~\ref{sec:fundamental-lemma}, if $n\ge\Theta$
\begin{align*}
\bigcup E_n
&\subset
\bigcup \big\{
A_{(\eta,\upsilon)}(n):
(\eta,\upsilon)\,\,\mbox{is such that}\,\,
\eta+\upsilon\ge\frac{\beta\rho}{20\pi}\cdot n
\big\}.
\end{align*}
Thus by Corollary~\ref{cor:c.1} we can estimate
\begin{align*}
\lambda\Big( \bigcup E_n \Big)
&\le
n^3 e^{o(\Theta)n}\cdot \hspace{-0.7cm}
\sum_{\genfrac{}{}{0pt}{}
{\eta+\upsilon\ge \beta\rho n/(20\pi)}
{\eta\ge k_0, \upsilon\ge s(\tau)}}
\hspace{-0.7cm}
e^{(4\vartheta-1)\frac{\pi}{\beta}(\eta+\upsilon)}
\le 
n^3 e^{o(\Theta)n}\cdot\hspace{-0.7cm}
\sum_{\Delta\ge \beta\rho n/(20\pi)} \hspace{-0.7cm}
\Delta e^{(4\vartheta-1)\frac{\pi}{\beta}\Delta}
\\
&\le C \cdot n^3 e^{o(\Theta)n}\cdot
e^{(4\vartheta-1)\frac{\rho}{20}n}
\le
C\cdot e^{(4\vartheta-1)\frac{\rho}{100}n}
\end{align*}
for some constant $C>0$ with
$\xi=(4\vartheta-1)\frac{\rho}{100}$ for a big enough
$\Theta$. This finishes the proof of the lemma.
\end{proof}

\begin{lemma}
  \label{le:n3}
  If $n$ is big enough, $\rho$ small enough (depending only
  on $\sigma$ and $A$ from Lemma~\ref{le:distorcao}) and
  $x\in I\setminus\cup E_n$, then $\big| (f^n)'(x)
  \big|\ge\sigma^{n/3}$.
\end{lemma}

\begin{proof}
  Let us take $x\in \omega\in\cP_n\setminus E_n$ and let
  $0<r_1<\dots<r_k<n$ be the consecutive returns (either
  essential or inessential) of the first $n$ iterates of the
  orbit of $\omega$, and $p_1,p_2,\dots,p_k$ the respective
  binding periods. We also set $q_i=r_{i+1}-(r_i+p_i+1)$ the
  free periods and possibly escape times between consecutive
  returns, for $i=1,\dots,k-1$, $q_0=r_1$ and
  $q_{k}=n-(r_k+p_k+1)$ if $n > r_k+p_k$ or $q_{k+1}=0$
  otherwise.

We split the argument in the following two cases. If $n >
r_k+p_k$ then
\[
\big|
(f^n)'(x)
\big|
=
\prod_{i=0}^k \Big(
\big|
(f^{q_i})'(f^{r_i+p_i+1}(x))
\big|
\cdot
\big|
(f^{p_i+1})'(f^{r_i}(x))
\big|
\Big)
\ge
\sigma_0^{\sum_{i=0}^{k+1}q_i}\cdot
A_0^{k}\cdot \sigma^{\sum_{i=1}^k (p_i+1)/3}
\ge \sigma^{n/3},
\]
since $\sigma_0>\tilde\sigma>\sigma$ and $A_0>1$ by
Lemma~\ref{le:maisdistorcao}(c), and also by
(\ref{eq:expansaolivre}) we may assume that at escape times
the expansion rate is at least $\sigma$.

On the other hand, if $n\le r_k+p_k$ then using
Lemma~\ref{le:distorcao}, Lemma~\ref{l3.1}(2) and that
$\omega\in\cP_n\setminus E_n$
\begin{align*}
\big|(f^n)'(x)\big|
&=
\big|(f^{r_k})'(x)\big|
\cdot
\big|f'(f^{r_k}(x))\big|
\cdot
\big|(f^{n-r_k-1})'(f^{r_k+1}(x))\big|
\\
&\ge
\big|(f^{r_k})'(x)\big|
\cdot
C^{-1} |x_{l_{k}}|^{\alpha-2}e^{-\rho n}
\cdot
\frac1{A}\cdot \big| (f^{n-r_k-1})'(x_l)\big|
\\
&\ge
(CA)^{-1}\cdot\sigma^{r_k}\cdot 
e^{-\rho n(\frac{\alpha-2}{20}-\frac{s(\tau)}{\rho n}) } 
\cdot \sigma^{n-r_k-1}
\\
&\ge 
\exp\left(
n\Big(\frac{n-1}n\log\sigma
-\frac{\log (CA)}{n}
-\rho\big( \frac{\alpha-2}{20}-\frac{s(\tau)}{\rho n} \big) 
\right)
\ge\sigma^{n/3},
\end{align*}
for $\rho>0$ small enough and $n$ big enough, where $x_l$ is
the critical point associated to $f^{r_k}(\omega)$ and we
have used also the calculation for the previous case to
estimate $|(f^{r_k})'(x)|$.
\end{proof}

Finally since 
\begin{align*}
  \{ x\in I: \mathcal E(x)>n \} \subseteq 
  \bigcup_{k>n} \left(
  \bigcup\Big\{\omega\in\cP_k:
  \big|(f^k)'(x)\big|<\sigma^{k/3}, x\in\omega \Big\}
  \right)
\end{align*}
we conclude from Lemmas~\ref{le:there-are-constants}
and~\ref{le:n3} that there are $C_2,\xi_2>0$ such that
\begin{align*}
  \lambda\Big(\big\{ x\in I : \mathcal E(x)>n
  \big\}\Big)\le\sum_{k>n} \lambda(\cup E_k) \le C_2\cdot
  e^{-\xi_2\cdot n},
\end{align*}
concluding the proof of Theorem~\ref{thm:expansionestimate}.


\section{Exponential bound on derivatives along critical
  orbits}
\label{sec:extra-excl-param}

Here we explain how to obtain the bound \eqref{eq:CPR0} for
the parameters in the set $S$.
First we claim that it is enough to obtain the following
bound for a sufficiently small value of $\flat>0$
\begin{equation}
  \label{eq:CPR2}
  \cC^\flat_n(z_k)
= 
\sum_{j=0}^{n-1} -\log \dist_\flat\big( f^j_\mu(z_k)
  ,\cC\big)\le 
  \hat M\cdot n,
\end{equation}
where $\hat M>0$ is a big constant and this bound holds for
all $n\ge1$ which is \emph{not a bound time} for $z_k$, for every
critical value $z_k$, $|k|\ge k_0$. Indeed, fixing $n$, $k$
and $\flat>0$, if \eqref{eq:CPR2} holds then we can write
\begin{align*}
\sum_{j=0}^{n-1}
-\log \dist\big( f^j_\mu(z_k)  ,\cC\big)
&\le 
\hspace{-0.7cm}
\sum_{\genfrac{}{}{0pt}{}{\dist(f^j(z_k),\cC)<\flat}{0\le
  j < n}}
\hspace{-0.7cm}
 -\log \dist_\flat\big( f^j_\mu(z_k)
  ,\cC\big)
+
\hspace{-0.7cm}
\sum_{\genfrac{}{}{0pt}{}{\dist(f^j(z_k),\cC)\ge\flat}{0\le
  j < n}}
\hspace{-0.7cm}
-\log\flat
\\
&\le
\hat M\cdot n  -\log\flat\cdot\#\{0\le
  j < n : \dist(f^j(z_k),\cC)\ge\flat  \}
\\
&\le
\big( \hat M -\log\flat \big)\cdot n,
\end{align*}
proving the claim. 

For any given $\mu\in S$ and $|k|\ge k_0$, let $n$ satisfy
$1= r_0<r_1 < \dots < r_u \le n$, where $r_1,\dots,r_u$ are
the essential return situations of the orbit of $z_k$, in
the sense of the construction performed in \cite{PRV98}. We
denote by $(l_i,s_i)$ the depth corresponding to $r_i$ for
$i=1,\dots,u$, and set $(l_0,s_0)$ to be such that
$z_k(\mu)=f(x_k)\in I(l_0,s_0)$. To
prove~\eqref{eq:CPR2} it suffices to obtain the following
pair of relations
\begin{equation}
  \label{eq:CPRV3}
\cC^\flat_n(z_k)
\le
B\sum_{i=0}^{u_n(k)} ( |l_i| + |s_i| )
\le
B\cdot\frac{n}2,
\end{equation}
where $u_n(k)$ denotes the number of return situations of
the orbit of $z_k$ up to the n$th$ iterate, for every time
$n\ge1$ which is \emph{not a bound time} for the orbit of
$z_k$, for all critical values $z_k, \, |k|\ge k_0$, and $B$
is a positive constant.  Obviously both inequalities in
(\ref{eq:CPRV3}) together imply~\eqref{eq:CPR2} with $\hat
M=B/2$.

\subsection{The left hand side from the right hand side}
\label{sec:left-hand-side}

To obtain~(\ref{eq:CPRV3}) \emph{we first assume that the right
hand side inequality has been} proved \emph{for all $n\ge1$,
$|k|\ge k_0$ and $\mu\in S$ and} deduce \emph{the left hand side
inequality in the same setting by induction on the number
$n$ of iterates, as follows. }

According to \cite[Section 4]{PRV98} we have that for a
fixed $\alpha<\gamma<1$
\begin{align*}
  \epsilon> \big| f^i_\mu\big(z_k(\mu)\big) \big|
  \ge 
  \big|z_k(\mu)\big|^{\gamma^{i-1}}
  \quad\text{and}\quad
  |(f_\mu^i)^\prime(z_k(\mu))|\ge\epsilon^{(\gamma-1)i}
\end{align*}
for every $1\le i < j_0=j_0(k,\mu)$, where $j_0$ is the
first iterate of $f_\mu$ such that
$|f_{\mu}^j(z_k(\mu))|>\epsilon$.  This threshold $j_0$ is
uniformly bounded from above (by $L$, say) for all $|k|\ge
k_0$ and $\mu\in S$. Consequently, since all the above
iterates are in the region of expansion between critical
points, there are $|l_i|\ge k_0$ and $s_i$ such that
\begin{align*}
  f^i_\mu\big(z_k(\mu)\big)\in I(l_i,s_i)
  \quad\mbox{with}\quad
  |s_i|\le s(\tau),
  \quad i=0,\dots,j_0-1.
\end{align*}
Hence by Remark~\ref{rmk:partitiondimension}
\begin{align*}
  \dist\Big(f^i_\mu\big(z_k(\mu)\big),\cC\Big) \ge a_2
  e^{-\frac{\pi}{\beta}(|l_i|+|s_i|)} \ge a_2
  e^{-\frac{\pi}{\beta} s(\tau)} e^{-\frac{\pi}{\beta}|l_i|}
\end{align*}
and since
\begin{align*}
  \big|  f_\mu^j(z_k(\mu)) \big|
  \le
  \dist\big( I(l_i-1,s) , 0 )
  \le
  (a_2+\hat x) e^{\frac{\pi}{\beta}} e^{-\frac{\pi}{\beta}|l_i|}
\end{align*}
we have
\[
\dist\Big(f^i_\mu\big(z_k(\mu)\big),\cC\Big)
\ge 
\frac{a_2 e^{-\frac{\pi}{\beta} s(\tau)}}
{e^{\pi/\beta}(a_2+\hat x)} \big|f^i_\mu\big(z_k(\mu)\big)\big|
\ge
K(\tau) \big| z_k(\mu) \big|^{\gamma^{i-1}}
\]
for $i=1,\dots,j_0-1$, with $0<K(\tau)<1$. For $i=0$ we obtain
\begin{align*}
  \dist\big(z_k(\mu),\cC\big)
  \ge 
  a_2 e^{-\frac{\pi}{\beta}(|l_0|+|s_0|)} 
  \ge
  a_2 e^{-\frac{\pi}{\beta} s(\tau)} e^{-\frac{\pi}{\beta}|l_0|}
  \ge
  K(\tau) \big| z_k(\mu) \big|.
\end{align*}
Because $-\log\dist_\flat(\cdot,\cC) \le
-\log\dist(\cdot,\cC)$ we get the bound
\begin{align}\label{eq:boundepth}
  \sum_{i=0}^{j_0-1} -\log\dist_\flat
  \Big(f^i_\mu\big(z_k(\mu)\big),\cC\Big)
  &\le
  -j_0\cdot\log K(\tau)
  -\Big(1+\sum_{i=1}^{j_0-1}\gamma^{i-1}\Big)\log\big| z_k(\mu) \big|.
\end{align}
Now using again Remark~\ref{rmk:partitiondimension} we have
that
\begin{align*}
  \big| z_k(\mu) \big|
  \ge
  \dist(I(l_0,s_0),\cC)
  \ge
  |\hat x -a_2|e^{-\frac{\pi}{\beta}|l_0|}
  >
  |\hat x -a_2|e^{-\frac{\pi}{\beta}(|l_0|+|s_0|)}
\end{align*}
and so \eqref{eq:boundepth} is bounded by
\begin{align*}
  -j_0&\log K(\tau) - C\log|\hat x -a_2|+C\frac{\pi}{\beta}(|l_0|+|s_0|)
  \\
  &=
  (|l_0|+|s_0|)\left(
    C\frac{\pi}{\beta}-\frac{j_0\log K(\tau)+\log|\hat x -a_2|}{|l_0|+|s_0|}
    \right)
  \\
  &\le
  (|l_0|+|s_0|)\left(
    C\frac{\pi}{\beta} -\frac{j_0\log K(\tau)+\log|\hat x -a_2|}{k_0+1}
    \right)
    \le \hat B\cdot(|l_0|+|s_0|)
\end{align*}
for a constant $\hat B>0$ which depends on $\tau$. Since
these initial iterates are all free, we have shown that the
left hand side of (\ref{eq:CPRV3}) holds for the $j_0-1$
initial iterates of the orbit of every critical value for
any $B\ge \hat B$.

Now assume that the left hand side of (\ref{eq:CPRV3}) is
true for a free time $n-1$ of the orbit of a critical value
$z_k$ for a fixed $\mu\in S$. 

If $n$ is a free time for
$z_k$ and the last essential return situation was not deep,
i.e. the depth was smaller than $\Theta$, then the depth of
$f^n(z_k)$ is also smaller than $\Theta$ by arguments
akin to Lemma~\ref{lem:essen-geq-inessen}. Hence this
situation does not contribute to the sum $\cC^\flat_n(z_k)$.

If $n$ is a free time for
$z_k$ and the last essential return situation was deep, then
let $r\le n-1$ be the last essential return time with depth
$(l_r,s_r)$ and $t> n$ be the next return situation
(either essential, inessential or escape). We claim that
\begin{equation}
  \label{eq:claim!}
  \sum_{j=r}^{t-1} -\log\dist_\flat
\Big(f^j_\mu\big(z_k(\mu)\big),\cC\Big)
\le
B\cdot (|l_r|+|s_r|),
\end{equation}
which shows that the induction can be carried out up to the
iterate $t-1\ge n$. 

To prove the claim, note that similar arguments to those
proving Lemma~\ref{le:logsoma} show that the part of the sum
(\ref{eq:claim!}) corresponding to free times after the
binding period is bounded by $B_2\cdot(|l_r|+|s_r|)$. Hence
we get (\ref{eq:claim!}) as long as $B\ge B_2$.

Finally, if $n$ is a return situation for $z_k$, either
essential, inessential or an escape, we let $p$ be the
binding period corresponding to the return and consider the
next free iterate $n+p$ of the orbit of $z_k$.

If $n$ is an escape, then $p=0$ and we are done by
Remark~\ref{rmk:partitiondimension}.  Otherwise $0<p\le
\frac{2\rho}{\log\sigma}\cdot n < n$ by \cite[Lemma
5.3]{PRV98} and so we can use the induction hypothesis to
get (\ref{eq:dist0}) with $t_i=n$, $h=1,\dots,p$, $x=z_k$
and $x_l=x_{l_n}$, the critical point which will shadow the
orbit of $z_k$ during the binding period. Hence we obtain as
before the bound (\ref{eq:dist1}).

Note that if $|l_n|+|s_n|\le\Theta$, then the truncated
distance is always $1$ and we are done, by the same
arguments in the proof of Lemma~\ref{le:sumbinding}.  If we
have a deep return then, analogously, we need to consider
the next free time $t\ge p$ of the orbit of the bound
critical value $z_{l_n}$ in order to properly use the
induction hypothesis. We have $t\le
(1+2\rho/\log\sigma)\cdot p$ again by \cite[Lemma
5.3]{PRV98} from which we get
\[
\sum_{j=n+1}^{n+p}
\hspace{-0.3cm}
-\log\dist_\flat
\Big(f^j_\mu\big(z_k(\mu)\big),\cC\Big)
\le
C(\tau) + \Big(1+\frac{2\rho}{\log\sigma}\Big)
\cdot p
\le
B\cdot
(|l_n|+|s_n|)
\]
by induction and assuming that the right hand side of
(\ref{eq:CPRV3}) holds. We also used the upper bound in Lemma
\ref{le:maisdistorcao}(a) and assumed that
$
B> C(\tau)/\Theta
+
\frac{2\pi}{\beta\log\sigma}
\cdot
\Big(1+\frac{2\rho}{\log\sigma}\Big).
$
Thus if $B$ is sufficiently big, then the inductive step can be
performed in every situation. 

\emph{This shows that the left hand
side inequality in (\ref{eq:CPRV3}) is true \emph{if} the right
hand side inequality holds.}

\subsection{The right hand side inequality}
\label{sec:right-hand-side}

Now we explain why we can assume that parameters $\mu\in S$
satisfy the right hand side of (\ref{eq:CPRV3}) for all
critical values. The condition (30) used in
\cite[p. 463]{PRV98} 
\begin{equation}
  \label{eq:cond30}
B(n,\omega,k) = \sum_{j=1}^{n-1} p(j,\omega,k) < \frac{n}2  
\end{equation}
to test whether a given interval $\omega$ of parameters should be
excluded or not, can be replaced by the following
\begin{equation}
  \label{eq:newcond}
  \sum_{i=0}^{u_n(k)} ( |l_i| + |s_i| )
  \le
  \cdot\frac{n}2
\end{equation}
without loss since $B(n,\omega,k)\le C\sum_{i=0}^{u_n(k)} (
|l_i| + |s_i| )$ by Claim (3) in \cite[p. 478]{PRV98}.
Indeed it is this last inequality which is used in the
arguments proving \cite[Lemma 5.7]{PRV98} establishing the
exponential bound on the measure of the set of excluded
parameters.  Therefore repeating the algorithm presented
there step by step with the new condition (\ref{eq:newcond})
instead of (\ref{eq:cond30}) leads to the construction of a
positive Lebesgue measure set $S$ satisfying
Theorem~\ref{th:ref-annals} and (\ref{eq:CPRV3}) for all
$n\ge1$, every $|k|\ge k_0$ and for every $\mu\in S$.  This
concludes the proof of (\ref{eq:CPR2}).


\section{Constants depend uniformly on initial parameters}
\label{sec:const-depend-unif}

We finally complete the proof of
Corollary~\ref{cor:contSRBEntropy} by explicitly showing the
dependence of the constants used in the estimates on
Sections~\ref{sec:idea-construction} to
\ref{sec:fast-expansion-most}.

In the statements of the lemmas and propositions in the
aforementioned sections we stated explicitly the direct
dependence of the constants appearing in each claim from
earlier statements. For constants which depend only on $\hat
f$ we used the plain letter $C$. 

It is straightforward to see that every
constant depends on values that ultimately rest on the
choice of initial values for $\sigma,\sigma_0$ and 
$k_0$  and on the the choice of $\rho$
and $\tau$, which are taken to be small enough and where
$0<\rho<\tau$ is the unique restriction, used solely in the
proof of Proposition~\ref{pr:bindingaway}.
Note that by definition $\epsilon=\epsilon(k_0)$ and
that $k_0=k_0(\tau)$ according to \eqref{eq:condhaty}.
Thus $\tau$ can be made as small as needed.

Hence by choosing $1<\sigma<\sqrt{\tilde\sigma}<\sigma_0$
and a small $\delta$, we may then take $0<\rho<\tau$ as
small as we need to obtain a small $\ep>0$ (and $k_0$ big
enough, as a consequence, see Remark~\ref{rmk:k_0BIG}), and
then find $\flat>0$ in order that $\Theta=\Theta(\flat)$ be
big enough so that the constants $C_1,C_2,\xi_1,\xi_2$, and
consequently $C_3,\xi_3$ on the statements of
Section~\ref{sec:introduction}, are defined depending only
on $\alpha,\beta$, which depend only on $\hat f$. So $C_1,
C_2, C_3, \xi_1, \xi_2, \xi_3$ depend on
$\sigma,\sigma_0,\rho$ and $\tau$, but \emph{do not depend
  on} $\mu\in S$.
This concludes the proof of
Corollary~\ref{cor:contSRBEntropy}.


\end{document}